\UseRawInputEncoding
\documentclass[11pt,reqno]{amsart}
\usepackage{mathrsfs}
\usepackage{amsmath}
\usepackage{amssymb, color}
\usepackage{amsthm}
\usepackage{epsfig}
\usepackage{graphicx}
\usepackage{titletoc}
\usepackage{fullpage}
\usepackage{palatino,mathpazo}
\numberwithin{equation}{section}

\newtheorem{theorem}{Theorem}[section]
\newtheorem{proposition}[theorem]{Proposition}
\newtheorem{remark}[theorem]{Remark}
\newtheorem{lemma}[theorem]{Lemma}
\newtheorem{corollary}[theorem]{Corollary}

\newcommand{\R}{\mathbb{R}}

\newcommand{\dd}{\,\mathrm{d}}

\theoremstyle{definition}

\theoremstyle{remark}

\def\XXint#1#2#3{{\setbox0=\hbox{$#1{#2#3}{\int}$}
        \vcenter{\hbox{$#2#3$}}\kern-.5\wd0}}

\def\XXint#1#2#3{{\setbox0=\hbox{$#1{#2#3}{\int}$ }
        \vcenter{\hbox{$#2#3$ }}\kern-.6\wd0}}

\renewcommand{\t }{\tau }

\begin{document}
\title[Bifurcation analysis of Stokes waves with piecewise smooth vorticity in deep water]
{Bifurcation analysis of Stokes waves with piecewise smooth
vorticity in deep water}

\author[C.-F. Gui]{Changfeng Gui}
\address{\noindent Changfeng ~Gui,~Department of Mathematics, Faculty of Science and Technology, University of Macau, Taipa, Macau,
China.} \email{changfenggui@um.edu.mo}

\author[J. Wang]{Jun Wang}
\address{\noindent Jun ~Wang,~School of Mathematical Sciences, Jiangsu University,
    Zhenjiang, Jiangsu, 212013, P.R. China.}
\email{wangj2011@ujs.edu.cn}

\author[W. Yang]{Wen Yang}
\address{\noindent Wen ~Yang,~Department of Mathematics, Faculty of Science and Technology,
University of Macau, Taipa, Macau, P.R. China}
\email{wenyang@um.edu.mo}

\author[Y. Zhang]{Yong Zhang}
\address{\noindent Yong ~Zhang,~School of Mathematical Sciences, Jiangsu University, Zhenjiang, 212013, P.R. China.}
\email{zhangyong@ujs.edu.cn}

\thanks{Corresponding author: Yong Zhang, E-mail address: zhangyong@ujs.edu.cn. On behalf of all authors, the corresponding author states that there is no conflict of interest. \\
The research of C.F. Gui is supported by University of Macau research grants CPG2024-00016-FST, CPG202500032-FST, SRG2023-00011-FST, MYRG-GRG2023-00139-FST-UMDF, UMDF Professorial Fellowship of Mathematics, Macao SAR FDCT0003/2023/RIA1 and Macao SAR FDCT0024/2023/RIB1, NSFC No. 12531010. The research of
J. Wang is partially supported by National Key R$\&$D Program of
China 2022YFA1005601 and National Natural Science Foundation of
China 12371114. The research of W. Yang is partially supported by National
Key R\&D Program of China 2022YFA1006800, NSFC No. 12531010, No.
12171456 and No. 12271369, FDCT No. 0070/2024/RIA1, UMDF No. TISF/2025/006/FST, No. MYRG-GRG2025-00051-FST, No. MYRG-GRG2024-00082-FST and Startup
Research Grant No. SRG2023-00067-FST. The research of
Y. Zhang is partially supported by National Natural Science
Foundation of China No. 12301133 and the Postdoctoral Science
Foundation of China (No. 2023M741441, No. 2024T170353) and Jiangsu
Education Department (No. 23KJB110007)}

\begin{abstract}
In this paper, we establish the existence of Stokes waves with
piecewise smooth vorticity in a two-dimensional, infinitely deep
fluid domain. These waves represent traveling water waves
propagating over sheared currents in a semi-infinite cylinder, where
the vorticity may have jump discontinuities across internal
streamlines. The analysis is carried
out by applying a hodograph transformation, which reformulates the
original free boundary problem into an abstract elliptic boundary
value problem. Compared to previously studied steady water waves,
the present setting introduces several novel features: the presence
of an internal interface, an unbounded spatial domain, and a
non-Fredholm linearized operator. To address these difficulties, we
introduce a height function formulation, casting the problem as a
transmission problem with suitable transmission conditions. A
singular bifurcation approach is then employed, combining global
bifurcation theory with Whyburn's topological lemma.  The singular
limit is taken in a natural deep-water topology based on the surface
trace and the derivatives of the height function.  An exact mean-flux
identity controls the closing zero Fourier mode uniformly as the
regularization is removed.  This yields a connected
global continuum of exact waves together with explicit norm,
ellipticity, surface-obliqueness, and parameter-boundary alternatives.


\end{abstract}

\maketitle

{\scriptsize \tableofcontents }

\section{Introduction}

This work presents a rigorous construction of two-dimensional
periodic steady water waves in deep water with piecewise smooth
vorticity, propagating under gravity. Unlike most existing
mathematical treatments, we incorporate deep-water settings and
the vorticity jump produces an internal shear interface and the associated
transmission structure, without introducing density stratification. Indeed, these two topics
are important and promising research directions in water wave theory
as mentioned by Constantin in his book \cite[subsection
3.6]{Constantin}:

\emph{There are interesting possible extensions of the presented
theory outside the realm of gravity water waves propagating over a
flat bed:}

\emph{$\bullet$ Allowing discontinuous vorticity (the typical
example being a sudden jump in the vorticity) is technically
challenging since in this setting one has to investigate weak
solutions to nonlinear elliptic partial differential equations with
nonlinear boundary conditions;}

\emph{$\bullet$  The theory of rotational deep water waves (infinite
depth, with the velocity field and the vorticity decaying deep down
to capture the realistic assumption that the wave motion dies out)
is in its early stages.}

To the best of our knowledge, previous studies have addressed the
first aspect exclusively for finite-depth flat beds in
\cite{ConstantinS1,MartinM}. The second aspect was primarily
investigated by Hur \cite{Hur06,Hur11}, but under the restrictive
assumption of smooth vorticity. In the present work, we make a
comprehensive treatment of both aspects by establishing a global
continuum of periodic gravity waves in deep water with a piecewise
smooth vorticity distribution. The main difficulty
is that we are working in a domain where the top boundary is unknown
and the bottom is unbounded. In addition, the discontinuous
vorticity distribution would introduce a new interface inside fluid.
These difficulties are overcome by adopting the height function
formulation, employing a singular bifurcation argument, and applying
the Whyburn lemma (refer to subsection 1.2 for details).

\begin{figure}
\centering
\IfFileExists{1.jpg}{\includegraphics[width=9cm]{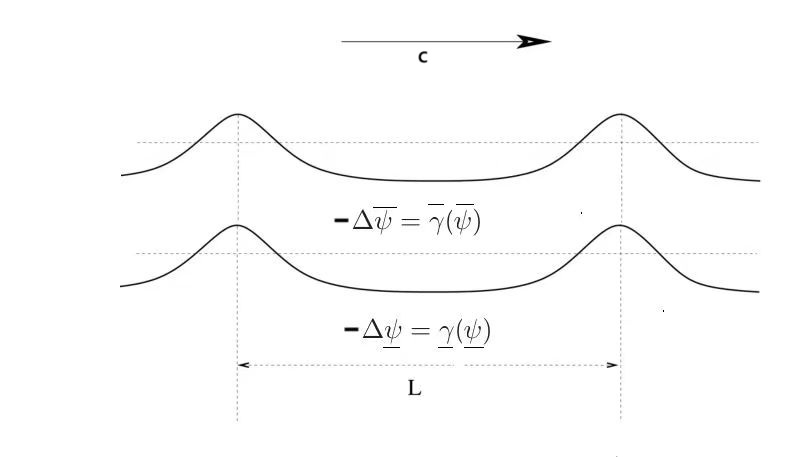}}{\fbox{\parbox[c][4cm][c]{8cm}{\centering Schematic of the fluid domain}}}
\caption{The schematic of the problem.} \label{fig1}
\end{figure}\par

\subsection{The historical background}

We first recall briefly the background of the problem. During the
twentieth century, most rigorous work concerned irrotational flows,
for which the velocity is the gradient of a harmonic potential and
complex-analytic tools play a central role. At the beginning of the
present century, the seminal work \cite{ConstantinS} established the
existence of water waves with arbitrary smooth vorticity
distributions by means of the hodograph transformation. However,
this approach requires the gradient of the stream function to be
non-vanishing, thereby precluding stagnation points and critical
layers in the resulting waves. Moreover, the free surface must be
the graph of a function, which excludes overturning profiles.
Based on methods developed in \cite{ConstantinS}, more and more rigorous
existence results for smooth rotational flows have been established
mathematically, which models other complex scenarios such as
incorporating fluid stratification \cite{ChuDE,ChuW,Walsh09}, the capillary
effects of surface tension \cite{CaoQ,Wahlen1} or their
combination \cite{Walsh14a}, the presence of interface
\cite{ChuE,ChuY,Matioc} and even accommodating the unbounded bottom
\cite{Hur06,Hur11}. The extension of the existence result in
\cite{ConstantinS1} to allow for discontinuous vorticity in bounded
domain is used to model more general steady flows, which leads to
considering the problem in weak sense. Subsequently, the existence
of capillary-gravity waves of small-amplitude propagating at
constant speed over a flat bed with a discontinuous vorticity was
considered in \cite{MartinM}, where authors use the height function
formulation associating to a transmission problem due to a jump of
the vorticity. In fact, the transmission problems are common in some
mathematical physical models, such as multiphase flows
\cite{Escher,Sinambela} and the Muskat problems
\cite{Cordoba,Cheng,MatiocM}. However, to the best of our knowledge,
there are fewer results limited to unbounded domains, which is the
gap we will address in this paper.

In addition, several other transform techniques are applicable to this class of free boundary problems. For instance, Wahl\'{e}n \cite{Wahlen2} utilized a flattening technique to stabilize the free surface, successfully constructing small-amplitude rotational waves with constant vorticity that allow for critical layers. This approach was recently extended by Varholm \cite{Varholm} to establish large-amplitude rotational waves with arbitrary vorticity distributions. However, this formulation encounters difficulties in proving suitable nodal properties in the presence of internal stagnation points. Moreover, it precludes the existence of overhanging wave profiles. We also would like to mention the recent work of Dai et al. \cite{DaiXZ}, who achieved secondary bifurcation for electrohydrodynamic waves with vorticity while allowing stagnation points by flattening transformation.

On the other hand, Constantin et al.
\cite{ConstantinSV,ConstantinV} used conformal mappings to represent
the fluid domain as the image of a strip and to construct both small-
and large-amplitude waves with constant vorticity.  This formulation
does not impose an a priori graph condition on the physical domain or
the stream function and therefore permits critical layers,
stagnation points, and overhanging profiles.  Related conformal
formulations have been extended to stratified waves \cite{Haziot},
capillary-gravity waves \cite{WangXZ}, and electrohydrodynamic waves
\cite{DaiFZ}.

Most recently, Wahlen and Weber \cite{Erik,Erik1} have also employed a conformal change of variables to establish the existence of large-amplitude capillary-gravity or gravity water waves, accommodating stagnation points, critical layers, and overhanging profiles. Their key advancement lies in removing assumptions on the vorticity distribution--unlike the constant vorticity restriction inherent in earlier conformal mapping approaches \cite{ConstantinSV,ConstantinV}.

\subsection{The plan of the paper}

Now we will outline the structure of the paper and explain the main
mathematical difficulties and how we approach them.

In Section 2, we introduce several reformulations of the problem as
in \cite{ConstantinS}. More precisely, we use the hodograph
transformation to map the fluid domain onto a semi-infinite periodic
strip. In the new coordinates, the steady Euler system becomes a
quasilinear elliptic PDE with oblique top boundary conditions. Then
we state our main result of this paper.

In Section 3, we rewrite the height formulation for piecewise smooth
vorticity as the transmission problem \eqref{eq3.2}, imposing
suitable interface conditions as in \cite{MartinM} for the jump
vorticity function. Then we consider the laminar flows (i.e. the
trivial solutions of (\ref{eq3.2})), the transmission problem
(\ref{eq3.2}) would reduce to ODEs (\ref{eq3.3}), whose solutions
are given as (\ref{eq3.4}) and (\ref{eq3.5}), where we introduce a
parameter $\lambda$ related to waves speed $c$.

The treatment of Section 4 is influenced by
\cite{ConstantinS,ConstantinSV,Hur11,MartinM}.  For the present
infinite-depth problem we follow the singular-bifurcation strategy of
\cite{Hur11}, but for each fixed regularization parameter
$\varepsilon>0$ we continue the local branch by the real-analytic
global bifurcation theorem of Dancer--Buffoni--Toland in the corrected
form stated in \cite[Theorem 6 and Remarks 7--8]{ConstantinSV}.
The jump of the vorticity produces a genuine transmission interface,
which is handled by the height formulation and the transmission
conditions introduced in Section 3.  Since the unbounded lower
cylinder destroys the Fredholm property of the exact linearization,
we introduce a family of Fredholm regularizations, as in
\cite{Hur11}.  For every fixed $\varepsilon>0$, the regularized map is
real analytic, its linearization is Fredholm of index zero, and the
compactness hypothesis required by analytic continuation follows from
the properness statement proved below.  Thus no generalized Fredholm
degree, orientation, or additional spectral crossing-number theory is
needed.  The only global alternative for the distinguished analytic
route is the corrected periodic-versus-escape alternative; the
periodic possibility will be excluded locally, when needed, by the
oriented nodal structure of the two local half-routes.

The heart of Section 5 is the singular limit.  The zero Fourier mode
of the regularized problem has decay length $O(\varepsilon^{-1/2})$,
so its absolute height cannot be required to decay uniformly as
$\varepsilon\downarrow0$.  We instead derive an exact mean-flux identity
and a threshold-uniform half-line estimate controlling the derivative
of the mean and the product $\varepsilon$ times the mean.  This leads
to a natural deep-water topology consisting of the upper-layer trace
and the lower-layer height derivatives.  In this topology the singular
solution set is compact on bounded admissible families and Whyburn's
lemma yields an exact connected continuum.

Section 6 translates the resulting global alternative into the physical
variables.  We state explicitly the possible norm, ellipticity,
surface-obliqueness and parameter-boundary mechanisms; the lower
parameter and surface-obliqueness mechanisms are shown to force large
vertical stretching.

\section{Equivalent formulations and main results}

In this section, we introduce several reformulations of the problem
that will make it convenient to state our main results.

\subsection{Governing equations in velocity field formulation}

Let's recall the governing equations for two-dimensional steady
Stokes waves in infinite depth. These are periodic waves over a
rotational, inviscid and incompressible fluid. Choose Cartesian
coordinates $(X,Y)$ such that $X$-axis points to the horizontal and
$Y$-axis points to the vertical. Assume that the free surface is
given by $Y=\eta(t,X)$, $(u(t,X,Y), v(t,X,Y))$ is the velocity field
of the flow and $P=P(t,X,Y)$ is the pressure. All of these functions
depend on $(X-ct)$ and $Y$ in steady periodic travelling waves,
where $c$ represents the speed of waves. For convenience, let
$x=X-ct, y=Y$ and consider the problem in
$\Omega_\eta=\{(x,y)|-\infty<y<\eta(x)\}$.

The incompressibility gives that the vector field $(u,v)$ is
divergence free
\begin{eqnarray}
u_{x}+v_{y}=0. \label{eq2.1}
\end{eqnarray}
Combining conservation of momentum with the boundary conditions, the
velocity formulation is the following nonlinear problem:
\begin{eqnarray}
\left\{\begin{array}{llll}
{(u-c) u_{x}+  v u_{y}=-P_{x}} & {\text { in } \Omega_{\eta}}, \\
{(u-c) v_{x}+ v v_{y}=-P_{y}-g} & { \text { in } \Omega_{\eta}}, \\
{P=P_{atm}} & {\text { on } y=\eta(x)}, \\
{v=(u-c) \eta_{x}} & {\text { on } y=\eta(x)}, \\
{(u,v)\rightarrow (0,0)} & {\text { as } y\rightarrow
-\infty},\end{array}\right. \label{eq2.2}
\end{eqnarray}
where $P_{atm}$ is the constant atmosphere pressure, and $g$ is the
gravitational acceleration at the Earth's surface. We suppose that
the flow is free from stagnation points, that is
\begin{eqnarray}
 u<c \label{eq2.3}
\end{eqnarray}
throughout the fluid, which implies that flows are unidirectional.

\subsection{Governing equation in stream function formulation}
 To reformulate the problem (\ref{eq2.2}) into a simpler one, we may introduce a pseudo-stream function $\psi=\psi(x,y)$ satisfying
\begin{eqnarray}
\psi_{x}=-v, ~~~~~~~\psi_{y}=u-c. \label{eq2.4}
\end{eqnarray}
The level sets of $\psi$ can be regarded as streamlines of the flow,
thus we normalize $\psi=0$ on the free boundary $y=\eta(x)$.
In addition, under the assumption
(\ref{eq2.3}), there exists a vorticity function $\gamma$ defined on
$[0, \infty)$ such that
\begin{eqnarray}
-\Delta \psi=\gamma(\psi).  \label{eq2.5}
\end{eqnarray}
From the Euler equation (\ref{eq2.2}) we obtain Bernoulli's law,
which states that
$$
E=\frac{1}{2}((u-c)^{2}+v^{2})+gy+P-\Gamma(-\psi)
$$
where
\begin{eqnarray}
 \Gamma(p)=\int^p_0 \gamma(-s)ds \label{eq2.6}
 \end{eqnarray}
is bounded for $p\in(-\infty,0]$, and $E$ is the hydraulic head,
which is constant along each streamline. Denote
$$
\Gamma_{\inf}:=\inf_{(-\infty,0]}\Gamma, \qquad
\Gamma_{\infty}:=\lim_{p\to-\infty}\Gamma(p).
$$
The limit exists by the decay assumption on $\gamma$, and
$-\Gamma_{\inf}\geq0$. Evaluating
Bernoulli's law on the free surface $y=\eta(x)$, we obtain
\begin{eqnarray}
|\nabla\psi|^{2}+2g\eta=0~~~~~~\text { on }~~y=\eta(x),
\label{eq2.7}
\end{eqnarray}
where, after a vertical translation of the $y$-coordinate, we
normalize the Bernoulli constant
$Q=2(E|_{\eta}-P_{atm})$ to be zero.

Summarizing these considerations gives
\begin{eqnarray}
\left\{\begin{array}{ll}{- \Delta\psi=\gamma(\psi)} & {\text { in } \Omega_{\eta}}, \\
{|\nabla\psi|^{2}+2gy=0} & { \text { on } y=\eta(x)}, \\
{\psi=0} & {\text { on } y=\eta(x)}, \\
{\nabla\psi\rightarrow (0,-c)} & {\text { as } y\rightarrow
-\infty}.\end{array}\right. \label{eq2.8}
\end{eqnarray}
The assumption (\ref{eq2.3}) would transform into
\begin{eqnarray}
\psi_{y}<0~~~~~~\text { in } \Omega_{\eta}. \label{eq2.9}
\end{eqnarray}
It is not difficult to find that (\ref{eq2.9}) forbids the presence
of stagnation points except the surface stagnation points.

\subsection{Governing equation in height function formulation}

The assumption (\ref{eq2.9}) enables us to introduce the
Dubreil-Jacotin's transformation by
\begin{eqnarray}
q=x,~~~~~p=-\psi(x,y), \label{eq2.10}
\end{eqnarray}
which transforms the fluid domain
\begin{eqnarray}
{\Omega_{\eta}=\{(x,y):~-\infty<y<\eta(x)\}} \nonumber
\end{eqnarray}
into rectangular domain
\begin{eqnarray}
D=\{(q,p):~ -\infty< p< 0\}. \nonumber
\end{eqnarray}
The function $\gamma$ in (\ref{eq2.8}) can be written as
\begin{eqnarray}
\gamma=\gamma(-p). \label{eq2.11}
\end{eqnarray}
Define the height function by
\begin{eqnarray}
h(q,p):=y. \label{eq2.12}
\end{eqnarray}
A direct calculation gives
\begin{eqnarray}
\psi_{y}=-\frac{1}{h_{p}},~~~~ \psi_{x}=\frac{h_{q}}{h_{p}},
\label{eq2.13}
\end{eqnarray}
\begin{eqnarray}
\partial_{q}=\partial_{x}+h_{q}\partial_{y},~~~~ \partial_{p}=h_{p}\partial_{y}. \label{eq2.14}
\end{eqnarray}
It follows from (\ref{eq2.9}) and a simple computation that
\begin{eqnarray}
h_{q}=\frac{v}{u-c},~~~~ h_{p}=\frac{1}{c-u}>0. \label{eq2.15}
\end{eqnarray}
Thus we can rewrite the governing equations in terms of the height
function $h$ by
\begin{eqnarray}
\left\{\begin{array}{lll}{\left(1+h_{q}^{2}\right) h_{p p}-2 h_{q} h_{p} h_{q p}+h_{p}^{2} h_{q q}=-\gamma(-p)h_{p}^{3}}  & {\text { in } -\infty<p<0}, \\
{1+h_{q}^{2}+2ghh_{p}^{2}=0} & {\text { on } p=0}, \\
{\nabla h=(h_q, h_p)\rightarrow (0, \frac{1}{c})} & {\text { as }
p\rightarrow -\infty}.\end{array}\right. \label{eq2.16}
\end{eqnarray}
To construct Stokes waves, we seek height functions $h$ that are even
and $2\pi$-periodic in the $q$-variable.

\subsection{Main results}
In this paper, we construct solutions of problem
(\ref{eq2.15})-(\ref{eq2.16}) in the case when the vorticity
function is a piecewise smooth function. More precisely, we suppose
that there exists a finite number $p_0\in (-\infty, 0)$ such that
\begin{eqnarray} \label{eq2.17}
\gamma \in C^{1,\alpha}([0,-p_0))\cap C^{1,\alpha}([-p_0, \infty)),
\end{eqnarray}
which means that at $p=p_0$, the vorticity function has a jump. Our
main result is the following theorem:

\begin{theorem}\label{thm2.1}
Suppose that $\alpha\in(0,1)$ and
\[
 \gamma\in C^{1,\alpha}([0,-p_0))
 \cap C^{1,\alpha}([-p_0,\infty)),
 \qquad
 \gamma(s)=O(s^{-2-r})\quad(s\to\infty)
\]
for some $r>0$. Assume also the local-bifurcation condition
\begin{equation}\label{eq2.bifcondition}
\int_{-\infty}^{0}
\left\{
2\bigl(2\Gamma(p)-2\Gamma_{\inf}\bigr)^{3/2}
+
\bigl(2\Gamma(p)-2\Gamma_{\inf}\bigr)^{1/2}
\right\}e^{2p}\,\dd p<g.
\end{equation}
For the endpoint H\"older control used only in the cylindrical-end
Fredholm parametrix, assume in addition
\begin{equation}\label{eq2.asymptotic-holder}
 \lim_{S\to\infty}
 \|\gamma\|_{C^{1,\alpha}([S,S+2])}=0.
\end{equation}
Then there exists a connected set
$\mathcal{K}$ of solutions $(c,h)$ of
(\ref{eq2.15})--(\ref{eq2.16}).  If
$w=h-H(\cdot\,;\lambda)$ denotes the perturbation of the
corresponding laminar height, then, for every point of
$\mathcal K$, $w$ belongs piecewise to $C^{3,\alpha}$.  The set $\mathcal K$ is connected in the
strong physical topology $\mathbb R\times X_{\rm s}$ defined in
\eqref{eq5dw:Xstrong}.  Moreover,

\begin{itemize}
  \item [(1)] $\mathcal{K}$ contains a trivial solution that corresponds to a
laminar flow solution $H(p)$ given as (\ref{eq3.4}) and
(\ref{eq3.5});
  \item [(2)] there exists a sequence of solutions $(c_k, h_k)\subset\mathcal{K}$,
for which either
\begin{equation*}
(a)~~ c_k\rightarrow \infty;\quad
\text{or}\quad (b)~~ \lim_{k\rightarrow
\infty}\sup_{\overline{D}}\partial_p h_k\rightarrow \infty.
\end{equation*}
\end{itemize}
\end{theorem}

We make a few remarks about this result.

\begin{remark}

\begin{itemize}
  \item [(1)] Our primary contribution is to construct a global
connected continuum of Stokes waves for a broad class of piecewise
smooth vorticities in deep water.  Along this continuum either the
wave speed is unbounded or the flow approaches horizontal stagnation
in the precise sense stated in part~(2) of Theorem~\ref{thm2.1}.
The proof uses a singular global-bifurcation argument and Whyburn's
lemma; no separate assertion that the free-surface amplitude itself
is unbounded is made.

\item [(2)] In Theorem \ref{thm2.1}, we only assume a single discontinuity
of vorticity function. In fact, we can extend our theory to finitely
many discontinuities by supposing
$$
\gamma\in C^{1,\alpha}([0,-p_0))\cap C^{1,\alpha}([-p_0,
-p_1))\cap\cdot\cdot\cdot\cap C^{1,\alpha}([-p_{n-1},-p_n))\cap
C^{1,\alpha}([-p_n, \infty))
$$
for $-\infty<p_n<p_{n-1}<\cdot\cdot\cdot<p_1<p_0<0$.

\item [(3)] Alternative (a) means that the wave speed is unbounded.
If alternative (b) occurs while the wave speeds remain bounded, then,
using
\[
 h_{k,p}=\frac{1}{c_k-u_k}
\]
in the hodograph variables, we obtain
\[
 \inf_D(c_k-u_k)\longrightarrow0.
\]
Equivalently, there exist points $(q_k,p_k)\in D$ such that
\[
 u_k(q_k,p_k)-c_k\longrightarrow0.
\]
Thus the sequence comes arbitrarily close to horizontal stagnation
somewhere in the fluid; no pointwise convergence of $u_k(q,p)$ at a
fixed point is asserted.  This is consistent with the limiting
behavior of Stokes waves with smooth vorticity and bounded depth
\cite{Constantin,ConstantinS}.
\end{itemize}
\end{remark}

\section{The related transmission problem and laminar flow}
In this section, we decompose the height function $h$ and the vorticity
function $\gamma$ into two piecewise functions as follows:
\begin{eqnarray}
h=\left\{\begin{array}{ll}{\overline{h}(q,p)},  & {\text { for } p_0\leq p\leq0}, \\
{\underline{h}(q,p)}, & {\text { for } -\infty<p\leq
p_0},\end{array}\right.
\gamma=\left\{\begin{array}{ll}{\overline{\gamma}(-p)},  & {\text { for } p_0<p\leq0}, \\
{\underline{\gamma}(-p)}, & {\text { for } -\infty<p\leq
p_0}.\end{array}\right. \label{eq3.1}
\end{eqnarray}
We associate with \eqref{eq2.16} and the piecewise vorticity
\eqref{eq2.17} the following transmission problem:
\begin{eqnarray}
\left\{\begin{array}{lll}
{\left(1+\overline{h}_{q}^{2}\right) \overline{h}_{p p}-2 \overline{h}_{q} \overline{h}_{p} \overline{h}_{q p}+\overline{h}_{p}^{2} \overline{h}_{q q}+\overline{\gamma}(-p)\overline{h}_{p}^{3}}=0  & {\text { in } p_0<p<0}, \\
{1+\overline{h}_{q}^{2}+2\overline{h}g\overline{h}_{p}^{2}=0} & {\text { on } p=0}, \\
{\left(1+\underline{h}_{q}^{2}\right) \underline{h}_{p p}-2 \underline{h}_{q} \underline{h}_{p} \underline{h}_{q p}+\underline{h}_{p}^{2} \underline{h}_{q q}+\underline{\gamma}(-p)\underline{h}_{p}^{3}}=0  & {\text { in } -\infty<p<p_0}, \\
{\overline{h}=\underline{h}} & {\text { on } p=p_0}, \\
{\overline{h}_p=\underline{h}_p} & {\text { on } p=p_0}, \\
{\nabla \underline{h}=(\underline{h}_q, \underline{h}_p)\rightarrow
(0, \frac{1}{c})} & {\text { as } p\rightarrow
-\infty}.\end{array}\right. \label{eq3.2}
\end{eqnarray}
It follows from \cite{MartinM} that if
$(\overline{h},\underline{h})\in
C^{3,\alpha}(\mathbb{R}\times[p_0,0])\times
C^{3,\alpha}(\mathbb{R}\times(-\infty,p_0])$ is a solution of
(\ref{eq3.2}), then the function $h: \overline{D}\rightarrow
\mathbb{R}$ defined by (\ref{eq3.1}) belongs to
$C^{1,\alpha}(\overline{D})\cap
C^{3,\alpha}(\mathbb{R}\times(-\infty,p_0])\cap
C^{3,\alpha}(\mathbb{R}\times[p_0,0])$ and solves (\ref{eq2.16})
with piecewise vorticity function (\ref{eq2.17}). The solution $h$
solves the last two boundary conditions of (\ref{eq2.16}) in
classical sense and solves the first main equation of (\ref{eq2.16})
also in classical sense for $p\in(-\infty,p_0)\cup(p_0,0)$, but
solves the first main equation of (\ref{eq2.16}) almost everywhere
in $D$ in the following weak sense
$$
\int_{D}\left\{\frac{h_q}{h_p}\psi_q+
\left(\Gamma-\frac{1+h_q^2}{2h_p^2}\right)\psi_p\right\}dqdp=0,
~~~\text{for~all}~~~\psi\in C_0^1(D)
$$

Before investigating the nontrivial solutions of operator equation
(\ref{eq3.2}), we first consider the trivial solutions, that is the
laminar flow solutions of (\ref{eq3.2}). These solutions describe
water waves with a flat surface and parallel streamlines. To this
end, we introduce an additional parameter $\lambda$ into the problem
(\ref{eq3.2}). In the following, we denote the laminar flow
solutions by $(\overline{H},\underline{H})$, which depends only on
the variable $p$. That is to say, $(\overline{H},\underline{H})$
solves the system
\begin{eqnarray}
\left\{\begin{array}{llll}
{\overline{H}_{pp}+\gamma(-p)\overline{H}_p^3=0},  & {\text { in } p_0<p<0}, \\
{1+2g\overline{H}(0)\overline{H}_{p}^{2}(0)=0},  \\
{\underline{H}_{pp}+\gamma(-p)\underline{H}_p^3=0},  & {\text { in } -\infty<p<p_0}, \\
{\overline{H}(p_0)=\underline{H}(p_0)}, \\
{\overline{H}_p(p_0)=\underline{H}_p(p_0)}, \\
{\underline{H}_p\rightarrow \frac{1}{c}}, & {\text { as }
p\rightarrow -\infty}.\end{array}\right. \label{eq3.3}
\end{eqnarray}
By observation, there exists $\lambda>-2\Gamma_{\inf}\geq 0$ such
that it holds that
\begin{eqnarray}
\overline{H}(p,\lambda)=\int_{0}^{p}\frac{1}{\sqrt{\lambda+2\Gamma(s)}}ds-\frac{\lambda}{2g},~~~~
p\in[p_0,0] \label{eq3.4}
\end{eqnarray}
and
\begin{eqnarray}
\underline{H}(p,\lambda)=\int_{0}^{p}\frac{1}{\sqrt{\lambda+2\Gamma(s)}}ds-\frac{\lambda}{2g},~~~~
p\in(-\infty,p_0]. \label{eq3.5}
\end{eqnarray}
We would like to mention that the parameter $\lambda$ in
(\ref{eq3.4}) and (\ref{eq3.5}) is the same one due to the boundary
condition on $p=p_0$. In particular, the parameter $\lambda$ can be
explicitly expressed by $\lambda=\overline{H}_p^{-2}(0)$. It follows
from the last boundary condition that the speed of wave propagation
is determined by $\lambda$, that is
$$
c^2=\lambda+2\Gamma_{\infty}.
$$

\section{The bifurcation of the approximating problem}
\subsection{The global bifurcation of the approximating problem}

In this subsection, we will first introduce an appropriate function
space to recast the problem (\ref{eq3.2}) into an abstract operator
equation. Since we seek periodic solutions, it is sufficient to
consider a domain of one wavelength. Let
$$
R_1:=\{(q,p):-\pi\leq q\leq \pi,~~ p_0<p<0~~ \text{with}~~ q=\pm
\pi~~ \text{identified}\}
$$
and
$$
R_2:=\{(q,p):-\pi\leq q\leq \pi,~~ -\infty<p<p_0~~ \text{with}~~
q=\pm \pi~~ \text{identified}\}.
$$
For $k\in\mathbb N_0$, let $C^{k,\alpha}_{\rm per,0}(\overline R_2)$
denote the closed subspace of periodic uniform H\"older functions for
which
\begin{equation}\label{eq:vanishing-holder-end}
 \lim_{P\to-\infty}\sup_{s<P}
 \|f\|_{C^{k,\alpha}(\mathbb S\times[s-1,s+1])}=0.
\end{equation}
Thus the subscript $0$ records convergence to zero in the full local
H\"older norm, rather than only pointwise convergence of the
derivatives.  This endpoint convention is needed in the Fredholm
parametrix below.
Define
$$
X:=\{(\overline{f},\underline{f})\in
C^{3,\alpha}_{per}(\overline{R_1})\times
C^{3,\alpha}_{\rm per,0}(\overline{R_2}): \overline{f}=\underline{f},
\overline{f}_p=\underline{f}_p~\text{on}~p=p_0\}
$$
and
$$
Y_1:=\{(\overline{f},\underline{f})\in
C^{1,\alpha}_{per}(\overline{R_1})\times
C^{1,\alpha}_{\rm per,0}(\overline{R_2})\}
$$
and
$Y_2:=C^{2,\alpha}_{per}(\mathbb{S})$, where the subscript "per"
means evenness and $2\pi$-periodicity in the $q$ variable.
 Here $\mathbb{S}$ denotes the $2\pi$-circle
on $p=0$. In order to tackle the existence of solutions for
(\ref{eq3.2}) by bifurcation theory, we let
\begin{eqnarray}
h=\left\{\begin{array}{ll}{\overline{h}(q,p)=\overline{H}(p)+\overline{w}(q,p)},  & {\text { for } p_0<p\leq0}, \\
{\underline{h}(q,p)=\underline{H}(p)+\underline{w}(q,p)}, & {\text {
for } -\infty<p\leq p_0},\end{array}\right. \label{eq4.1}
\end{eqnarray}
and introduce the operator $F: (-2\Gamma_{\inf}, \infty)\times
X\rightarrow Y:=Y_1\times Y_2$ with
\begin{eqnarray}
F(\lambda,\overline{w},\underline{w})=(F_1(\lambda,\overline{w}),F_2(\lambda,\underline{w}),F_3(\lambda,\overline{w}))=0
\label{eq4.2}
\end{eqnarray}
by the following formulations
\begin{eqnarray}
F_1(\lambda,\overline{w})&=&\left(1+\overline{w}_q^2\right)\overline{w}_{pp}-2(a^{-1}(\lambda)+\overline{w}_p)\overline{w}_q\overline{w}_{qp}+\left(a^{-1}(\lambda)+\overline{w}_p\right)^2\overline{w}_{qq}\nonumber\\
&~&+\gamma(-p)\left(a^{-1}(\lambda)+\overline{w}_p\right)^3-\gamma(-p)a^{-3}(\lambda)(1+\overline{w}_q^2),\label{eq4.3}
\end{eqnarray}
\begin{eqnarray}
F_2(\lambda,\underline{w})&=&\left(1+\underline{w}_q^2\right)\underline{w}_{pp}-2(a^{-1}(\lambda)+\underline{w}_p)\underline{w}_q\underline{w}_{qp}+\left(a^{-1}(\lambda)+\underline{w}_p\right)^2\underline{w}_{qq}\nonumber\\
&~&+\gamma(-p)\left(a^{-1}(\lambda)+\underline{w}_p\right)^3
-\gamma(-p)a^{-3}(\lambda)(1+\underline{w}_q^2),\label{eq4.4}
\end{eqnarray}
and
\begin{eqnarray}
F_3(\lambda,\overline{w})=\frac{\lambda}{2}
\left\{1+(2g\overline{w}-\lambda)
\left(\lambda^{-\frac{1}{2}}+\overline{w}_p\right)^2
+\overline{w}_q^2\right\}\Big|_{p=0}.
\label{eq4.5}
\end{eqnarray}
with $a(\lambda)=a(p;\lambda)=\sqrt{\lambda+2\Gamma(p)}$.
The factor $\lambda/2$ is a harmless normalization: $\lambda>0$
throughout the parameter domain, so it does not change the zero set.
It is included in the definition of the operator in order that the
third component of the linearization agree exactly, rather than only
up to a nonzero scalar multiple, with the boundary operator used in
the range and transversality calculations below.

To carry out the bifurcation analysis, denote by
$\partial_{(\overline{w},\underline{w})}F(\lambda,
\overline{w},\underline{w})$ the Fr\'{e}chet derivative of $F$ at
$(\lambda,\overline{w},\underline{w})\in \mathbb{R}\times X$.
Direct differentiation gives
$$
\partial_{(\overline{w},\underline{w})}F(\lambda, \overline{w},\underline{w})=(L_1(\lambda,\overline{w}), L_2(\lambda,\underline{w}), L_3(\lambda,\overline{w}))\in \mathbb{L}(X, Y),
$$
where
\begin{eqnarray}
L_1(\lambda,\overline{w})[u]&=&\left(1+\overline{w}_q^2\right)u_{pp}-2\left(a^{-1}(\lambda)+\overline{w}_p\right)\overline{w}_qu_{qp}+\left(a^{-1}(\lambda)+\overline{w}_p\right)^2u_{qq} \nonumber\\
&~&+\left[ -2\overline{w}_q\overline{w}_{qp}+2\left(a^{-1}(\lambda)+\overline{w}_p\right)\overline{w}_{qq}+3\gamma(-p)\left(a^{-1}(\lambda)+\overline{w}_p\right)^2 \right]u_p \nonumber\\
&~&+\left[
2\overline{w}_q\overline{w}_{pp}-2\left(a^{-1}(\lambda)+\overline{w}_p\right)\overline{w}_{qp}-2\gamma(-p)a^{-3}(\lambda)\overline{w}_q
\right]u_q, \label{eq4.6}
\end{eqnarray}
\begin{eqnarray}
L_2(\lambda,\underline{w})[v]&=&\left(1+\underline{w}_q^2\right)v_{pp}-2\left(a^{-1}(\lambda)+\underline{w}_p\right)\underline{w}_qv_{qp}+\left(a^{-1}(\lambda)+\underline{w}_p\right)^2v_{qq} \nonumber\\
&~&+\left[ -2\underline{w}_q\underline{w}_{qp}+2\left(a^{-1}(\lambda)+\underline{w}_p\right)\underline{w}_{qq}+3\gamma(-p)\left(a^{-1}(\lambda)+\underline{w}_p\right)^2 \right]v_p \nonumber\\
&~&+\left[
2\underline{w}_q\underline{w}_{pp}-2\left(a^{-1}(\lambda)+\underline{w}_p\right)\underline{w}_{qp}-2\gamma(-p)a^{-3}(\lambda)\underline{w}_q
\right]v_q \label{eq4.7}
\end{eqnarray}
and
\begin{eqnarray}
L_3(\lambda,\overline{w})[u]
&=&\lambda g\left(\lambda^{-\frac{1}{2}}+\overline{w}_p\right)^2u
+\lambda(2g\overline{w}-\lambda)
\left(\lambda^{-\frac{1}{2}}+\overline{w}_p\right)u_p
+\lambda\overline{w}_qu_q\Big|_{p=0}.
\label{eq4.8}
\end{eqnarray}
for $(u,v)\in X$. In the infinite cylinder, the linearized operator
$\partial_{(\overline{w},\underline{w})}F(\lambda,
\overline{w},\underline{w})$ of the exact problem is not Fredholm.
The obstruction is the zero Fourier mode of the limiting operator at
$p=-\infty$, exactly as in the smooth-vorticity deep-water problem
\cite[Section 2.4]{Hur11}; in particular, the limiting operator is
not invertible and the range of the exact linearization is not
closed.  We therefore follow \cite{Hur11} and study the Fredholm
regularizations
\begin{eqnarray}
 F^{\varepsilon}(\lambda,\overline{w},\underline{w})=0, \label{eq4.9}
\end{eqnarray}
where
\begin{eqnarray}
F^{\varepsilon}(\lambda,\overline{w},\underline{w})=(F_1(\lambda,\overline{w})-\varepsilon
a^{-3}(\lambda)\overline{w},F_2(\lambda,\underline{w})-\varepsilon
a^{-3}(\lambda) \underline{w},F_3(\lambda,\overline{w})).
\label{eq4.10}
\end{eqnarray}
For the approximate problems we use the corrected analytic global
bifurcation theorem for real-analytic Fredholm maps; see
\cite[Theorem 6 and Remarks 7--8]{ConstantinSV}, based on the theory
of Dancer and Buffoni--Toland \cite{Dancer,BuffoniT}.  We record the
local real-analytic reparametrization in (C4) and the corrected
bounded-closed global alternative.  The loop argument below uses
(C4) exactly as it is used in the nodal analysis of
\cite[Section 4]{ConstantinSV}: an incoming distinguished arc cannot
be continued by retracing the same local analytic half-arc.  We do
not use discreteness of the singular parameters, a nonvanishing
derivative of the bifurcation parameter, or an integer-period
criterion for repeated points.  No additional degree or spectral
crossing theory is involved.

\begin{theorem}[Analytic global bifurcation]\label{thm4.1}
Let $X$ and $Y$ be Banach spaces, let $\mathcal O\subset\mathbb R\times X$
be open, and let $F:\mathcal O\to Y$ be real analytic. Assume that:
\begin{itemize}
\item[(H1)] $F(\lambda,0)=0$ whenever $(\lambda,0)\in\mathcal O$;
\item[(H2)] for some $(\lambda_*,0)\in\mathcal O$,
\[
 \dim\ker \partial_\psi F(\lambda_*,0)
 =\operatorname{codim}\operatorname{Ran}\partial_\psi F(\lambda_*,0)=1,
\]
with
\[
 \ker \partial_\psi F(\lambda_*,0)=\operatorname{span}\{\psi_*\},
\]
and the transversality condition
\[
 \partial_{\lambda\psi}F(\lambda_*,0)\psi_*
 \notin\operatorname{Ran}\partial_\psi F(\lambda_*,0)
\]
holds;
\item[(H3)] $\partial_\psi F(\lambda,\psi):X\to Y$ is Fredholm of
index zero at every $(\lambda,\psi)\in\mathcal O$ satisfying
$F(\lambda,\psi)=0$;
\item[(H4)] for every bounded closed set $Q\subset\mathcal O$,
\[
 Q\cap\{(\lambda,\psi)\in\mathcal O:F(\lambda,\psi)=0\}
\]
is compact.
\end{itemize}
Then there exists a continuous global route
\[
 \mathscr R=\{z(s)=(\lambda(s),\psi(s)):s\in\mathbb R\}
 \subset\mathcal O
\]
of solutions of $F(\lambda,\psi)=0$ such that:
\begin{itemize}
\item[(C1)] $z(0)=(\lambda_*,0)$;
\item[(C2)] $\psi(s)=s\psi_*+o(s)$ in $X$ as $s\to0$;
\item[(C3)] in a neighborhood of $(\lambda_*,0)$ the nontrivial
solutions are precisely $z(s)$ with $0<|s|<s_0$;
\item[(C4)] $\mathscr R$ has a real-analytic reparametrization
locally around each of its points;
\item[(C5)] one of the following alternatives occurs:
\begin{itemize}
\item[(a)] for every bounded closed set $Q\subset\mathcal O$,
there is $S_Q>0$ such that
\[
 z(s)\notin Q\qquad\text{whenever }|s|>S_Q;
\]
\item[(b)] there exists $T>0$ such that
\[
 z(s+T)=z(s)\qquad\text{for all }s\in\mathbb R.
\]
\end{itemize}
\end{itemize}
Moreover, the distinguished route with properties (C1)--(C5) is
unique up to reparametrization.
\end{theorem}

\begin{remark}\label{rem:analytic-global-alternative}
The compactness hypothesis in (H4) is slightly stronger than the
minimal exhaustion hypothesis in \cite[Theorem 6]{ConstantinSV}, but
it is exactly what Lemma \ref{lem4.5} provides in our application.
Accordingly, the corrected conclusion in (C5)(a) may be stated for
every bounded closed $Q\subset\mathcal O$; compare
\cite[Remark 8]{ConstantinSV}.  The conclusion in (C5)(a) must not be replaced
by the stronger and generally unjustified assertion that either
$\|z(s)\|\to\infty$ or
$\operatorname{dist}(z(s),\partial\mathcal O)\to0$.  Also, the
analytic route is a distinguished continuation of the local curve;
it need not coincide with the maximal connected component of the
entire solution set.  In (C4), ``reparametrization locally around a
point'' has the distinguished-arc meaning of
\cite[Theorem 6(C4) and Remark 7]{ConstantinSV}: the incoming and
outgoing distinguished arcs at that point possess one local
real-analytic uniformizing parameter.  Consequently those two sides
cannot be the same analytic half-arc traversed in opposite
directions.  This is the local analytic no-retracing consequence used
below, and it is precisely the use of (C4) made in
\cite[Section 4]{ConstantinSV}.  No assertion that the singular
parameter set is discrete, and no condition that
$\lambda'(s)\not\equiv0$, is being added to Theorem
\ref{thm4.1}.
\end{remark}

In order to use Theorem \ref{thm4.1}, we work on the following
open subset of $\mathbb R\times X$:
\[
\mathcal O_\delta=\left\{
(\lambda,\overline w,\underline w)
\in(-2\Gamma_{\inf},\infty)\times X:
\begin{aligned}
&\lambda>-2\Gamma_{\inf}+\delta,\\
&\inf_{\overline R_1}\bigl(a^{-1}(p;\lambda)+\overline w_p\bigr)>\delta,\\
&\inf_{\overline R_2}\bigl(a^{-1}(p;\lambda)+\underline w_p\bigr)>\delta,\\
&\sup_{q\in\mathbb S}\left(\overline w(q,0)
-\frac{2\lambda-\delta}{4g}\right)<0
\end{aligned}
\right\}.
\]
The use of the infima is essential on the unbounded lower cylinder:
a pointwise strict inequality need not define an open condition in the
global H\"older topology.
For any $\delta>0$, $\varepsilon\geq0$, and
$(\lambda,\overline{w},\underline{w})\in \mathcal{O}_{\delta}$, the
operator
\[
 D_wF^{\varepsilon}(\lambda,w)
 =\bigl(L_1-\varepsilon a^{-3}I,
 L_2-\varepsilon a^{-3}I,L_3\bigr):X\longrightarrow Y
\]
is continuous, where $I$ denotes the identity map. Moreover,
H\"older spaces are Banach algebras and inversion is analytic on
sets bounded away from zero. Hence $F^{\varepsilon}:
\mathcal{O}_{\delta}\rightarrow Y$ is real analytic. In the following, our goal is to
construct for each $\varepsilon$ a global connected set of
nontrivial solutions to (\ref{eq4.10}) by Theorem \ref{thm4.1}. We
first show the Fredholm property.

\begin{lemma}({\bf Fredholm property})\label{lem4.2}
Let $\alpha\in(0,1)$ and suppose
\[
 \gamma\in C^{1,\alpha}([0,-p_0))
 \cap C^{1,\alpha}([-p_0,\infty)),
 \qquad \gamma(s)=O(s^{-2-r})\quad(s\to\infty)
\]
for some $r>0$. Assume also \eqref{eq2.asymptotic-holder}.  For
each $\delta>0$, $\varepsilon>0$, and
$(\lambda,\overline{w},\underline{w})\in \mathcal{O}_{\delta}$, the
linear operator
\[
 D_wF^{\varepsilon}(\lambda,w)
 =\bigl(L_1(\lambda,\overline w)-\varepsilon a^{-3}I,
 L_2(\lambda,\underline w)-\varepsilon a^{-3}I,
 L_3(\lambda,\overline w)\bigr):X\longrightarrow Y
\]
is Fredholm of index zero.
\end{lemma}
\begin{proof}
Set
\[
 \mathcal L^\varepsilon
 :=\partial_{(\overline w,\underline w)}F^\varepsilon
 (\lambda,\overline w,\underline w).
\]
The inequalities defining $\mathcal O_\delta$ imply ellipticity in
both layers and obliqueness at $p=0$.  Moreover,
$\overline w=\underline w$ on $p=p_0$ implies equality of all
tangential derivatives there, and the additional condition
$\overline w_p=\underline w_p$ implies that the principal
coefficients have the same traces from the two sides.  Hence
\[
 u^+=u^-,\qquad u_p^+=u_p^-
 \quad\text{on }p=p_0
\]
are perfect-transmission conditions satisfying the complementing
condition.

We give the cylindrical-end construction in detail because the
endpoint H\"older scale requires one important distinction.  For a
coefficient $c$ on the lower end, write
\[
 \|c\|_{1,\alpha;P}:=
 \sup_{s<P}\|c\|_{C^{1,\alpha}
 (\mathbb S\times[s-1,s+1])}.
\]
The definition of $X$, the decay of $\gamma$, and
\eqref{eq2.asymptotic-holder} imply that every coefficient of
$\mathcal L^\varepsilon$ converges, in this unit-strip norm, to the
corresponding coefficient of
\begin{equation}\label{eq:limit-operator-eps}
 \mathcal L_\infty^\varepsilon
 =\partial_{pp}+a_\infty^{-2}\partial_{qq}
  -\varepsilon a_\infty^{-3},
 \qquad a_\infty=(\lambda+2\Gamma_\infty)^{1/2}.
\end{equation}
Thus, if $\vartheta(P)$ denotes the operator-norm size of the
coefficient difference on $p<P$, then
\begin{equation}\label{eq:tail-coefficient-small}
 \vartheta(P)\longrightarrow0\qquad(P\to-\infty).
\end{equation}
Here and below the norms are the piecewise uniform H\"older norms
defining $X$ and $Y$.

The full-cylinder operator in \eqref{eq:limit-operator-eps} is an
isomorphism from periodic $C^{3,\alpha}$ to periodic
$C^{1,\alpha}$.  Indeed, its $k$th Fourier coefficient is
\[
 \partial_{pp}-\rho_k^2,
 \qquad
 \rho_k^2=k^2a_\infty^{-2}+\varepsilon a_\infty^{-3},
\]
and $\rho_k\geq\sqrt{\varepsilon a_\infty^{-3}}>0$.  The Green
kernel $e^{-\rho_k|p-s|}/(2\rho_k)$ and the standard periodic
Schauder estimate give a bounded inverse, denoted by
$G_\infty^\varepsilon$.

Choose $P<p_0-4$ and nested smooth cutoffs whose transition regions
are contained in $\mathbb S\times[P-3,P]$.  On the bounded
transmission cylinder above $p=P-4$, impose an auxiliary Dirichlet
condition at $p=P-4$.  The Agmon--Douglis--Nirenberg boundary
estimate and the transmission estimate of Ladyzhenskaya
\cite{Agmon,Lady} give a two-sided core parametrix.  Patch it with
$G_\infty^\varepsilon$ by the chosen cutoffs.  Extension across the
artificial cross-section is taken before applying
$G_\infty^\varepsilon$, and the outer cutoff removes the extension
region afterwards.  This produces bounded left and right
parametrices (which need not be the same), and the two compositions
have the form
\begin{equation}\label{eq:small-compact-parametrix}
 \begin{aligned}
  Q_{P,L}\mathcal L^\varepsilon
    &=I+S_{P,X}-K_{P,X},\\
  \mathcal L^\varepsilon Q_{P,R}
    &=I+S_{P,Y}-K_{P,Y}.
 \end{aligned}
\end{equation}
The terms $K_{P,X}$ and $K_{P,Y}$ are compact.  They consist of the
regularizing remainders of the bounded-core parametrices and
commutators with the cutoffs.  Every such commutator has order at
most one, is supported in a fixed bounded cylinder, and is compact
from $C^{3,\alpha}$ to $C^{1,\alpha}$ by a local H\"older embedding
with a positive loss of differentiability.

In contrast, the terms containing the second-order coefficient
difference are not declared compact.  They are precisely the
operators $S_{P,X}$ and $S_{P,Y}$.  These terms occur only in the
tail leg of the construction.  The boundedness of the translated
extension and cutoff operators and of the full-cylinder inverse gives
\begin{equation}\label{eq:small-parametrix-errors}
 \|S_{P,X}\|_{\mathbb L(X)}+
 \|S_{P,Y}\|_{\mathbb L(Y)}
 \leq C\vartheta(P).
\end{equation}
where $C$ depends on the fixed cutoff profile but not on its
translation level $P$.  Taking the transition sufficiently deep and using
\eqref{eq:tail-coefficient-small}, we arrange that both norms in
\eqref{eq:small-parametrix-errors} are smaller than $1/2$.  Therefore
$I+S_{P,X}$ and $I+S_{P,Y}$ are invertible by Neumann series.  From
\eqref{eq:small-compact-parametrix} we obtain
\[
 (I+S_{P,X})^{-1}Q_{P,L}\mathcal L^\varepsilon
  =I-(I+S_{P,X})^{-1}K_{P,X},
\]
and
\[
 \mathcal L^\varepsilon Q_{P,R}(I+S_{P,Y})^{-1}
  =I-K_{P,Y}(I+S_{P,Y})^{-1}.
\]
These are genuine left and right parametrices modulo compact
operators.  Atkinson's theorem now proves that
$\mathcal L^\varepsilon:X\to Y$ is Fredholm.  Notice that this
argument absorbs the tail principal-part error as a small bounded
operator; compact support is used only for the lower-order
commutators.  This is the transmission analogue of the
cylindrical-end construction in \cite[Lemma 3.1]{Hur11}.

It remains to compute the index.  Consider
\[
\mathcal L_{\rm ref}(u^+,u^-)
=
\left(
\Delta u^+-u^+,\,
\Delta u^--u^-,\,
-u_p^+-u^+\big|_{p=0}
\right)
\]
with the same transmission conditions.  For an arbitrary right-hand
side in $Y$, Fourier expansion reduces the problem to the uniquely
solvable equations
\[
 u_k''-(k^2+1)u_k=f_k,
 \qquad u_k(-\infty)=0,
 \qquad -u_k'(0)-u_k(0)=g_k.
\]
The Green-function representation, followed by the piecewise
Schauder estimate at $p=p_0$, gives a unique element of $X$; the
mode estimates are uniform because $\sqrt{k^2+1}\geq1$.  Hence the
inverse is bounded from $Y$ to $X$.  Thus
$\mathcal L_{\rm ref}:X\to Y$ is an isomorphism and
$\operatorname{ind}\mathcal L_{\rm ref}=0$.

Finally set
\[
 \mathcal L_t=(1-t)\mathcal L_{\rm ref}
 +t\mathcal L^\varepsilon,
 \qquad 0\leq t\leq1.
\]
The convex combination of the two principal matrices remains
uniformly positive definite.  The coefficient of $u_p$ in the top
boundary operator of $\mathcal L^\varepsilon$ is
\[
 \lambda(2g\overline w-\lambda)
   (\lambda^{-1/2}+\overline w_p)<0,
\]
while the corresponding reference coefficient is $-1$; hence
obliqueness persists along the whole homotopy.  The transmission
conditions are unchanged.  At the cylindrical end the limit of
$\mathcal L_t$ is
\[
 \partial_{pp}+c_t\partial_{qq}-d_t,
 \qquad
 c_t=(1-t)+ta_\infty^{-2}>0,
 \quad
 d_t=(1-t)+t\varepsilon a_\infty^{-3}>0.
\]
Since $t\in[0,1]$ and $\varepsilon>0$, the lower bounds for $c_t$ and
$d_t$ are uniform in $t$.  The cutoff construction
\eqref{eq:small-compact-parametrix}--
\eqref{eq:small-parametrix-errors} therefore applies uniformly to
$\mathcal L_t$ after one sufficiently deep choice of $P$.  Every
$\mathcal L_t$ is Fredholm, and $t\mapsto\mathcal L_t$ is continuous
in operator norm.  Homotopy invariance of the index gives
\[
 \operatorname{ind}\mathcal L^\varepsilon
 =\operatorname{ind}\mathcal L_{\rm ref}=0.
\]
\end{proof}

It is known that the linearization of $F^{\varepsilon}$ at the
trivial solution $(\lambda,0,0)$ is
$$\partial_{(\overline{w},\underline{w})}F^{\varepsilon}(\lambda, 0,
0)=(L_1(\lambda,0)-\varepsilon a^{-3}(\lambda)I,
L_2(\lambda,0)-\varepsilon a^{-3}(\lambda)I, L_3(\lambda,0)),$$
where
\begin{eqnarray}
(L_1(\lambda,0)-\varepsilon
a^{-3}(\lambda)I)[u]=u_{pp}+a^{-2}(\lambda)u_{qq}+3\gamma(-p)a^{-2}(\lambda)u_p-\varepsilon
a^{-3}(\lambda)u, \label{eq4.14}
\end{eqnarray}
\begin{eqnarray}
(L_2(\lambda,0)-\varepsilon
a^{-3}(\lambda)I)[v]=v_{pp}+a^{-2}(\lambda)v_{qq}+3\gamma(-p)a^{-2}(\lambda)v_p-\varepsilon
a^{-3}(\lambda) v \label{eq4.15}
\end{eqnarray}
and
\begin{eqnarray}
L_3(\lambda,0)[u]= gu-\lambda^{\frac{3}{2}}u_p\mid_{p=0}
\label{eq4.16}
\end{eqnarray}
for $(u,v)\in X$. A necessary condition for bifurcation at a trivial
solution $(\lambda,0,0)$ is that
$\partial_{(\overline{w},\underline{w})}$ $F^{\varepsilon}(\lambda,
0, 0)$ from $X$ to $Y$ is not injective, which means that the
following problem
\begin{eqnarray}
\left\{\begin{array}{llll}
{u_{pp}+a^{-2}(\lambda)u_{qq}+3\gamma(-p)a^{-2}(\lambda)u_p-\varepsilon a^{-3}(\lambda)u=0} & {\text { in } R_1}, \\
{v_{pp}+a^{-2}(\lambda)v_{qq}+3\gamma(-p)a^{-2}(\lambda)v_p-\varepsilon a^{-3}(\lambda) v=0} & { \text { in } R_2}, \\
{gu(q,0)-\lambda^{\frac{3}{2}}u_p(q,0)=0}\end{array}\right.
\label{eq4.17}
\end{eqnarray}
admits a nontrivial solution in $X$. Now we give the kernel space of
operator
$\partial_{(\overline{w},\underline{w})}F^{\varepsilon}(\lambda, 0,
0)$ as follows.

For $\varepsilon\geq 0$, define the quadratic form
\begin{equation}
\label{eq:quadratic-form}
\mathfrak q_{\varepsilon,\lambda}[\Phi]
:=
-g\Phi(0)^2
+
\int_{-\infty}^{0}
\left(
a^3(p;\lambda)\Phi_p^2+\varepsilon\Phi^2
\right)\dd p
\end{equation}
and the weighted norm
\begin{equation}
\label{eq:weighted-norm}
\|\Phi\|_{a,\lambda}^2
:=
\int_{-\infty}^{0}a(p;\lambda)\Phi^2\,\dd p.
\end{equation}
The associated principal Rayleigh value is
\begin{equation}
\label{eq:principal-Rayleigh-value}
\mu^\varepsilon(\lambda)
:=
\inf_{\substack{\Phi\in H^1(-\infty,0)\\ \Phi\not\equiv 0}}
\frac{\mathfrak q_{\varepsilon,\lambda}[\Phi]}
{\|\Phi\|_{a,\lambda}^2}.
\end{equation}

\begin{lemma}[Kernel of the linearized approximate problem]
\label{lem4.3}
Assume that
\[
\gamma\in
C^{1,\alpha}\bigl([0,-p_0)\bigr)
\cap
C^{1,\alpha}\bigl([-p_0,\infty)\bigr),
\qquad
\alpha\in(0,1),
\]
that
\[
\gamma(s)=O(s^{-2-r})
\qquad\text{as }s\to\infty
\]
for some $r>0$, and that the local-bifurcation condition \eqref{eq2.bifcondition} holds.

Then there exists $\varepsilon_0>0$ such that, for every
$0<\varepsilon<\varepsilon_0$, there is a unique
\[
\lambda_*^\varepsilon\in(-2\Gamma_{\inf},\infty)
\]
for which
\begin{equation}
\label{eq:crossing}
\mu^\varepsilon(\lambda_*^\varepsilon)=-1.
\end{equation}
At $\lambda=\lambda_*^\varepsilon$, the kernel of the linearized operator is
one-dimensional:
\begin{equation}
\label{eq:kernel-formula}
\ker
\partial_{(w,\underline w)}
F^\varepsilon(\lambda_*^\varepsilon,0,0)
=
\operatorname{span}
\left\{
\left(
\phi^\varepsilon(p)\cos q,\,
\varphi^\varepsilon(p)\cos q
\right)
\right\}.
\end{equation}
Here
\[
\Psi^\varepsilon(p)
=
\begin{cases}
\phi^\varepsilon(p),&p\in[p_0,0],\\
\varphi^\varepsilon(p),&p\in(-\infty,p_0],
\end{cases}
\]
is a nontrivial eigenfunction satisfying
\begin{equation}
\label{eq:one-dimensional-eigenproblem}
\left\{
\begin{aligned}
-\bigl(a^3(p;\lambda_*^\varepsilon)\Psi_p^\varepsilon\bigr)_p
+\varepsilon\Psi^\varepsilon
&=
-a(p;\lambda_*^\varepsilon)\Psi^\varepsilon
&&\text{in }(-\infty,p_0)\cup(p_0,0),\\
(\lambda_*^\varepsilon)^{3/2}\Psi_p^\varepsilon(0)
&=g\Psi^\varepsilon(0),\\
[\Psi^\varepsilon]_{p=p_0}
=
[\Psi_p^\varepsilon]_{p=p_0}
&=0,\\
\Psi^\varepsilon(p),\ \Psi_p^\varepsilon(p)
&\longrightarrow0
&&\text{as }p\to-\infty.
\end{aligned}
\right.
\end{equation}
Moreover,
\begin{equation}
\label{eq:root-convergence}
\lambda_*^\varepsilon\longrightarrow\lambda_*^0
\qquad\text{as }\varepsilon\downarrow0,
\end{equation}
where $\lambda_*^0$ is the unique solution of
\[
\mu^0(\lambda)=-1
\]
in $(-2\Gamma_{\inf},\infty)$.
\end{lemma}

\begin{proof}
We divide the proof into five steps.

\medskip
\noindent
\textbf{Step 1: Variational realization of every negative Rayleigh value.}

Fix $\lambda>-2\Gamma_{\inf}$ and $\varepsilon\geq0$.  Then
\[
a(p;\lambda)\geq
a_{\min}(\lambda)
:=
\bigl(\lambda+2\Gamma_{\inf}\bigr)^{1/2}>0.
\]
Consequently, the denominator in
\eqref{eq:principal-Rayleigh-value} is equivalent to the
$L^2(-\infty,0)$ norm.  The trace inequality on the half-line gives,
for every $\eta>0$,
\[
|\Phi(0)|^2
\leq
\eta\|\Phi_p\|_{L^2(-\infty,0)}^2
+
C_\eta\|\Phi\|_{L^2(-\infty,0)}^2.
\]
It follows that the form
$\mathfrak q_{\varepsilon,\lambda}$ is closed and bounded from below on
$H^1(-\infty,0)$.

Suppose that $\mu^\varepsilon(\lambda)<0$.  Let $\{\Phi_n\}$ be a minimizing
sequence normalized by
\[
\|\Phi_n\|_{a,\lambda}=1.
\]
The trace estimate and the positive lower bound for $a$ show that
$\{\Phi_n\}$ is bounded in $H^1(-\infty,0)$.  After passing to a
subsequence,
\[
\Phi_n\rightharpoonup\Phi
\quad\text{weakly in }H^1(-\infty,0),
\]
and
\[
\Phi_n\longrightarrow\Phi
\quad\text{strongly in }L^2_{\mathrm{loc}}(-\infty,0]
\]
and at the trace point $p=0$.

The weak limit cannot vanish identically.  Indeed, if $\Phi\equiv0$,
then $\Phi_n(0)\to0$, whereas all integral terms in
\eqref{eq:quadratic-form} are nonnegative.  This would imply
\[
\liminf_{n\to\infty}
\mathfrak q_{\varepsilon,\lambda}[\Phi_n]\geq0,
\]
contrary to $\mu^\varepsilon(\lambda)<0$.

By weak lower semicontinuity,
\[
\mathfrak q_{\varepsilon,\lambda}[\Phi]
\leq
\liminf_{n\to\infty}
\mathfrak q_{\varepsilon,\lambda}[\Phi_n]
=
\mu^\varepsilon(\lambda),
\]
while
\[
0<\|\Phi\|_{a,\lambda}^2\leq1.
\]
Since the right-hand side is negative,
\[
\frac{\mathfrak q_{\varepsilon,\lambda}[\Phi]}
{\|\Phi\|_{a,\lambda}^2}
\leq
\mathfrak q_{\varepsilon,\lambda}[\Phi]
\leq
\mu^\varepsilon(\lambda).
\]
The opposite inequality follows from the definition of
$\mu^\varepsilon(\lambda)$.  Hence the infimum is attained.

The Euler--Lagrange equation is
\begin{equation}
\label{eq:euler-lagrange-general}
-\bigl(a^3\Psi_p\bigr)_p+\varepsilon\Psi
=
\mu^\varepsilon(\lambda)a\Psi
\quad\text{in }(-\infty,p_0)\cup(p_0,0),
\end{equation}
with
\[
a^3(0;\lambda)\Psi_p(0)=g\Psi(0).
\]
Since $\Gamma(0)=0$, we have $a^3(0;\lambda)=\lambda^{3/2}$.
The weak formulation also yields continuity of $\Psi$ and of
$a^3\Psi_p$ at $p=p_0$.  As $a$ is continuous and strictly positive,
this is equivalent to
\[
[\Psi]_{p=p_0}=[\Psi_p]_{p=p_0}=0.
\]
Piecewise elliptic regularity gives
\[
\Psi\in
C^{3,\alpha}\bigl((-\infty,p_0]\bigr)
\cap
C^{3,\alpha}\bigl([p_0,0]\bigr).
\]
Moreover, since the coefficients converge at $p=-\infty$ and the
spectral value is negative, the one-dimensional equation has an
exponential dichotomy.  The $H^1$ condition removes the growing
solution; consequently, for some $C,\rho>0$,
\[
 |\partial_p^j\Psi(p)|\leq Ce^{\rho p},
 \qquad p\leq p_0-1,\qquad 0\leq j\leq3.
\]
In particular the eigenfunction belongs to the function space used in
\eqref{eq:kernel-formula}.

\medskip
\noindent
\textbf{Step 2: The nonpositive spectral subspace has dimension at most one.}

For every nonzero $\Phi\in H^1(-\infty,0)$ satisfying $\Phi(0)=0$,
\[
\mathfrak q_{\varepsilon,\lambda}[\Phi]
=
\int_{-\infty}^{0}
\left(
a^3\Phi_p^2+\varepsilon\Phi^2
\right)\dd p>0
\]
for every $\varepsilon\geq0$.  Therefore, the maximal dimension of a subspace on
which $\mathfrak q_{\varepsilon,\lambda}$ is nonpositive is at most one.
Indeed, if a two-dimensional subspace $E$ were nonpositive, then the
linear trace map
\[
E\ni\Phi\longmapsto\Phi(0)\in\R
\]
would have a nontrivial kernel.  A nonzero element of that kernel
would have strictly positive quadratic form, a contradiction.

It follows in particular that, whenever
$\mu^\varepsilon(\lambda)<0$, this negative eigenvalue is simple and is the
only nonpositive eigenvalue of the generalized Sturm--Liouville
problem \eqref{eq:euler-lagrange-general}.

\medskip
\noindent
\textbf{Step 3: Existence and uniqueness of the crossing
$\mu^\varepsilon(\lambda)=-1$.}

Set
\[
\lambda_0:=-2\Gamma_{\inf},
\qquad
a_0(p):=
\bigl(2\Gamma(p)-2\Gamma_{\inf}\bigr)^{1/2}.
\]
Condition \eqref{eq2.bifcondition} implies, in particular, that
\[
\delta_0
:=
g-
\int_{-\infty}^{0}
\bigl(a_0^3(p)+a_0(p)\bigr)e^{2p}\,\dd p
>0.
\]
Choose
\[
0<\varepsilon_0<\delta_0.
\]
For $\lambda>\lambda_0$, evaluate the Rayleigh quotient at
$\Phi(p)=e^p$.  Since
\[
\int_{-\infty}^{0}e^{2p}\,\dd p=\frac12,
\]
we have
\begin{align}
&\mathfrak q_{\varepsilon,\lambda}[e^p]
+\|e^p\|_{a,\lambda}^2
\nonumber\\
&\qquad=
-g+
\int_{-\infty}^{0}
\left(
a^3(p;\lambda)+a(p;\lambda)
\right)e^{2p}\,\dd p
+\frac{\varepsilon}{2}.
\label{eq:left-trial-expression}
\end{align}
As $\lambda\downarrow\lambda_0$, the right-hand side of
\eqref{eq:left-trial-expression} converges to
\[
-\delta_0+\frac{\varepsilon}{2}.
\]
After decreasing $\varepsilon_0$ if necessary, there exists
$\lambda_->\lambda_0$, independent of
$0<\varepsilon<\varepsilon_0$, such that
\[
\mathfrak q_{\varepsilon,\lambda_-}[e^p]
+\|e^p\|_{a,\lambda_-}^2<0.
\]
Hence
\begin{equation}
\label{eq:left-sign}
\mu^\varepsilon(\lambda_-)<-1.
\end{equation}

Next set
\[
\lambda_+:=g-2\Gamma_{\inf}.
\]
For every $\Phi\in H^1(-\infty,0)$, we have
\[
a(p;\lambda_+)\geq\sqrt g.
\]
Moreover,
\[
\Phi(0)^2
=
2\int_{-\infty}^{0}\Phi\Phi_p\,\dd p
\leq
g^{-1/2}\int_{-\infty}^{0}\Phi^2\,\dd p
+
g^{1/2}\int_{-\infty}^{0}\Phi_p^2\,\dd p.
\]
Consequently,
\begin{align*}
\mathfrak q_{\varepsilon,\lambda_+}[\Phi]
+\|\Phi\|_{a,\lambda_+}^2
&=
-g\Phi(0)^2
+
\int_{-\infty}^{0}
\left(
a^3\Phi_p^2+(a+\varepsilon)\Phi^2
\right)\dd p\\
&\geq
\varepsilon\int_{-\infty}^{0}\Phi^2\,\dd p>0
\end{align*}
for every nonzero $\Phi$ when $\varepsilon>0$.  Since
$a(\cdot;\lambda_+)$ is bounded above, with
\[
 a_{\max,+}:=\sup_{p\leq0}a(p;\lambda_+)<\infty,
\]
we have
\[
 \int_{-\infty}^0\Phi^2\,\dd p
 \geq a_{\max,+}^{-1}\|\Phi\|_{a,\lambda_+}^2.
\]
Therefore,
\begin{equation}
\label{eq:right-sign}
\mu^\varepsilon(\lambda_+)
\geq-1+\frac{\varepsilon}{a_{\max,+}}>-1.
\end{equation}

On every compact interval contained in
$(-2\Gamma_{\inf},\infty)$, the coefficients of the quadratic forms
depend continuously on $\lambda$, and hence
$\lambda\mapsto\mu^\varepsilon(\lambda)$ is continuous.  The signs
\eqref{eq:left-sign}--\eqref{eq:right-sign} therefore yield at least
one
\[
\lambda_*^\varepsilon\in(\lambda_-,\lambda_+)
\]
such that \eqref{eq:crossing} holds.

We now prove uniqueness.  Whenever $\mu^\varepsilon(\lambda)<0$, Step~2
shows that it is a simple isolated eigenvalue.  We normalize its
eigenfunction by
\[
\int_{-\infty}^{0}a(p;\lambda)\Psi^2(p)\,\dd p=1.
\]
Standard perturbation theory for simple eigenvalues of closed
quadratic forms implies differentiability with respect to $\lambda$.
Since
\[
\partial_\lambda a=\frac1{2a},
\qquad
\partial_\lambda(a^3)=\frac32a,
\]
the Hellmann--Feynman formula gives
\begin{equation}
\label{eq:mu-derivative}
\frac{\dd}{\dd\lambda}\mu^\varepsilon(\lambda)
=
\frac32
\int_{-\infty}^{0}a\Psi_p^2\,\dd p
-
\frac{\mu^\varepsilon(\lambda)}2
\int_{-\infty}^{0}a^{-1}\Psi^2\,\dd p.
\end{equation}
Both terms on the right-hand side are strictly positive whenever
$\mu^\varepsilon(\lambda)<0$.  Thus
\[
\frac{\dd}{\dd\lambda}\mu^\varepsilon(\lambda)>0
\qquad\text{whenever }\mu^\varepsilon(\lambda)<0.
\]
In particular, the equation
$\mu^\varepsilon(\lambda)=-1$ has exactly one solution.  This
conclusion follows from the strict derivative formula
\eqref{eq:mu-derivative}; no unproved global shape of the eigenvalue
curve is being assumed.
\medskip
\noindent
\textbf{Step 4: Identification of the full kernel, including the
zero Fourier mode.}

Let $(u,\underline u)$ belong to the kernel of
$\partial_{(w,\underline w)}F^\varepsilon(\lambda_*^\varepsilon,0,0)$.  Since the
function space consists of even $2\pi$-periodic functions, we must
include the zero Fourier mode and write
\[
u(q,p)=\sum_{k=0}^{\infty}\phi_k(p)\cos(kq),
\qquad
\underline u(q,p)
=
\sum_{k=0}^{\infty}\varphi_k(p)\cos(kq).
\]
For each $k\geq0$, the piecewise function
\[
\Psi_k(p)
=
\begin{cases}
\phi_k(p),&p\in[p_0,0],\\
\varphi_k(p),&p\in(-\infty,p_0],
\end{cases}
\]
satisfies
\begin{equation}
\label{eq:k-mode}
-\bigl(a^3\Psi_{k,p}\bigr)_p+\varepsilon\Psi_k
=
-k^2a\Psi_k,
\qquad
\lambda^{3/2}\Psi_{k,p}(0)=g\Psi_k(0),
\end{equation}
together with the transmission and decay conditions.

At $\lambda=\lambda_*^\varepsilon$, the value $-1$ is the principal
eigenvalue.  By Step~2 it is simple and is the only nonpositive
eigenvalue.  The mode $k=0$ would require the eigenvalue $0$, while
every mode $k\geq2$ would require the eigenvalue
$-k^2<-1$.  Both possibilities are impossible.  Therefore only
$k=1$ occurs, and its one-dimensional eigenspace yields
\eqref{eq:kernel-formula}.

\medskip
\noindent
\textbf{Step 5: Convergence of the entire family
$\lambda_*^\varepsilon$ as $\varepsilon\downarrow0$.}

Let
\[
I:=[\lambda_-,\lambda_+]
\Subset(-2\Gamma_{\inf},\infty)
\]
and set
\[
a_*:=\inf_{\substack{\lambda\in I\\p\leq0}}
a(p;\lambda)
=
\bigl(\lambda_-+2\Gamma_{\inf}\bigr)^{1/2}>0.
\]
For every nonzero $\Phi\in H^1(-\infty,0)$ and every $\lambda\in I$,
\[
0
\leq
\frac{\mathfrak q_{\varepsilon,\lambda}[\Phi]}
{\|\Phi\|_{a,\lambda}^2}
-
\frac{\mathfrak q_{0,\lambda}[\Phi]}
{\|\Phi\|_{a,\lambda}^2}
=
\varepsilon
\frac{\displaystyle\int_{-\infty}^{0}\Phi^2\,\dd p}
{\displaystyle\int_{-\infty}^{0}a\Phi^2\,\dd p}
\leq
\frac{\varepsilon}{a_*}.
\]
Taking infima gives
\begin{equation}
\label{eq:uniform-mu-convergence}
0
\leq
\mu^\varepsilon(\lambda)-\mu^0(\lambda)
\leq
\frac{\varepsilon}{a_*}
\qquad
\text{for every }\lambda\in I.
\end{equation}
Thus $\mu^\varepsilon\to\mu^0$ uniformly on $I$.

The argument above for $\varepsilon=0$ still shows that every negative
Rayleigh value is attained, simple, and strictly increasing in
$\lambda$.  Moreover,
\[
\mu^0(\lambda_-)<-1,
\qquad
\mu^0(\lambda_+)\geq-1.
\]
Hence there exists a unique
$\lambda_*^0\in I$ satisfying
\[
\mu^0(\lambda_*^0)=-1.
\]

Let $\varepsilon_j\downarrow0$ and suppose that, along a subsequence,
$\lambda_*^{\varepsilon_j}\to\bar\lambda\in I$.  By
\eqref{eq:uniform-mu-convergence},
\[
-1
=
\mu^{\varepsilon_j}(\lambda_*^{\varepsilon_j})
\longrightarrow
\mu^0(\bar\lambda).
\]
Therefore $\mu^0(\bar\lambda)=-1$, and uniqueness gives
$\bar\lambda=\lambda_*^0$.  Since every convergent subsequence has
the same limit, the whole family converges:
\[
\lambda_*^\varepsilon\longrightarrow\lambda_*^0.
\]
This proves \eqref{eq:root-convergence} and completes the proof.
\end{proof}

\begin{remark}
The proof deliberately expands the kernel element as
\[
\sum_{k=0}^{\infty}\Psi_k(p)\cos(kq),
\]
rather than starting with $k=1$.  Evenness and periodicity do not, by
themselves, impose zero horizontal mean.  The exclusion of the
$k=0$ mode follows instead from the fact that the quadratic form has
at most one nonpositive spectral direction and that this direction
is already occupied by the eigenvalue $-1$.
\end{remark}

For later use, fix once and for all a number $\delta_*>0$ so small
that
\[
 2\delta_*<\min\left\{
 \lambda_*^0+2\Gamma_{\inf},\;
 2\lambda_*^0,\;
 \inf_{p\leq0}a(p;\lambda_*^0)^{-1}
 \right\}.
\]
Then $(\lambda_*^0,0,0)\in\mathcal O_{2\delta_*}$.  By
\eqref{eq:root-convergence}, after decreasing $\varepsilon_0$ if
necessary,
\[
 (\lambda_*^\varepsilon,0,0)\in\mathcal O_{\delta_*}
 \qquad(0<\varepsilon<\varepsilon_0).
\]
In the remainder of Sections 4--6, $\delta$ is always chosen with
$0<\delta<\delta_*$.

\begin{lemma}({\bf Transversality condition})\label{lem4.4}
Let $\lambda_*^{\varepsilon}$ and
$(\phi^{\varepsilon},\varphi^{\varepsilon})$ be as in
Lemma~\ref{lem4.3}.  With
\[
 \Psi^*(q,p)
 :=(\phi^{\varepsilon}(p)\cos q,
    \varphi^{\varepsilon}(p)\cos q),
\]
the transversality condition is
\[
 \partial_{\lambda w}F^{\varepsilon}
 (\lambda_*^{\varepsilon},0)\Psi^*
 \notin
 \operatorname{Ran}\partial_wF^{\varepsilon}
 (\lambda_*^{\varepsilon},0).
\]
\end{lemma}
\begin{proof}
We first identify a cokernel functional.  If
$((f_1,f_2),f_3)\in
\operatorname{Ran}\partial_wF^{\varepsilon}
(\lambda_*^{\varepsilon},0)$, then
\begin{eqnarray} \label{eq4.24}
\int_{R_1}a^3(\lambda)u^*(q,p)f_1
dqdp+\int_{R_2}a^3(\lambda)v^*(q,p)f_2
dqdp+\int_{\mathbb{S}\times\{0\}}u^*(q,p) f_3 dq=0,
\end{eqnarray}
where
$(u^*(q,p),v^*(q,p))=(\phi^{\varepsilon}(p)\cos(q),\varphi^{\varepsilon}(p)\cos(q))$.

Indeed, since  $((f_1,f_2),f_3)\in \text{Im}~
\partial_{(\overline{w},\underline{w})}F^{\varepsilon}(\lambda_*^{\varepsilon},
0, 0)$, then there exist a pair $(u,v)\in X$ such that
\begin{eqnarray}
\left\{\begin{array}{llll}
{u_{pp}+a^{-2}(\lambda)u_{qq}+3\gamma(-p)a^{-2}(\lambda)u_p-\varepsilon a^{-3}(\lambda)u=f_1} & {\text { in } R_1}, \\
{v_{pp}+a^{-2}(\lambda)v_{qq}+3\gamma(-p)a^{-2}(\lambda)v_p-\varepsilon a^{-3}(\lambda) v=f_2} & { \text { in } R_2}, \\
{gu(q,0)-\lambda^{\frac{3}{2}}u_p(q,0)=f_3}.\end{array}\right.
\label{eq4.25}
\end{eqnarray}
On the other hand, $(u^*(q,p),v^*(q,p))$ satisfy
\begin{eqnarray}
\left\{\begin{array}{llll}
{u^*_{pp}+a^{-2}(\lambda)u^*_{qq}+3\gamma(-p)a^{-2}(\lambda)u^*_p-\varepsilon a^{-3}(\lambda)u^*=0} & {\text { in } R_1}, \\
{v^*_{pp}+a^{-2}(\lambda)v^*_{qq}+3\gamma(-p)a^{-2}(\lambda)v^*_p-\varepsilon a^{-3}(\lambda)v^*=0} & { \text { in } R_2}, \\
{gu^*(q,0)-\lambda^{\frac{3}{2}}u^*_p(q,0)=0}.\end{array}\right.
\label{eq4.26}
\end{eqnarray}
Based on these facts, we use integration to find that
\begin{eqnarray}
&~& \int_{R_1}a^3(\lambda)u^*(q,p)f_1 dqdp+\int_{R_2}a^3(\lambda)v^*(q,p)f_2 dqdp+\int_{\mathbb{S}\times\{0\}}u^*(q,p) f_3 dq\nonumber\\
&~&= -\int_{R_1}a^3u_pu^*_p+au_qu^*_q +\varepsilon uu^* dqdp-\int_{R_2}a^3v_pv^*_p+av_qv^*_q+\varepsilon vv^* dqdp+g\int_{\mathbb{S}\times\{0\}}uu^* dq\nonumber\\
&~&=0, \nonumber
\end{eqnarray}
where the first equality being obtained by using $(\ref{eq4.25})$
with a good observation
$a^3f_1=\left(a^3u_p\right)_p+\left(au_q\right)_q-\varepsilon u$ and
$a^3f_2=\left(a^3v_p\right)_p+\left(av_q\right)_q-\varepsilon v$ and
the last equality follows from $(\ref{eq4.26})$. This proves the
claim.  The normalization in \eqref{eq4.5} causes no hidden term in
the mixed derivative.  Indeed,
\[
 \partial_\lambda\{(\lambda/2)D_w\mathcal B\}[u^*,v^*]
 =\tfrac12D_w\mathcal B[u^*,v^*]
  +(\lambda/2)\partial_\lambda D_w\mathcal B[u^*,v^*],
\]
and the first term vanishes by the homogeneous boundary condition.
Here $\mathcal B$ denotes the unscaled expression in braces in
\eqref{eq4.5}.

A direct differentiation now gives
\begin{eqnarray}
\partial_{\lambda(\overline{w},\underline{w})}F^{\varepsilon}(\lambda_*^{\varepsilon}, 0, 0)(u^*,v^*)=\left\{\begin{array}{llll}
{-a^{-4}(\lambda_*^{\varepsilon})u^*_{qq}-3a_p(\lambda_*^{\varepsilon})a^{-3}(\lambda_*^{\varepsilon})u^*_p+\frac{3}{2}\varepsilon a^{-5}(\lambda_*^{\varepsilon})u^*:=g_1} \\
{-a^{-4}(\lambda_*^{\varepsilon})v^*_{qq}-3a_p(\lambda_*^{\varepsilon})a^{-3}(\lambda_*^{\varepsilon})v^*_p+\frac{3}{2}\varepsilon a^{-5}(\lambda_*^{\varepsilon})v^*:=g_2} \\
{-\frac{3}{2}\sqrt{\lambda_*^{\varepsilon}}u^*_p(q,0):=g_3}.\end{array}\right.
\nonumber
\end{eqnarray}
At last, we just need to verify that $(g_1,g_2,g_3)$ does not
satisfy (\ref{eq4.24}). In fact, by a simple computation, we can
deduce that
\begin{eqnarray}
&~&\int_{R_1}a^3(\lambda_*^{\varepsilon})u^*(q,p)g_1 dqdp+\int_{R_2}a^3(\lambda_*^{\varepsilon})v^*(q,p)g_2 dqdp+\int_{\mathbb{S}\times\{0\}}u^*(q,p) g_3 dq\nonumber\\
&~&= -\int_{R_1}a^{-1}u^*_{qq}u^*+3a_pu^*_pu^* -\frac{3}{2}\varepsilon a^{-2}u^{*2} dqdp-\frac{3\sqrt{\lambda_*^{\varepsilon}}}{2}\int_{\mathbb{S}\times\{0\}}u^*_pu^* dq\nonumber\\
&~&-\int_{R_1}a^{-1}v^*_{qq}v^*+3a_pv^*_pv^* -\frac{3}{2}\varepsilon
a^{-2}v^{*2} dqdp.  \label{eq4.27}
\end{eqnarray}
It follows from (\ref{eq4.26}) again that
\begin{eqnarray}
&~&\int_{R_1}a_pu^*_pu^* dqdp+\int_{R_2}a_pv^*_pv^* dqdp=-\frac{\sqrt{\lambda_*^{\varepsilon}}}{2}\int_{\mathbb{S}\times\{0\}}u^*_pu^* dq\nonumber\\
&~&+\int_{R_1}\frac{\varepsilon}{2}a^{-2}u^{*2}+\frac{a}{2}(u^*_p)^2
+\frac{1}{2a}(u^*_q)^2 dqdp+
\int_{R_2}\frac{\varepsilon}{2}a^{-2}v^{*2}+\frac{a}{2}(v^*_p)^2
+\frac{1}{2a}(v^*_q)^2 dqdp . \nonumber\\ \label{eq4.28}
\end{eqnarray}
With (\ref{eq4.27}) and (\ref{eq4.28}) in hand, we can obtain that
\begin{eqnarray}
&~&\int_{R_1}a^3(\lambda_*^{\varepsilon})u^*(q,p)g_1 dqdp+\int_{R_2}a^3(\lambda_*^{\varepsilon})v^*(q,p)g_2 dqdp+\int_{\mathbb{S}\times\{0\}}u^*(q,p) g_3 dq\nonumber\\
&~&=-\int_{R_1}\frac{3a}{2}(u^*_p)^2 +\frac{1}{2a}(u^*_q)^2
dqdp-\int_{R_2}\frac{3a}{2}(v^*_p)^2 +\frac{1}{2a}(v^*_q)^2
dqdp<0,\nonumber
\end{eqnarray}
which proves the lemma.
\end{proof}

We next isolate the global estimate used in the properness argument.
This prevents the local compactness and the cylindrical-end estimate
from being conflated.

\begin{lemma}[Uniform global transmission estimate]
\label{lem:global-schauder-fixed-eps}
Let $\mathfrak A$ be a family of linear operators having the same
bulk, transmission, and top-boundary structure as the averaged
linearizations of $F^\varepsilon$ over a bounded subset of
$\mathcal O_\delta$, where $\varepsilon>0$ is fixed.  Assume that the
bulk coefficients are uniformly bounded in the piecewise uniform
$C^{1,\alpha}$ norm, the top coefficients are uniformly bounded in
$C^{2,\alpha}$, and the ellipticity, transmission-complementing, and
top-obliqueness constants have a common positive lower bound.  Then
there is a constant $C$, uniform for $\mathscr A\in\mathfrak A$, such
that
\begin{equation}\label{eq:global-transmission-schauder}
 \|U\|_X\le C\left(
  \|\mathscr A U\|_Y+\|U\|_{C^0(\overline R_1\cup\overline R_2)}
 \right),\qquad U\in X.
\end{equation}
\end{lemma}
\begin{proof}
On the bounded cylinder
$\overline R_1\cup(\mathbb S\times[p_0-2,p_0])$, the boundary and
transmission Schauder estimates \cite{Agmon,Lady} give
\eqref{eq:global-transmission-schauder}, restricted to this core,
with the right-hand side enlarged to a fixed neighboring strip.  On
the lower end, cover by the overlapping unit strips
\[
 S_j=\mathbb S\times[p_0-j-2,p_0-j],\qquad j\ge1.
\]
The periodic interior Schauder estimate on a fixed enlargement
$S_j^*$ gives
\[
 \|u^-\|_{C^{3,\alpha}(S_j)}
 \le C\left(
  \|\mathscr A_2u^-\|_{C^{1,\alpha}(S_j^*)}
  +\|u^-\|_{C^0(S_j^*)}
 \right).
\]
Translation of the strip does not change $C$, and the stated
coefficient bounds make it independent of $j$ and of
$\mathscr A$.  Taking the supremum over $j$ controls all local
H\"older seminorms on the end.  Pairs of points separated by distance
at least one are controlled by twice the corresponding supremum
norm, so these local estimates are equivalent to the global uniform
H\"older norm.  Combining the core and tail estimates proves
\eqref{eq:global-transmission-schauder}.
\end{proof}

We now record the properness needed for global continuation.
\begin{lemma}({\bf Proper property})\label{lem4.5}
Assume the hypotheses of Lemma \ref{lem4.2}.  For every
$\varepsilon>0$ and every bounded closed set
$Q\subset\mathcal O_\delta$, the set
\[
 Q\cap(F^\varepsilon)^{-1}(K)
\]
is compact in $\mathbb R\times X$ for each compact $K\subset Y$.
In particular, $F^\varepsilon$ is proper on every member of a bounded
closed exhaustion of $\mathcal O_\delta$.
\end{lemma}
\begin{proof}
Let
$z_k=(\lambda_k,\overline w_k,\underline w_k)\in Q$ and assume
$F^\varepsilon(z_k)=y_k\to y$ in $Y$.  Passing to a subsequence,
$\lambda_k\to\lambda$ and local transmission Schauder compactness
gives
\[
 (\overline w_k,\underline w_k)
 \longrightarrow(\overline w,\underline w)
 \quad\text{in }C^{3,\beta}_{\rm loc}(\overline R_1)
 \times C^{3,\beta}_{\rm loc}(\overline R_2),
 \qquad 0<\beta<\alpha.
\]
The limit satisfies
$F^\varepsilon(\lambda,\overline w,\underline w)=y$ locally.  We
first record the cylindrical-end Liouville step that is needed below.

\smallskip
\noindent\emph{Claim.}
If $p_k\to-\infty$, then, after passage to a subsequence,
\begin{equation}\label{eq:shifted-background-zero}
 \underline w_k(q,p+p_k)\longrightarrow0,
 \qquad
 \underline w(q,p+p_k)\longrightarrow0
 \quad\text{in }C^{2,\beta}_{\rm loc}(\mathbb S\times\mathbb R).
\end{equation}
Indeed, a locally convergent subsequence of the first family has a
bounded limit $W$ on $\mathbb S\times\mathbb R$.  Since
$\gamma(-p-p_k)\to0$, $a(p+p_k;\lambda_k)\to a_\infty$, and the
translated right-hand sides tend locally to zero, $W$ solves
\begin{equation}\label{eq4.30}
 (1+W_q^2)W_{pp}
 -2(a_\infty^{-1}+W_p)W_qW_{pq}
 +(a_\infty^{-1}+W_p)^2W_{qq}
 -\varepsilon a_\infty^{-3}W=0.
\end{equation}
The inequalities defining $Q\subset\mathcal O_\delta$ make this
equation uniformly elliptic.  If $\sup W>0$, translate once more in
the $p$-variable along a maximizing sequence.  Local compactness
then produces a bounded solution of \eqref{eq4.30} attaining a
positive interior maximum.  At that point $W_p=W_q=0$ and the
Hessian is nonpositive.  Hence the second-order part of
\eqref{eq4.30} is nonpositive, whereas
$-\varepsilon a_\infty^{-3}W<0$, a contradiction.  Thus
$\sup W\leq0$.  If instead $\inf W<0$, translate along a minimizing
sequence and pass once more to a locally convergent subsequence.
The resulting bounded solution attains a negative interior minimum.
There $W_p=W_q=0$ and the Hessian is nonnegative, so the second-order
part of \eqref{eq4.30} is nonnegative, while
$-\varepsilon a_\infty^{-3}W>0$.  The left-hand side of
\eqref{eq4.30} is therefore strictly positive, again a contradiction.
Consequently $\inf W\geq0$, and hence $W\equiv0$.  Notice that this
argument does not use any invariance of \eqref{eq4.30} under
$W\mapsto-W$.  The argument for the second family is identical.
This proves the claim.

We now prove global $C^0$ convergence.  Suppose, to the contrary,
that there are $(q_k,p_k)\in R_2$, $p_k\to-\infty$, and $\kappa>0$
such that
\[
 |\underline w_k(q_k,p_k)-\underline w(q_k,p_k)|\geq\kappa.
\]
Set
\[
 \nu_k(q,p)=\underline w_k(q,p+p_k)
             -\underline w(q,p+p_k),
 \qquad p<p_0-p_k.
\]
By the Banach-space mean-value formula, the difference of the two
bulk equations is
\[
 \mathscr A_k\nu_k=r_k
 \quad\text{in }\mathbb S\times(-\infty,p_0-p_k),
\]
where $\mathscr A_k$ is the average of the linearizations along the
segment joining the shifted functions.  More explicitly,
\[
 \mathscr A_k
 =\int_0^1D_wF_2^\varepsilon
 \bigl(\lambda_k,
 \underline w(\cdot,\cdot+p_k)+t\nu_k\bigr)\,\dd t,
\]
and $r_k$ is the sum of the translated difference $y_{2,k}-y_2$
and the parameter difference obtained by replacing $\lambda_k$ by
$\lambda$.  Thus $r_k\to0$ in $C^{0,\beta}_{\rm loc}$.  By the
claim, all first derivatives of both shifted background functions
tend locally to zero.  Consequently
\[
 \mathscr A_k\longrightarrow
 \partial_{pp}+a_\infty^{-2}\partial_{qq}
 -\varepsilon a_\infty^{-3},
 \qquad a_\infty=(\lambda+2\Gamma_\infty)^{1/2},
\]
locally on $\mathbb S\times\mathbb R$.  Since the shifted domains
exhaust the whole cylinder, a subsequence of $\nu_k$ converges to a
bounded solution $\nu_\infty$ of
\[
 (\nu_\infty)_{pp}+a_\infty^{-2}(\nu_\infty)_{qq}
 -\varepsilon a_\infty^{-3}\nu_\infty=0
 \quad\text{in }\mathbb S\times\mathbb R.
\]
After taking $q_k\to q_\infty$ modulo periodicity, local uniform
convergence gives
$|\nu_\infty(q_\infty,0)|\geq\kappa$, so the limit is nonzero.
On the other hand, its Fourier coefficients satisfy
\[
 \nu_j''-
 \bigl(j^2a_\infty^{-2}+\varepsilon a_\infty^{-3}\bigr)\nu_j=0.
\]
No nonzero solution of this equation is bounded on the whole real
line.  Hence $\nu_\infty=0$, a contradiction.  Therefore the
convergence is global in $C^0$.

It remains to upgrade the convergence to the endpoint $X$ norm.
We do this as a Cauchy argument, so no endpoint regularity of the
local limit is assumed in advance.  Put
$W_{k\ell}=(\overline w_k-\overline w_\ell,
\underline w_k-\underline w_\ell)$ and use the Banach-space
mean-value formula at the fixed parameter $\lambda_k$:
\[
 \mathscr A_{k\ell}W_{k\ell}=R_{k\ell},
 \qquad
 \mathscr A_{k\ell}:=
 \int_0^1D_wF^\varepsilon
 \bigl(\lambda_k,w_\ell+t(w_k-w_\ell)\bigr)\,\dd t.
\]
Here
\[
 R_{k\ell}=y_k-y_\ell-
 \bigl(F^\varepsilon(\lambda_k,w_\ell)
       -F^\varepsilon(\lambda_\ell,w_\ell)\bigr),
\]
and analyticity on bounded admissible sets gives
\[
 \|R_{k\ell}\|_Y
 \le C\bigl(\|y_k-y_\ell\|_Y
             +|\lambda_k-\lambda_\ell|\bigr)\longrightarrow0.
\]
We spell out why the averaged segment retains uniform admissibility.
The map
$\lambda\mapsto a(\,\cdot\,;\lambda)^{-1}$ is uniformly Lipschitz on
the compact parameter interval containing the sequence.  Hence, for
all sufficiently large $k,\ell$, replacing $\lambda_\ell$ by
$\lambda_k$ changes each of the two ellipticity inequalities and the
top inequality by at most $\delta/4$.  Thus $w_\ell$, evaluated at
$\lambda_k$, has margins at least $3\delta/4$.  The quantities
$a^{-1}+w_p$ and $\lambda_k-2gw(\cdot,0)$ are affine along the
segment $w_\ell+t(w_k-w_\ell)$, so every point of that segment has
margins at least $\delta/2$.  Consequently, boundedness of $Q$
supplies uniform coefficient norms for $\mathscr A_{k\ell}$ and the
whole averaged family has common ellipticity and top-obliqueness
constants.  Equality of the principal traces at $p=p_0$ is preserved
along the segment, so the perfect-transmission complementing constant
is uniform as well.
Lemma~\ref{lem:global-schauder-fixed-eps} and the already proved
global $C^0$ convergence therefore yield
\[
 \|W_{k\ell}\|_X
 \le C\left(\|R_{k\ell}\|_Y+
             \|W_{k\ell}\|_{C^0}\right)\longrightarrow0.
\]
Thus $w_k$ is Cauchy in $X$.  Its $X$-limit agrees with the local
limit $w$, belongs to $Q$ because $Q$ is closed, and satisfies
$F^\varepsilon(\lambda,w)=y$ by continuity.  This proves the
assertion.
\end{proof}

With these properties of $F^{\varepsilon}$ established, all the
hypotheses of Theorem \ref{thm4.1} are verified.  Indeed,
$F^\varepsilon$ is real analytic on $\mathcal O_\delta$; Lemma
\ref{lem4.2} gives the Fredholm index-zero property at every solution,
Lemmas \ref{lem4.3} and \ref{lem4.4} give the one-dimensional kernel,
one-dimensional cokernel (by index zero), and transversality at
$(\lambda_*^\varepsilon,0,0)$, and Lemma \ref{lem4.5} gives the
bounded-closed compactness required in (H4).  These are precisely the
hypotheses of the analytic global bifurcation theorem.  In particular,
no generalized Fredholm degree, odd-crossing-number computation, or
global spectral description of the linearized operators is required.

Let
\[
 z_\delta^\varepsilon(s)
 =\bigl(\lambda_\delta^\varepsilon(s),
        \overline w_\delta^\varepsilon(s),
        \underline w_\delta^\varepsilon(s)\bigr),
 \qquad s\in\mathbb R,
\]
be the distinguished analytic route supplied by Theorem
\ref{thm4.1}, normalized so that
\[
 z_\delta^\varepsilon(0)=(\lambda_*^\varepsilon,0,0)
\]
and the positive local parameter points in the direction
$\Psi^\varepsilon(p)\cos q$.  We write
\[
 \mathscr R_\delta^\varepsilon
 :=z_\delta^\varepsilon(\mathbb R),\qquad
 \mathcal K_{\delta,+}^\varepsilon
 :=z_\delta^\varepsilon([0,\infty)),\qquad
 \mathcal K_{\delta,-}^\varepsilon
 :=z_\delta^\varepsilon(( -\infty,0]),
\]
and, when no confusion is possible,
$\mathcal K_\delta^\varepsilon:=\mathscr R_\delta^\varepsilon$.
This notation refers to the distinguished analytic route and not to
the maximal connected component of the complete solution set.

\begin{theorem}\label{thm4.6}
Suppose that the vorticity satisfies the hypotheses of Theorem
\ref{thm2.1}.  For every $0<\delta<\delta_*$ and every sufficiently
small $0<\varepsilon<\varepsilon_0$, there exists a distinguished
global analytic route
$\mathscr R_\delta^\varepsilon$ through
$(\lambda_*^\varepsilon,0,0)$.  It has the local expansion described
above and has a real-analytic reparametrization locally around each
of its points, in the sense of (C4).  Moreover, one of the following
alternatives holds:
\begin{enumerate}
\item for every bounded closed set $Q\subset\mathcal O_\delta$,
there exists $S_Q>0$ such that
\[
 z_\delta^\varepsilon(s)\notin Q
 \qquad\text{for all }|s|>S_Q;
\]
\item the distinguished route is periodic: there exists $T>0$ such that
\[
 z_\delta^\varepsilon(s+T)=z_\delta^\varepsilon(s)
 \qquad\text{for all }s\in\mathbb R.
\]
\end{enumerate}
\end{theorem}

The periodic alternative will not be discarded at this stage.  What
is needed for the singular limit is the localized boundary-point
property of Theorem \ref{thm4.10}.  There, if the route were trapped
inside a fixed bounded set, Lemma \ref{lem:qind-approx} excludes a
nonzero $q$-independent point for all sufficiently small
$\varepsilon$, and the two opposite local nodal orientations then
rule out a periodic return to the bifurcation point.

\subsection{The bifurcation structure of the approximating problem}

In this subsection, we study the oriented nodal pattern inherited
from the eigenfunction of the linearized problem at the bifurcation
point $(\lambda^{\varepsilon}_*,0,0)$.  The positive and negative
local half-routes have opposite orientations.  On any bounded part of
an oriented route on which no nonzero $q$-independent solution is
encountered, the maximum principle preserves the corresponding nodal
pattern.  This fact will be used in Theorem \ref{thm4.10} to exclude
the periodic alternative of Theorem \ref{thm4.6} whenever the route
is assumed to remain in a fixed bounded neighborhood.  Since strict
pointwise inequalities are not open or closed in the whole unweighted
space at infinite depth, we use the strengthened nodal conditions
(\ref{eq4.31})--(\ref{eq4.34}).

\begin{figure}
\centering
\IfFileExists{3.jpg}{\includegraphics[width=9cm]{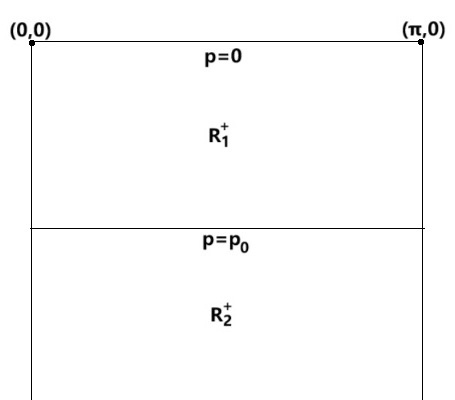}}{\fbox{\parbox[c][4cm][c]{8cm}{\centering Schematic of the nodal domains}}}
\caption{The nodal domain.} \label{fig3}
\end{figure}

Let's define
$$
R_1^+=(0,\pi)\times(p_0,0), \qquad R_2^+=(0,\pi)\times(-\infty,p_0),
$$
$$
\partial R_{1t}^+=(0,\pi)\times\{0\}, \qquad \partial R_{2t}^+=(0,\pi)\times\{p_0\},
$$
$$
\partial R_{1l}^+=\{(0,p):p\in(p_0,0)\}, \qquad \partial R_{2l}^+=\{(0,p):p\in(-\infty,p_0)\},
$$
$$
\partial R_{1r}^+=\{(\pi,p):p\in(p_0,0)\}, \qquad \partial R_{2r}^+=\{(\pi,p):p\in(-\infty,p_0)\}
$$
(see the Figure \ref{fig3}).
\par
Our goal is to show any nontrivial solution of (\ref{eq4.10}) in
$\mathcal{K}^{\varepsilon}_{\delta}$ possess the following nodal
pattern:
\begin{eqnarray}
w_q<0 \quad \text{in} \quad R_1^+\cup R_2^+ \cup \partial
R_{1t}^+\cup\partial R_{2t}^+, \label{eq4.31}
\end{eqnarray}
\begin{eqnarray}
w_{qq}<0~~\text{on}~~\partial R_{1l}^+\cup\partial R_{2l}^+, \quad
w_{qq}>0~~\text{on}~~\partial R_{1r}^+\cup\partial R_{2r}^+,
\label{eq4.32}
\end{eqnarray}
\begin{eqnarray}
w_{qq}(0,0)<0, \label{eq4.33}
\end{eqnarray}
\begin{eqnarray}
w_{qq}(\pi,0)>0. \label{eq4.34}
\end{eqnarray}
where
\begin{eqnarray}
w(q,p)=\left\{\begin{array}{ll}{\overline{w}(q,p)},  & {\text { for } p_0\leq p\leq0}, \\
{\underline{w}(q,p)}, & {\text { for } -\infty<p\leq
p_0}.\end{array}\right. \label{eq4.35}
\end{eqnarray}
Strict pointwise inequalities are not open in the whole
unweighted space $X$, because all derivatives vanish at infinite
depth.  They are, however, relatively open on the solution set of
the approximate equation.  The reason is the negative zeroth-order
term in the equation for $v=w_q$.  More precisely, if $P<p_0$ and
$v\leq0$ on $(0,\pi)\times\{P\}$, then the maximum principle on the
lower semi-cylinder, together with
\[
 v=0\quad\text{on }q=0,\pi,
 \qquad \lim_{p\to-\infty}v(q,p)=0,
\]
shows that $v\leq0$ on $(0,\pi)\times(-\infty,P]$.  If $v$ is not
identically zero, the strong maximum principle and the Hopf lemma
make the inequalities strict.  Thus a sign condition on one fixed
cross-section controls the entire cylindrical tail.  This is the
unbounded-domain maximum-principle device used in the nodal argument
of \cite[Section 4.2]{Hur11}; it avoids any unsupported relative
asymptotic expansion at $p=-\infty$.

Let $\mathcal N_\varepsilon^+$ denote the set of nontrivial solutions
satisfying \eqref{eq4.31}--\eqref{eq4.34}.  Let
$\mathcal T w(q,p)=w(q+\pi,p)$ and set
\[
 \mathcal N_\varepsilon^-:=\mathcal T\mathcal N_\varepsilon^+,
 \qquad
 \mathcal N_\varepsilon^{\pm}
 :=\mathcal N_\varepsilon^+\cup\mathcal N_\varepsilon^-.
\]
The two sets correspond to placing a crest or a trough at $q=0$.
For the limiting problem we use the analogous notation
$\mathcal N_0^{\pm}$, together with the uniform decay estimate of
Lemma \ref{lem5.2}.  Evenness and periodicity imply
\begin{eqnarray}
w_q=0 \quad \text{on} \quad \partial R_{1l}^+\cup\partial
R_{2l}^+\cup \partial R_{1r}^+\cup\partial R_{2r}^+. \label{eq4.36}
\end{eqnarray}

\begin{lemma}\label{lem4.7}
For all sufficiently small $s>0$, the positive local half-route
$z_\delta^\varepsilon(s)$ belongs to $\mathcal N_\varepsilon^+$.
For all sufficiently small $s<0$, the negative local half-route
belongs to $\mathcal N_\varepsilon^-$.  In particular, the two
punctured local half-routes have opposite nodal orientations.
\end{lemma}
\begin{proof}
Let
\[
 w(s,q,p)=s\Psi^\varepsilon(p)\cos q+r(s,q,p),
 \qquad \|r(s)\|_X=o(|s|),
\]
be the local analytic branch.  We first take $s>0$.  Choose the sign of the principal
eigenfunction so that $\Psi^\varepsilon>0$ on $(-\infty,0]$.
Fix $P<p_0$.  On the bounded transmission cylinder
$(0,\pi)\times[P,0]$, the local expansion gives
\[
 w_q(s,q,p)=-s\Psi^\varepsilon(p)\sin q+o(s)
 \quad\text{in }C^{2,\alpha}.
\]
Since $\Psi^\varepsilon$ has a positive minimum on $[P,0]$ and
$w_q/\sin q$ extends continuously to $q=0,\pi$, it follows that
\[
 w_q(s,q,p)<0
 \quad(0<q<\pi,\ P\leq p\leq0)
\]
for all sufficiently small $s>0$.

In the lower tail $D_P=(0,\pi)\times(-\infty,P)$, the function
$v=w_q$ satisfies the differentiated equation
\[
 \mathscr L_s v-\varepsilon a^{-3}v=0,
\]
where $\mathscr L_s$ is uniformly elliptic and has no positive
zeroth-order term.  On $q=0,\pi$ we have $v=0$, on $p=P$ we have
$v<0$, and $v\to0$ uniformly as $p\to-\infty$.  Applying the maximum
principle first on truncated cylinders and then letting the lower
boundary tend to $-\infty$ gives $v\leq0$ in $D_P$; the strong
maximum principle gives $v<0$.  This proves \eqref{eq4.31} on the
whole half-period.

The Hopf boundary lemma on $q=0$ and $q=\pi$, used separately in the
two layers and together with the transmission conditions at
$p=p_0$, yields
\[
 w_{qq}(0,p)<0,
 \qquad
 w_{qq}(\pi,p)>0,
\]
which is \eqref{eq4.32}.  At the two surface corners the local
expansion gives
\[
 w_{qq}(0,0)=-s\Psi^\varepsilon(0)+o(s)<0,
 \qquad
 w_{qq}(\pi,0)=s\Psi^\varepsilon(0)+o(s)>0.
\]
Thus \eqref{eq4.33}--\eqref{eq4.34} hold, and hence
$z_\delta^\varepsilon(s)\in\mathcal N_\varepsilon^+$ for small
$s>0$.  For $s<0$, the leading term is
$-|s|\Psi^\varepsilon(p)\cos q
=|s|\Psi^\varepsilon(p)\cos(q+\pi)$, so the same argument after the
half-period translation gives
$z_\delta^\varepsilon(s)\in\mathcal N_\varepsilon^-$.  We do not
include a sign condition on $w_{qqp}$ in the global nodal set, since
such a third-order corner sign is not needed later and is not
recovered by the Hopf--Serrin closedness argument.
\end{proof}

\begin{lemma}\label{lem4.8}
For $\sigma\in\{+,-\}$ let
\[
 \mathscr B_\varepsilon^\sigma
 :=\mathcal K_\delta^\varepsilon\cap
 \bigl(\mathbb R\times\mathcal N_\varepsilon^\sigma\bigr).
\]
Each $\mathscr B_\varepsilon^\sigma$ is relatively open among the
nontrivial $q$-dependent points of the distinguished route
$\mathcal K_\delta^\varepsilon$.  Moreover, if
$z_n\in\mathscr B_\varepsilon^\sigma$ and $z_n\to z$ in
$\mathbb R\times X$, where $z$ is a nontrivial $q$-dependent
solution, then $z\in\mathscr B_\varepsilon^\sigma$.  Consequently,
the nodal orientation is constant on every parameter interval of the
distinguished route containing only nontrivial $q$-dependent
solutions.  A relative boundary point of either oriented set is
$q$-independent (possibly a trivial point).
\end{lemma}
\begin{proof}
The two local half-branches belong to their corresponding oriented
sets by Lemma \ref{lem4.7}.  We first prove relative openness.  Fix a solution in
$\mathcal N_\varepsilon^+$; the case
$\mathcal N_\varepsilon^-$ follows by the half-period translation.
Choose $P<p_0$.  On the compact cross-section $p=P$, the quotient
$w_q/\sin q$ extends continuously to $[0,\pi]$ and is strictly
negative.  Hence the same sign holds for every sufficiently nearby
solution in $X$.  The maximum principle on
$(0,\pi)\times(-\infty,P]$ propagates this sign through the entire
lower tail, while on the bounded cylinder above $p=P$ all strict
inequalities persist by ordinary $C^{3,\alpha}$ continuity.  Thus
$\mathscr B_\varepsilon^+$ is relatively open among nontrivial
$q$-dependent solutions.

For the sequential closedness assertion, first take $\sigma=+$ and
write $v_n=(w_n)_q$ and $v=w_q$ for the limit.  Then $v\leq0$ in the
positive half-period and $v=0$ on $q=0,\pi$.  If $v\not\equiv0$, the
strong maximum principle gives strict negativity in each layer.  A
zero on the transmission interface would, by the Hopf lemma from the
two sides, imply
\[
 v_p^+(q,p_0)<0<v_p^-(q,p_0),
\]
contrary to the transmission condition $v_p^+=v_p^-$ at $p=p_0$.
The oblique boundary lemma on the open surface edge gives
$v<0$ there, and the Hopf lemma on the open vertical edges gives
\[
 v_q(0,p)<0,\qquad v_q(\pi,p)>0\qquad(p<0),
\]
with the same conclusion across $p=p_0$ by the two-sided
transmission argument.  It remains only to justify that strictness
cannot be lost at the two surface corners.  Serrin's corner lemma,
applied in a quarter-disc at $(0,0)$, gives the alternative
\[
 v_q(0,0)<0\quad\hbox{or}\quad v_{pq}(0,0)>0.
\]
At $(\pi,0)$ the reflected alternative is
\[
 v_q(\pi,0)>0\quad\hbox{or}\quad v_{pq}(\pi,0)<0.
\]
These alternatives are strict because $v<0$ in the half-period and
$v=0$ on its vertical sides.  To eliminate the second alternative
when the first derivative vanishes, put
\[
 B:=\lambda^{-1/2}+w_p>0\qquad\hbox{on }p=0.
\]
The surface condition differentiated once with respect to $q$ is
\[
 gB^2v+(2gw-\lambda)Bv_p+w_qv_q=0.
\]
At either corner, evenness and periodicity give
$v=w_q=0$, and the displayed equation gives $v_p=0$.  Differentiating
it once more tangentially and evaluating at the corner therefore
yields
\[
 gB^2v_q+(2gw-\lambda)Bv_{pq}+v_q^2=0.
\]
If $v_q=0$, then the original surface condition gives
$(2gw-\lambda)B=-B^{-1}\ne0$, so
the last identity forces $v_{pq}=0$, contradicting
the corresponding strict corner alternative.  Consequently
\[
 v_q(0,0)<0,\qquad v_q(\pi,0)>0.
\]
Thus all of \eqref{eq4.31}--\eqref{eq4.34} are recovered.  No loss
can escape to $p=-\infty$, because the tail maximum principle uses
the uniform condition $v_n\to0$ there.  Hence
$z\in\mathscr B_\varepsilon^+$.  The proof for $\sigma=-$ is the
same after the half-period translation.
If $v\equiv0$, the limiting solution is $q$-independent, which is
precisely the only possible boundary mechanism left by the
argument.
\end{proof}

\begin{lemma}\label{lem4.9}
Assume that there exists a sequence
\[
 (\lambda_k,\overline w_k,\underline w_k)
 \in\mathcal K_\delta^\varepsilon
 \cap\mathcal N_\varepsilon^{\pm}
\]
of nontrivial solutions converging in $\mathbb R\times X$ to a
trivial solution $(\lambda,0,0)$.  Then
$\lambda=\lambda_*^\varepsilon$.
\end{lemma}
\begin{proof}
Let
$\{(\lambda_k,\overline w_k,\underline w_k)\}$ be nontrivial
solutions in $\mathcal K_\delta^\varepsilon\cap\mathcal N_\varepsilon^{\pm}$
converging in $\mathbb R\times X$ to $(\lambda,0,0)$, and let $w_k$
denote the corresponding piecewise function.  Since
$(w_k)_q<0$ in the half-period, $(w_k)_q\not\equiv0$.  Define
\[
 v_k:=\frac{(w_k)_q}
 {\|(w_k)_q\|_{C^{2,\alpha}(\overline R_1)
 \times C^{2,\alpha}(\overline R_2)}}.
\]
Put
\[
 A_k:=a(p;\lambda_k)^{-1}+(w_k)_p.
\]
Differentiating the two bulk equations and the surface condition
with respect to $q$ gives, in each layer,
\begin{equation}\label{eq:vk-nonlinear}
 (1+(w_k)_q^2)(v_k)_{pp}
 -2A_k(w_k)_q(v_k)_{pq}
 +A_k^2(v_k)_{qq}
 +b_{1,k}(v_k)_p+b_{2,k}(v_k)_q
 -\varepsilon a(p;\lambda_k)^{-3}v_k=0,
\end{equation}
where
\begin{align*}
 b_{1,k}&=-2(w_k)_q(w_k)_{qp}
 +2A_k(w_k)_{qq}
 +3\gamma(-p)A_k^2,\\
 b_{2,k}&=2(w_k)_q(w_k)_{pp}
 -2A_k(w_k)_{qp}
 -2\gamma(-p)a(p;\lambda_k)^{-3}(w_k)_q.
\end{align*}
On $p=0$,
\begin{equation}\label{eq:vk-boundary}
 gB_k^2v_k+(2gw_k-\lambda_k)B_k(v_k)_p
 +(w_k)_q(v_k)_q=0,
 \qquad
 B_k:=\lambda_k^{-1/2}+(w_k)_p.
\end{equation}
Moreover $v_k^+=v_k^-$ and $(v_k)_p^+=(v_k)_p^-$ on $p=p_0$,
and $v_k,Dv_k\to0$ as $p\to-\infty$.

The homogeneous Schauder estimate for
\eqref{eq:vk-nonlinear}--\eqref{eq:vk-boundary}, with coefficients
converging to those at $(\lambda,0,0)$, gives
\[
 1=\|v_k\|_{C^{2,\alpha}}
 \leq C\|v_k\|_{C^0}.
\]
Thus $\|v_k\|_{C^0}\geq C^{-1}$.  If all points at which this
$C^0$ mass is attained escaped to $p=-\infty$, the translation
argument of Lemma \ref{lem4.2} would produce a nonzero bounded
solution on $\mathbb S\times\mathbb R$ of
\[
 V_{pp}+a_\infty^{-2}V_{qq}
 -\varepsilon a_\infty^{-3}V=0,
\]
which is impossible.  Hence a fixed finite subcylinder carries a
uniform amount of $C^0$ mass.  Local Schauder compactness therefore
gives, after passing to a subsequence,
\[
 v_k\longrightarrow v
 \quad\text{in }C^{2,\beta}_{\rm loc}(\overline R_1)
 \times C^{2,\beta}_{\rm loc}(\overline R_2),
 \qquad 0<\beta<\alpha,
\]
with $v\not\equiv0$.  The global decay follows from the equation and
the fixed-$\varepsilon$ tail estimate.  Passing to the limit in
\eqref{eq:vk-nonlinear}--\eqref{eq:vk-boundary}, we obtain
\begin{equation}\label{eq:linear-v-limit}
 \left\{
 \begin{aligned}
 &v_{pp}+a^{-2}v_{qq}+3\gamma(-p)a^{-2}v_p
   -\varepsilon a^{-3}v=0
   &&\text{in }R_1\cup R_2,\\
 &gv-\lambda^{3/2}v_p=0
   &&\text{on }p=0,\\
 &[v]_{p=p_0}=[v_p]_{p=p_0}=0,
 \qquad v,Dv\to0 &&\text{as }p\to-\infty,
 \end{aligned}
 \right.
\end{equation}
where $a=a(p;\lambda)$.  Since $v_k<0$ in the half-period, the
strong maximum principle and the transmission Hopf argument give
\[
 v(q,p)<0\qquad(0<q<\pi,\ -\infty<p\leq0).
\]

The function $v$ is odd and $2\pi$-periodic, so
\[
 v(q,p)=\sum_{j=1}^\infty m_j(p)\sin(jq).
\]
For each $j\geq1$, equation \eqref{eq:linear-v-limit} yields
\[
 -(a^3m_j')'+\varepsilon m_j=-j^2a m_j,
 \qquad
 \lambda^{3/2}m_j'(0)=gm_j(0),
\]
with the transmission and decay conditions.  In particular,
\[
 m_1(p)=\frac2\pi\int_0^\pi v(q,p)\sin q\,dq<0,
\]
so $m_1\not\equiv0$ and $-1$ is a generalized eigenvalue.  Hence
\[
 \mu^\varepsilon(\lambda)\leq-1.
\]
If $\mu^\varepsilon(\lambda)<-1$, let $\Phi>0$ be the principal
eigenfunction corresponding to $\mu^\varepsilon(\lambda)$.  The
self-adjoint Sturm--Liouville equations for $\Phi$ and $m_1$ imply
the weighted orthogonality relation
\[
 \int_{-\infty}^0 a(p;\lambda)\Phi(p)m_1(p)\,dp=0.
\]
This is impossible because $a>0$, $\Phi>0$, and $m_1<0$.  Therefore
$\mu^\varepsilon(\lambda)=-1$.  The strict monotonicity and
uniqueness established in Lemma \ref{lem4.3} give
$\lambda=\lambda_*^\varepsilon$.
\end{proof}

The next two lemmas separate the exact problem from its Fredholm
approximations.  For the exact equation every $q$-independent
solution is laminar.  For a fixed $\varepsilon>0$ the regularized
problem may possess additional $q$-independent solutions; what is
needed for the singular limit, and what is true, is that such
solutions cannot remain in a fixed bounded subset as
$\varepsilon\downarrow0$.

\begin{lemma}[Exact $q$-independent solutions]\label{lem:qind-exact}
Let $(\lambda,\overline w,\underline w)\in\mathcal O_\delta$ solve
$F(\lambda,\overline w,\underline w)=0$.  If the corresponding
piecewise function $w$ is independent of $q$, then $w\equiv0$.
\end{lemma}
\begin{proof}
Put $h=H(\cdot;\lambda)+w$.  In each layer,
\[
 h_{pp}+\gamma(-p)h_p^3=0,
 \qquad h_p>0.
\]
Consequently
\[
 (h_p^{-2})_p=2\gamma(-p)=(a^2)_p.
\]
The transmission conditions make $h_p^{-2}-a^2$ continuous across
$p=p_0$, while the condition at infinite depth gives
$h_p^{-2}-a^2\to0$.  Hence $h_p^{-2}=a^2$ on $(-\infty,0]$.
Since both quantities are positive, $h_p=a^{-1}=H_p$.  Thus $w_p=0$, so $w\equiv C$ for some constant $C$.  The
surface condition at $p=0$ now reads
\[
 1+(2gC-\lambda)\lambda^{-1}=0,
\]
and therefore $C=0$.  Hence $w\equiv0$.
\end{proof}

\begin{lemma}[Exclusion of bounded regularized shear solutions]
\label{lem:qind-approx}
Let $Q$ be a bounded set with $\overline Q\subset\mathcal O_\delta$.
There exists $\varepsilon_Q>0$ such that, for
$0<\varepsilon<\varepsilon_Q$, the set $Q$ contains no nonzero
$q$-independent solution of $F^\varepsilon=0$.
\end{lemma}
\begin{proof}
Suppose otherwise.  Then $\varepsilon_n\downarrow0$ and there are
nonzero $q$-independent solutions
$(\lambda_n,W_n)\in Q$.  Here $W_n$ denotes the piecewise function
formed from the two layers.  Set
\[
 A_n(p)=a(p;\lambda_n)^{-1}.
\]
Factoring the cubic term in the bulk equation gives
\begin{equation}\label{eq:qind-linear}
 W_n''+b_n(p)W_n'
 -\varepsilon_n a(p;\lambda_n)^{-3}W_n=0,
\end{equation}
where
\[
 b_n(p)=\gamma(-p)
 \left(3A_n^2+3A_nW_n'+(W_n')^2\right).
\]
The boundedness of $Q$ and the decay assumption on $\gamma$ imply
\[
 |b_n(p)|\leq C|\gamma(-p)|,
 \qquad
 \sup_n\|b_n\|_{L^1(-\infty,0)}<\infty.
\]

At $p=0$ the nonlinear boundary condition reads
\[
 2gW_n(0)=\lambda_n-
 \left(\lambda_n^{-1/2}+W_n'(0)\right)^{-2}.
\]
If $W_n(0)=0$, then $W_n'(0)=0$ and uniqueness for
\eqref{eq:qind-linear} gives $W_n\equiv0$, a contradiction.  Thus
$W_n(0)\neq0$.  By the mean-value theorem there is a number $\xi_n$
between $\lambda_n^{-1/2}$ and
$\lambda_n^{-1/2}+W_n'(0)$ such that
\begin{equation}\label{eq:qind-boundary-ratio}
 \frac{W_n'(0)}{W_n(0)}=g\xi_n^3.
\end{equation}
Because $\overline Q\subset\mathcal O_\delta$ and $Q$ is bounded,
both endpoints defining $\xi_n$ are bounded below by a positive
constant.  Therefore
\begin{equation}\label{eq:qind-ratio-lower}
 \frac{W_n'(0)}{W_n(0)}\geq c_Q>0.
\end{equation}

Equation \eqref{eq:qind-linear} has negative zeroth-order
coefficient.  The one-dimensional maximum principle, together with
$W_n(-\infty)=0$ and \eqref{eq:qind-boundary-ratio}, shows that
$W_n$ has the sign of $W_n(0)$ on the whole half-line.  If
\[
 M_n(p)=\exp\left(\int_0^p b_n(s)\,\dd s\right),
\]
then
\[
 (M_nW_n')'
 =\varepsilon_n a(p;\lambda_n)^{-3}M_nW_n.
\]
Since $M_n$ is uniformly bounded above and below and
$W_n'(-\infty)=0$, $W_n'$ has the same sign as $W_n$.  Hence the
normalized function
\[
 u_n(p)=\frac{W_n(p)}{W_n(0)}
\]
satisfies
\[
 0<u_n\leq1,\qquad u_n'\geq0,\qquad u_n(0)=1.
\]

The coefficients are locally precompact, and the equation gives
local $C^2$ bounds for $u_n$.  Passing to a subsequence,
$u_n\to u$ in $C^1_{\rm loc}(-\infty,0]$ and
$b_n\to b$ locally, where $b\in L^1(-\infty,0)$.  The limit obeys
\[
 u''+bu'=0,
 \qquad 0\leq u\leq1,
 \qquad u'\geq0,
 \qquad u(0)=1.
\]
Writing $M(p)=\exp(\int_0^p b)$ gives $Mu'=C$.  If $C>0$, then,
since $M$ and $M^{-1}$ are bounded on $(-\infty,0]$,
\[
 u(p)=1-C\int_p^0M(s)^{-1}\,\dd s
\]
becomes negative for $p\ll-1$, a contradiction.  Thus $C=0$ and
$u'(0)=0$.  On the other hand,
\eqref{eq:qind-ratio-lower} gives $u_n'(0)\geq c_Q$, contradicting
$C^1$ convergence at $p=0$.
\end{proof}

Lemma \ref{lem4.8} shows that, along an oriented half-route, the
strict nodal pattern can be lost only through a $q$-independent
solution, while Lemma \ref{lem4.9} identifies the parameter of every
nodal approach to the trivial line.  The analytic alternative in
Theorem \ref{thm4.6}, together with Lemma \ref{lem:qind-approx}, now
gives the localized boundary-point statement required for the
Whyburn limit.  We first isolate the global nodal argument which
excludes the periodic alternative.  This is the analogue, in the present
height/transmission setting, of the nodal loop exclusion in
\cite[Section 4]{ConstantinSV}.

\begin{lemma}[Local analytic opposite-germ passage]
\label{lem:no-retracing-return}
Let $z:\mathbb R\to\mathscr R$ be the distinguished continuation in
Theorem~\ref{thm4.1}, let $b=z(0)$, and suppose that
\[
 \tau:=\min\{t>0:z(t)=b\}
\]
exists.  Assume that $z(t)$ is nontrivial for $0<t<\tau$.  If the
outgoing germ $z((0,\eta))$ is one of the two nontrivial local
Crandall--Rabinowitz half-germs at $b$, then, for some $\eta_\tau>0$,
the incoming germ $z((\tau-\eta_\tau,\tau))$ is the opposite
half-germ.
\end{lemma}
\begin{proof}
By (C3), there is a neighborhood $V$ of $b$ in which the nontrivial
part of the distinguished route consists of exactly two local
Crandall--Rabinowitz half-germs, say $\mathscr G_+$ and
$\mathscr G_-$.  Shrinking $V$ if necessary, the outgoing segment
lies in $\mathscr G_+$.  Continuity and the definition of the first
return allow $\eta_\tau>0$ to be chosen so that the incoming segment
lies in $V\setminus\{b\}$.  It lies entirely in one of
$\mathscr G_+$ and $\mathscr G_-$: otherwise its connected image
would pass through their only common point $b$ before time $\tau$.

Suppose that the incoming segment also lies in $\mathscr G_+$.  Both
segments approach $b$ through the same local half-arc, and, after
they are shortened, their images contain a common nontrivial
relatively open analytic subarc.  The distinguished continuation
would therefore return to $b$ by retracing the very half-arc along
which it left.  This is precisely the possibility excluded by the
local real-analytic reparametrization in (C4): in the local
uniformizing parameter, the two sides joined by the continuation are
the two different sides of the distinguished analytic arc, not the
same half-arc traversed twice.  This is the same C4 argument used to
exclude equality of the initial and incoming local arcs in
\cite[Section 4]{ConstantinSV}.

Thus the incoming segment cannot lie in $\mathscr G_+$ and must lie
in $\mathscr G_-$.  Notice that this proof invokes neither a discrete
singular set nor any nonvanishing condition on the derivative of the
physical parameter.
\end{proof}

\begin{lemma}[Global nodal exclusion of a periodic distinguished route]
\label{lem:global-nodal-loop}
Fix $\varepsilon>0$ and let $Q\subset\mathcal O_\delta$ be bounded and
closed.  Suppose that $Q$ contains no nonzero $q$-independent solution
of $F^\varepsilon=0$.  Then the distinguished route through
$(\lambda_*^\varepsilon,0,0)$ cannot be periodic while remaining in
$Q$.
\end{lemma}
\begin{proof}
Write $b=(\lambda_*^\varepsilon,0,0)$ and suppose that
$z:\mathbb R\to Q$ is periodic, with some period $T>0$.  Orient its
parameter so that the positive local germ at $b$ is the
Crandall--Rabinowitz germ in $\mathcal N_\varepsilon^+$.  The local
expansion (C2)--(C3) gives $\eta_0>0$ such that $z(t)\ne b$ for
$0<t\leq\eta_0$.  Since $z(T)=b$, the closed return set in
$[\eta_0,T]$ is nonempty and has a minimum.  Thus there is a first
return time
\[
 T_0:=\min\{t>0:z(t)=b\}>0.
\]

We claim that
\begin{equation}\label{eq:periodic-route-positive-orientation}
 z(t)\in\mathcal N_\varepsilon^+\qquad(0<t<T_0).
\end{equation}
This holds initially by Lemma~\ref{lem4.7}.  If $t_*<T_0$ were the
right endpoint of the component of
$\{t\in(0,T_0):z(t)\in\mathcal N_\varepsilon^+\}$ containing the
initial local segment, choose $t_n\uparrow t_*$.  If $z(t_*)$ were
nontrivial and $q$-dependent, the sequential closedness in
Lemma~\ref{lem4.8} would keep $z(t_*)$ in
$\mathcal N_\varepsilon^+$, and relative openness would extend the
same orientation past $t_*$, a contradiction.  Thus $z(t_*)$ would
be $q$-independent.  The hypothesis on $Q$ excludes a nonzero such
state, while Lemma~\ref{lem4.9} identifies a trivial nodal limit with
$b$, contradicting the definition of $T_0$.  Therefore no such
$t_*<T_0$ exists, proving
\eqref{eq:periodic-route-positive-orientation}.

The hypotheses of Lemma~\ref{lem:no-retracing-return} now hold on
$(0,T_0)$.  It follows that the incoming germ at $T_0$ is the local
half-germ opposite to the outgoing positive germ.  Lemma~\ref{lem4.7}
places that opposite germ in $\mathcal N_\varepsilon^-$.  Hence
$z(t)\in\mathcal N_\varepsilon^-$ for $t<T_0$ sufficiently close to
$T_0$, contradicting
\eqref{eq:periodic-route-positive-orientation}.  This argument uses
only the two local half-germs in (C3), the local analytic
no-retracing conclusion of (C4), and the open--closed nodal argument;
it makes no assertion about the singular set of the linearization or
about the entire solution set being a finite analytic graph.
\end{proof}

\begin{theorem}\label{thm4.10}
Let $U$ be a bounded open subset of $\mathcal O_\delta$ such that
\[
 (\lambda_*^0,0,0)\in U,
 \qquad \overline U\subset\mathcal O_\delta.
\]
For all sufficiently small $\varepsilon>0$,
$(\lambda_*^\varepsilon,0,0)\in U$ and
\[
 \partial U\cap\mathcal K_{\delta,+}^\varepsilon
 \cap\bigl(\mathbb R\times\mathcal N_\varepsilon^+\bigr)
 \neq\varnothing.
\]
Consequently,
\[
 \partial U\cap\mathcal K_\delta^\varepsilon
 \cap\bigl(\mathbb R\times\mathcal N_\varepsilon^{\pm}\bigr)
 \neq\varnothing.
\]
\end{theorem}
\begin{proof}
The convergence $\lambda_*^\varepsilon\to\lambda_*^0$ gives the
first assertion.  By Lemma \ref{lem:qind-approx}, after decreasing
$\varepsilon$ if necessary, $\overline U$ contains no nonzero
$q$-independent solution of the approximate problem.

Suppose, for contradiction, that
\[
 \partial U\cap\mathcal K_{\delta,+}^\varepsilon=\varnothing.
\]
Since the positive half-route starts at
$(\lambda_*^\varepsilon,0,0)\in U$ and is continuous, the assumption
$\partial U\cap\mathcal K_{\delta,+}^\varepsilon=\varnothing$ forces
it to remain in $U$ for all positive values of its global parameter.
Because $\overline U\subset\mathcal O_\delta$, the set
$\overline U$ is an admissible bounded closed set in the escape
alternative of Theorem \ref{thm4.6}.  That alternative is therefore
impossible, and the analytic route must be periodic.

The periodic alternative is now ruled out by the preceding global
nodal lemma.  Indeed, the periodic distinguished route contains a
closed loop through $(\lambda_*^\varepsilon,0,0)$, its positive local
germ is the Crandall--Rabinowitz half-branch contained in
$\mathcal N_\varepsilon^+$ by Lemma~\ref{lem4.7}, and
$\overline U$ contains no nonzero $q$-independent regularized
solution.  Applying Lemma~\ref{lem:global-nodal-loop} with
$Q=\overline U$ gives a contradiction.  Thus the periodic
alternative in Theorem~\ref{thm4.6} is impossible while the positive
route is trapped in $\overline U$.

Hence $\mathcal K_{\delta,+}^\varepsilon$ meets $\partial U$.
Let $s_U>0$ be its first exit time.  On $(0,s_U)$ the positive nodal
orientation is preserved.  Indeed, at a first endpoint of the
positive oriented interval, Lemma~\ref{lem4.8} leaves only a
$q$-independent state.  A nonzero such state is excluded in
$\overline U$; a trivial state is $b$ by Lemma~\ref{lem4.9}, and
Lemma~\ref{lem:no-retracing-return} would then put the incoming germ
in $\mathcal N_\varepsilon^-$, contradicting its positive
orientation.  The sequential closedness in Lemma~\ref{lem4.8} and
the same exclusion of $q$-independent states show that the exit point
itself lies in $\mathcal N_\varepsilon^+$.  Therefore
\[
 \partial U\cap\mathcal K_{\delta,+}^\varepsilon
 \cap(\mathbb R\times\mathcal N_\varepsilon^+)\ne\varnothing,
\]
which proves the theorem.
\end{proof}

\section{The singular limit in the natural deep-water topology}
\label{sec:singular-dw}

The Fredholm regularization used in Section~4 is an auxiliary device.
For fixed $\varepsilon>0$ its zero Fourier mode decays at the rate
$O(\sqrt\varepsilon)$ at the infinite bottom.  Consequently, fixed-$\varepsilon$
decay of the absolute height perturbation cannot be made uniform as
$\varepsilon\downarrow0$.  The singular limit must therefore be taken in a
topology adapted to the physical deep-water condition, which prescribes the
derivatives $h_p$ and $h_q$ rather than an absolute additive height level.

For a $2\pi$-periodic function $f$ write
\[
 \langle f\rangle(p):=\frac1{2\pi}\int_{-\pi}^{\pi}f(q,p)\,\dd q.
\]
For the piecewise function $w$ represented by
$(\overline w,\underline w)$, set for every $p\le0$
\[
 m(p):=\langle w\rangle(p),\qquad
 \widetilde w(q,p):=w(q,p)-m(p),\qquad
 A(q,p):=a(p;\lambda)^{-1}+w_p(q,p).
\]
Thus these quantities are defined on both layers and have common traces
at $p=p_0$.  When we work on the lower tail $p<p_0$, the symbols
$w,\widetilde w,A$ mean their lower-layer representatives
$\underline w,\underline{\widetilde w},a^{-1}+\underline w_p$.

\subsection{Uniform transverse decay and the mean-flux identity}

We first record the precise uniform Phragm\'en--Lindel\"of estimate
used below.  The point of the statement is that its decay exponent
depends only on the ellipticity ratio, the drift bound, and the width
of the cylinder.  In particular, it is unaffected by a nonpositive
zero-order term which may converge to zero.

\begin{lemma}[Uniform Phragm\'en--Lindel\"of estimate]
\label{lem:uniform-cylinder-PL}
Let $\ell>0$, $P_0\in\mathbb R$, and
\[
 \mathscr LU=a^{qq}U_{qq}+2a^{qp}U_{qp}+a^{pp}U_{pp}
 +b^qU_q+b^pU_p+dU
\]
on $\mathcal C=(0,\ell)\times(-\infty,P_0)$.  Suppose that
the principal matrix is symmetric and, for some $\kappa\in(0,1)$
and $B>0$,
\begin{equation}\label{eq:uniform-PL-coefficients}
 \kappa |\xi|^2\le a^{ij}(q,p)\xi_i\xi_j
 \le\kappa^{-1}|\xi|^2,
 \qquad |b(q,p)|\le B,
 \qquad d(q,p)\le0.
\end{equation}
Let $U$ be a bounded classical solution of $\mathscr LU=0$ which
has one sign, satisfies
\[
 U(0,p)=U(\ell,p)=0,
 \qquad
 \sup_{0\le q\le\ell}|U(q,p)|\longrightarrow0
 \quad(p\to-\infty),
\]
and, on the finite cross-section, obeys
\begin{equation}\label{eq:uniform-PL-cross-section}
 |U(q,P_0)|\le M_0\operatorname{dist}(q,\{0,\ell\}).
\end{equation}
Then there are $C,\sigma>0$, depending only on
$\kappa,B,$ and $\ell$, such that
\begin{equation}\label{eq:uniform-PL-conclusion}
 |U(q,p)|\le CM_0 e^{\sigma(p-P_0)}
 \qquad(0\le q\le\ell,\ p\le P_0).
\end{equation}
The same constants work for every family satisfying
\eqref{eq:uniform-PL-coefficients}.
\end{lemma}
\begin{proof}
Set $B_1=B+2\kappa^{-1}$.  The one-dimensional Dirichlet
extremal problem
\begin{equation}\label{eq:uniform-PL-one-dimensional}
 \kappa\phi''+B_1|\phi'|+\lambda_D\phi=0
 \quad\hbox{in }(0,\ell),
 \qquad \phi(0)=\phi(\ell)=0,
\end{equation}
has a positive principal pair $(\lambda_D,\phi)$ with
$\lambda_D>0$.  For completeness, it can be obtained by shooting
on $(0,\ell/2)$ for
$\kappa\phi''+B_1\phi'+\lambda\phi=0$, with
$\phi(0)=0$, $\phi'(0)=1$, choosing the first value of $\lambda$
for which $\phi'(\ell/2)=0$, and then reflecting about
$q=\ell/2$.  The resulting function is strictly concave on both
half-intervals.  The one-dimensional Hopf lemma and normalization
$\max\phi=1$ give
\begin{equation}\label{eq:uniform-PL-phi-distance}
 c_0\operatorname{dist}(q,\{0,\ell\})
 \le\phi(q)\le
 C_0\operatorname{dist}(q,\{0,\ell\}).
\end{equation}

Choose $0<\sigma\le1$ so small that
\begin{equation}\label{eq:uniform-PL-sigma-choice}
 \kappa^{-1}\sigma^2+B\sigma\le\frac{\lambda_D}{2}
\end{equation}
and put $Z(q,p)=e^{\sigma(p-P_0)}\phi(q)$.  Since
$\phi''<0$ away from its midpoint, ellipticity and
\eqref{eq:uniform-PL-one-dimensional} give
\[
 \begin{aligned}
 e^{-\sigma(p-P_0)}\mathscr LZ
 ={}&a^{qq}\phi''+(2\sigma a^{qp}+b^q)\phi'
 +(a^{pp}\sigma^2+b^p\sigma+d)\phi\\
 \le{}&\kappa\phi''+B_1|\phi'|
 +(\kappa^{-1}\sigma^2+B\sigma)\phi
 \le-\frac{\lambda_D}{2}\phi.
 \end{aligned}
\]
The inequality at the midpoint follows by approximation, or in the
viscosity sense, and this is sufficient for comparison.

Assume first that $U\ge0$.  By
\eqref{eq:uniform-PL-cross-section} and
\eqref{eq:uniform-PL-phi-distance}, choose
$C_1=C_1(\kappa,B,\ell)$ so that
$U(q,P_0)\le C_1M_0\phi(q)$.  Fix $\eta>0$ and truncate at
$p=-R$.  For all sufficiently large $R$, the end condition gives
$U(q,-R)\le\eta$.  The function
\[
 U-C_1M_0Z-\eta
\]
is nonpositive on the boundary of the truncated rectangle, while
its image under $\mathscr L$ is nonnegative because
$\mathscr LZ<0$ and $d\le0$.  The weak maximum principle therefore
gives $U\le C_1M_0Z+\eta$.  Let first $R\to\infty$ and then
$\eta\downarrow0$.  This proves
\eqref{eq:uniform-PL-conclusion}; applying the same argument to
$-U$ treats the other sign.
\end{proof}

\begin{lemma}[Uniform decay of the nonzero Fourier modes]
\label{lem5.2}
Fix $\delta>0$ and $M>0$.  There exist constants $C,\sigma>0$,
independent of $0\leq\varepsilon\leq\varepsilon_0$, with the following
property.  Let $(\lambda,\overline w,\underline w)$ be an oriented exact
solution or an oriented regularized solution satisfying the admissibility
margins and
\[
 |\lambda|+\|\overline w\|_{C^{3,\alpha}(\overline R_1)}
 +\|\underline w_p\|_{C^{2,\alpha}(\overline R_2)}
 +\|\underline w_q\|_{C^{2,\alpha}(\overline R_2)}\le M.
\]
Then, for $p\le p_0-2$,
\begin{equation}\label{eq5dw:wq-tail}
 \sum_{i+j\le2}\bigl|\partial_p^i\partial_q^j\underline w_q(q,p)\bigr|
 \le Ce^{\sigma p}.
\end{equation}
Consequently,
\begin{equation}\label{eq5dw:osc-tail}
 \|\widetilde w\|_{C^{2,\alpha}(S_s)}
 +\|\widetilde w_p\|_{C^{1,\alpha}(S_s)}
 +\|\widetilde w_{pp}\|_{C^{\alpha}(S_s)}
 \le Ce^{\sigma s},\qquad s<p_0-2,
\end{equation}
where $S_s=\mathbb S\times[s-1,s+1]$.
\end{lemma}
\begin{proof}
Put $v=\underline w_q$ and $A=a^{-1}+\underline w_p$ on the lower
layer.  Horizontal differentiation of the bulk equation gives
\begin{equation}\label{eq5dw:linear-v}
 (\mathscr L-\varepsilon a^{-3})v=0.
\end{equation}
Here
\[
 \mathscr Lv=(1+\underline w_q^2)v_{pp}
 -2A\underline w_qv_{pq}+A^2v_{qq}+b_1v_p+b_2v_q,
\]
with
\[
 b_1=-2\underline w_q\underline w_{pq}
     +2A\underline w_{qq}+3\gamma(-p)A^2,
 \quad
 b_2=2\underline w_q\underline w_{pp}
     -2A\underline w_{pq}-2\gamma(-p)a^{-3}\underline w_q.
\]
The admissibility margins and the displayed derivative bound give a
common ellipticity ratio $\kappa$ and a common bound $B_0$ for the
first-order coefficients; the zero-order coefficient is
nonpositive.  On the half strip
$S^+=(0,\pi)\times(-\infty,p_0)$, orientation gives $v<0$ and
$v=0$ on $q=0,\pi$.

Take $P_0=p_0-1$.  Since $v=0$ at $q=0,\pi$, the assumed
$C^{2,\alpha}$ bound and the mean-value theorem give
\[
 |v(q,P_0)|\le C\operatorname{dist}(q,\{0,\pi\}).
\]
The endpoint condition in $X_{\rm dw}$ (and, for a regularized
solution, the stronger endpoint condition in $X$) gives
$v(\cdot,p)\to0$ uniformly as $p\to-\infty$.  All hypotheses of
Lemma~\ref{lem:uniform-cylinder-PL} are therefore satisfied by
$-v$ on the positive half-period and, by reflection, by $v$ on the
negative half-period.  We obtain
\[
 |v(q,p)|\le Ce^{\sigma p},\qquad p\le p_0-2.
\]
The constants depend only on $\delta,M$ and the fixed data.  The
negative half-period follows from evenness and periodicity.

Finally cover the lower end by unit strips.  Interior and side-boundary
Schauder estimates for the homogeneous linear equation
\eqref{eq5dw:linear-v}, whose coefficient norms are uniform, give
\[
 \|v\|_{C^{2,\alpha}(S_s)}
 \le C\|v\|_{C^0(S_s^*)},
\]
where $S_s^*$ is a fixed enlargement of $S_s$.  After decreasing
$\sigma$ this proves \eqref{eq5dw:wq-tail}.  The functions
$\widetilde w$, $\widetilde w_p$ and $\widetilde w_{pp}$ have zero
horizontal mean.  Poincar\'e's inequality in $q$, together with the
bounds just obtained for $w_q,w_{pq},w_{ppq}$, yields
\eqref{eq5dw:osc-tail}.
\end{proof}

\begin{lemma}[Exact mean-flux identity]
\label{lem5.3-tail}
For every regularized solution $F^\varepsilon=0$ define
\begin{equation}\label{eq5dw:Jdef}
 \mathfrak J(p):=
 \frac{a(p;\lambda)^2}{2}
 -\left\langle\frac{1+w_q^2}{2A^2}\right\rangle(p).
\end{equation}
Then
\begin{equation}\label{eq5dw:Jprime}
 \mathfrak J'(p)=
 \varepsilon a(p;\lambda)^{-3}
 \left\langle\frac{w}{A^3}\right\rangle(p),
\end{equation}
and the top boundary condition gives
\begin{equation}\label{eq5dw:Jtop}
 \mathfrak J(0)=g\,m(0).
\end{equation}
Moreover, on every bounded admissible family there are functions
$b_1,b_2$ and remainders $f_1,f_2$ such that
\begin{equation}\label{eq5dw:slow-system}
 \begin{aligned}
 m'&=b_1\mathfrak J+f_1,\\
 \mathfrak J'&=\varepsilon b_2m+\varepsilon f_2,
 \end{aligned}
\end{equation}
where
\begin{equation}\label{eq5dw:slow-coeff}
 0<c\le b_1(p),b_2(p)\le C,
 \qquad |f_1(p)|+|f_2(p)|\le Ce^{\sigma p}.
\end{equation}
If $\alpha=b_1^{-1}$, then
\begin{equation}\label{eq5dw:alphaf1}
 |(\alpha f_1)'(p)|\le Ce^{\sigma p}.
\end{equation}
The constants are independent of $\varepsilon$.
\end{lemma}
\begin{proof}
Divide each layer equation by $A^3$.  In either layer a direct calculation gives
\[
 \partial_p\!\left(
   \frac{a^2}{2}-\frac{1+w_q^2}{2A^2}\right)
 +\partial_q\!\left(\frac{w_q}{A}\right)
 =\varepsilon a^{-3}\frac{w}{A^3}.
\]
Averaging in $q$ proves \eqref{eq5dw:Jprime}.  This identity is also
valid across the vorticity interface: $w$ and $w_p$ are continuous at
$p=p_0$, $a$ is continuous there, and hence no interface distribution is
created when the two layer identities are pasted together.  At $p=0$,
\[
 \frac{1+\overline w_q^2}{(\lambda^{-1/2}+\overline w_p)^2}
 =\lambda-2g\overline w,
 \qquad a(0;\lambda)^2=\lambda,
\]
so averaging yields \eqref{eq5dw:Jtop}.

On the lower layer write $\underline w_p=m'+\widetilde w_p$ and put
\[
 \overline A(p):=a^{-1}+m'=\langle A(\cdot,p)\rangle.
\]
Define
\[
 \Phi_p(r):=\frac{a^2}{2}-\frac1{2(a^{-1}+r)^2}
\]
and
\[
 R_1(p):=\frac1{2\overline A^2}
 -\left\langle
   \frac{1+\underline w_q^2}{2(\overline A+\widetilde w_p)^2}
  \right\rangle.
\]
Then $\mathfrak J=\Phi_p(m')+R_1$.  Since both
$a^{-1}$ and $\overline A$ are bounded above and below on a bounded
admissible family,
\[
 \Phi_p(m')=d_1(p)m',\qquad
 d_1(p)=\int_0^1(a^{-1}+\theta m')^{-3}\,\dd\theta,
 \qquad 0<c\le d_1\le C.
\]
Thus
\[
 m'=b_1\mathfrak J+f_1,
 \qquad b_1=d_1^{-1},\qquad f_1=-d_1^{-1}R_1.
\]
By Lemma~\ref{lem5.2}, Taylor expansion in $\widetilde w_p$ and
$\underline w_q$ gives
$|R_1|+|R_1'|\le Ce^{\sigma p}$.  The identity
\[
 \alpha f_1=-R_1,\qquad \alpha=b_1^{-1}=d_1,
\]
is important: it proves \eqref{eq5dw:alphaf1} without differentiating
$b_1$.

Finally, from \eqref{eq5dw:Jprime},
\[
 b_2:=a^{-3}\langle A^{-3}\rangle,
 \qquad
 f_2:=a^{-3}\langle\widetilde w A^{-3}\rangle.
\]
The admissibility margins give $0<c\le b_2\le C$, while
Poincar\'e's inequality and \eqref{eq5dw:osc-tail} give
$|f_2|\le Ce^{\sigma p}$.  This proves \eqref{eq5dw:slow-system}--
\eqref{eq5dw:alphaf1}.
\end{proof}

\begin{lemma}[Uniform half-line estimate at a closing gap]
\label{lem5dw:slow-resolvent}
Fix $P_0<p_0$, $M_0>0$, and $0<c_0<C_0$.  Let
$0\le\varepsilon\le\varepsilon_0$ and let $y$ be a stable solution of
\begin{equation}\label{eq5dw:abstract-slow}
 (\alpha y')'-\varepsilon b y=f\qquad\text{on }(-\infty,P_0],
\end{equation}
where
\[
 c_0\le\alpha(p),b(p)\le C_0,
 \qquad
 \int_{-\infty}^{P_0}(1+|p|)|f(p)|\,\dd p\le M_0.
\]
Assume that $\alpha$ and $b$ are continuous and that $f$ is locally
integrable; no modulus of continuity enters the constant below.
Here stable means $y,y'\to0$ for $\varepsilon>0$; for
$\varepsilon=0$ it means $y'\to0$ and $y=o(|p|)$ as
$p\to-\infty$.  Assume also
\[
 |y(P_0)|+|\alpha(P_0)y'(P_0)|\le M_0.
\]
Then, for $T\ge1$,
\begin{equation}\label{eq5dw:slow-resolvent-est}
 \sup_{p<P_0-T}\bigl(|y'(p)|+\varepsilon|y(p)|\bigr)
 \le C\left[
 \frac{M_0}{1+T}
 +\int_{-\infty}^{P_0-T/2}(1+|p|)|f(p)|\,\dd p\right],
\end{equation}
where $C$ is independent of $\varepsilon$.
\end{lemma}
\begin{proof}
Set
\[
 t=\int_{P_0}^{p}\alpha(s)^{-1}\,\dd s,
 \qquad Y(t)=y(p(t)),\qquad k=\sqrt\varepsilon.
\]
Then $t\le0$, $|t|\asymp|P_0-p|$, and
\[
 Y''-k^2c(t)Y=F(t),
 \qquad c(t)=\alpha(p(t))b(p(t)),
 \qquad F(t)=\alpha(p(t))f(p(t)),
\]
with $0<c_-\le c(t)\le c_+$.  Let $u_k$ be the positive solution of
the homogeneous equation that is stable at $-\infty$ and normalized by
$u_k(0)=1$.  Put
\[
 a_-:=k\sqrt{c_-},\qquad a_+:=k\sqrt{c_+}.
\]
We first construct this solution with constants uniform at the
threshold.  For $N>1$, solve on $[-N,0]$ with the initial Robin
condition
\[
 u_{k,N}'(-N)=a_-u_{k,N}(-N)>0
\]
and rescale so that $u_{k,N}(0)=1$.  Positivity and convexity follow
from $u_{k,N}''=k^2cu_{k,N}$.  The identities
\[
 \bigl(u_k'-a_-u_k\bigr)'+a_-
 \bigl(u_k'-a_-u_k\bigr)
   =k^2\bigl(c(t)-c_-\bigr)u_k\ge0,
\]
and
\[
 \bigl(u_k'-a_+u_k\bigr)'+a_+
 \bigl(u_k'-a_+u_k\bigr)
   =k^2\bigl(c(t)-c_+\bigr)u_k\le0
\]
applied first on $[-N,0]$ show that
\[
 a_-u_{k,N}(t)\le u_{k,N}'(t)
 \le a_+u_{k,N}(t).
\]
The normalized functions are locally equibounded and equicontinuous.
Letting $N\to\infty$ gives a positive solution $u_k$ on the
half-line satisfying
\[
 k\sqrt{c_-}\le r_k(t)\le k\sqrt{c_+}.
\]
Here $r_k=u_k'/u_k$.  In particular,
$u_k(t)\le e^{a_-t}$ for $t\le0$, so $u_k,u_k'\to0$ at
$-\infty$.  The Wronskian of two normalized stable solutions is
constant and vanishes at $-\infty$, which proves uniqueness.  This
also proves existence without invoking a $k$-dependent exponential
dichotomy theorem.

Therefore, for $t\le s\le0$,
\[
 \frac{u_k(t)}{u_k(s)}
 =\exp\!\left(-\int_t^s r_k(\tau)\,\dd\tau\right)
 \le e^{-\sqrt{c_-}\,k(s-t)},
 \qquad
 |u_k'(t)|\le \sqrt{c_+}\,k u_k(t).
\]
Thus the estimates below hold with
$\gamma=\sqrt{c_-}$ and a constant independent of
$0<k\le\sqrt{\varepsilon_0}$.  Moreover, on every fixed compact
interval, $r_k\to0$ and hence $u_k\to1$, $u_k'\to0$ as $k\to0$.
For $k=0$ set $u_0\equiv1$.

With
\[
 v_k(t)=u_k(t)\int_0^t u_k(\tau)^{-2}\,\dd\tau
\]
we have $v_k(0)=0$ and
$u_kv_k'-u_k'v_k=1$.  The stable Green kernel having zero value at
$t=0$ is therefore
\[
 G_k(t,s)=
 \begin{cases}
 u_k(t)v_k(s),&t\le s,\\
 v_k(t)u_k(s),&s\le t.
 \end{cases}
\]
For $t\le s$, monotonicity and the logarithmic-derivative bounds give
\[
 |v_k(s)|
 \le |s|u_k(s)^{-1},
\]
and hence
\begin{equation}\label{eq:green-separated-left}
 \begin{aligned}
 |G_k(t,s)|&\le |s|e^{-\gamma k(s-t)},\\
 |\partial_tG_k(t,s)|
 &\le Ck|s|e^{-\gamma k(s-t)}.
 \end{aligned}
\end{equation}
For $s\le t$, the same bounds yield
\[
 |G_k(t,s)|
 \le C\min\{|t|,k^{-1}\}e^{-\gamma k(t-s)}.
\]
Furthermore $v_k<0$, $v_k''=k^2cv_k<0$, and $v_k'(0)=1$, so
$v_k'>0$ on $(-\infty,0]$.  Using
\[
 v_k'(t)u_k(s)
 =\frac{u_k(s)}{u_k(t)}\left[
 1-r_k(t)\int_t^0
 \left(\frac{u_k(t)}{u_k(\tau)}\right)^2\dd\tau
 \right]
\]
and $0\le r_k(t)\int_t^0(u_k(t)/u_k(\tau))^2\dd\tau\le C$, we obtain
\begin{equation}\label{eq:green-separated-right}
 |\partial_tG_k(t,s)|
 \le Ce^{-\gamma k(t-s)},\qquad s\le t.
\end{equation}
Equations \eqref{eq:green-separated-left}--
\eqref{eq:green-separated-right} imply the global estimates
\[
 |G_k(t,s)|+|\partial_tG_k(t,s)|\le C(1+|s|),
\]
and, when $t\le s$,
\[
 |\partial_tG_k(t,s)|\le Ck|s|e^{-\gamma k(s-t)},\qquad
 k^2|G_k(t,s)|\le Ck^2|s|e^{-\gamma k(s-t)}.
\]
At $k=0$, $u_0=1$, $v_0(t)=t$, and
\[
 G_0(t,s)=\begin{cases}s,&t\le s,\\t,&s\le t,\end{cases}
\]
so all estimates pass continuously to the threshold; the stable
condition $Y'(-\infty)=0$ selects this kernel.

Variation of constants now gives
\[
 Y(t)=Y(0)u_k(t)+\int_{-\infty}^0G_k(t,s)F(s)\,\dd s
\]
and the global kernel estimate gives $\|Y\|_\infty\le CM_0$
uniformly in $k$.  For $t<-T$, split the integral into
$s<-T/2$ and $s\ge-T/2$.  The derivative and $k^2$ times the value of
the deep part are bounded by the weighted tail of $F$.  On the
shallow part we have $t\le s$ and use
\eqref{eq:green-separated-left} together with
\[
 \sup_{k\ge0}ke^{-\gamma kT/2}\le\frac{C}{1+T},\qquad
 \sup_{k\ge0}k^2e^{-\gamma kT/2}\le\frac{C}{(1+T)^2}.
\]
The homogeneous term $Y(0)u_k(t)$ satisfies the same bounds because
$|u_k'(t)|\le Cke^{-\gamma k|t|}$ and
$u_k(t)\le e^{-\gamma k|t|}$.  Consequently
\[
 \sup_{t<-T}\bigl(|Y'(t)|+k^2|Y(t)|\bigr)
 \le C\left[
  \frac{M_0}{1+T}
  +\int_{-\infty}^{-T/2}(1+|s|)|F(s)|\,\dd s
 \right].
\]
Since $dt/dp=\alpha^{-1}$ and $|t|\asymp|P_0-p|$, this is precisely
\eqref{eq5dw:slow-resolvent-est}, after changing $C$ by a quantity
depending only on $P_0,c_0,C_0$, and $\varepsilon_0$.
\end{proof}

\begin{lemma}[Uniform physical tail for the slow mode]
\label{lem5dw:mean-tail}
Fix $\delta>0$ and $M>0$.  There is a function
$\omega_{\delta,M}(P)\downarrow0$ as $P\to\infty$, independent of
$0\le\varepsilon\le\varepsilon_0$, such that every oriented exact or
regularized solution satisfying the admissibility margins and
\[
 |\lambda|+\|\overline w\|_{C^{3,\alpha}(\overline R_1)}
 +\|\underline w_p\|_{C^{2,\alpha}(\overline R_2)}
 +\|\underline w_q\|_{C^{2,\alpha}(\overline R_2)}\le M
\]
obeys
\begin{equation}\label{eq5dw:mean-tail}
 \sup_{p<-P}\bigl(|m'(p)|+|\mathfrak J(p)|
 +\varepsilon|m(p)|\bigr)\le\omega_{\delta,M}(P).
\end{equation}
Moreover,
\begin{equation}\label{eq5dw:derivative-tail}
 \sup_{s<-P}\left(
 \|\underline w_p\|_{C^{2,\alpha}(S_s)}
 +\|\underline w_q\|_{C^{2,\alpha}(S_s)}
 +\varepsilon\|\underline w\|_{C^{1,\alpha}(S_s)}
 \right)\le\widetilde\omega_{\delta,M}(P),
\end{equation}
where $\widetilde\omega_{\delta,M}(P)\to0$.
\end{lemma}
\begin{proof}
Fix $P_0<p_0-2$.  The trace at $p=p_0$, the derivative bounds and
\eqref{eq5dw:Jdef} give uniform bounds for $m(P_0)$,
$m'(P_0)$ and $\mathfrak J(P_0)$.  By the first equation of
\eqref{eq5dw:slow-system},
\[
 \mathfrak J=\alpha(m'-f_1),\qquad \alpha=b_1^{-1}.
\]
Substitution into the second equation gives
\begin{equation}\label{eq5dw:mean-scalar}
 (\alpha m')'-\varepsilon b_2m
 =(\alpha f_1)'+\varepsilon f_2.
\end{equation}
The coefficients $\alpha,b_2$ are uniformly positive and bounded.
More explicitly, if
\[
 f_\varepsilon:=(\alpha f_1)'+\varepsilon f_2,
\]
then \eqref{eq5dw:slow-coeff}--\eqref{eq5dw:alphaf1} give
\begin{equation}\label{eq:mean-forcing-first-moment}
 |f_\varepsilon(p)|\le Ce^{\sigma p},\qquad
 \sup_{0\le\varepsilon\le\varepsilon_0}
 \int_{-\infty}^{P_0}(1+|p|)|f_\varepsilon(p)|\,\dd p\le C,
\end{equation}
and, moreover,
\[
 \sup_{0\le\varepsilon\le\varepsilon_0}
 \int_{-\infty}^{-P}(1+|p|)|f_\varepsilon(p)|\,\dd p
 \longrightarrow0\qquad(P\to\infty).
\]
All constants depend only on the fixed data, $\delta$, and $M$.
Thus the first-moment hypothesis in
Lemma~\ref{lem5dw:slow-resolvent} is verified uniformly over the
entire bounded admissible family, rather than separately for each
solution.  For $\varepsilon>0$ the original
regularized space supplies $m,m'\to0$ at $-\infty$.  For
$\varepsilon=0$, the deep-water derivative condition gives
$m'\to0$, and reconstruction from a fixed finite level gives
$m(p)=o(|p|)$.  Lemma~\ref{lem5dw:slow-resolvent} applied to
\eqref{eq5dw:mean-scalar} proves the bounds for $m'$ and
$\varepsilon m$.  The formula
$\mathfrak J=\alpha(m'-f_1)$ gives the same tail for $\mathfrak J$.

Lemma~\ref{lem5.2} controls the oscillatory part, hence
$\underline w_p=m'+\widetilde w_p\to0$ uniformly on the lower tail.
We next justify the endpoint H\"older norms in
\eqref{eq5dw:derivative-tail}.  On each enlarged unit interval
$I_s^*=[s-2,s+2]$, the scalar equation
\eqref{eq5dw:mean-scalar}, the uniform $C^{1,\alpha}$ coefficient
bounds, and the already established supremum estimate give the
one-dimensional interior estimate
\begin{equation}\label{eq:mean-endpoint-strip}
 \|m'\|_{C^{1,\alpha}(I_s)}
 +\varepsilon\|m\|_{C^{1,\alpha}(I_s)}
 \le C\left(
  \sup_{I_s^*}(|m'|+\varepsilon|m|)
  +\|f_\varepsilon\|_{C^\alpha(I_s^*)}
 \right).
\end{equation}
Here $I_s=[s-1,s+1]$.  The construction of $f_1,f_2$ from the
oscillatory part and Lemma~\ref{lem5.2} yield
$\|f_\varepsilon\|_{C^\alpha(I_s^*)}\le Ce^{\sigma s}$.
Consequently the right-hand side of
\eqref{eq:mean-endpoint-strip} tends to zero uniformly in
$\varepsilon$ as $s\to-\infty$.  This supplies, in particular, the
local $C^\alpha$ smallness of $\underline w_p$ and of
$\varepsilon\underline w$ that is needed below; it is not obtained
from an endpoint compact embedding.

Solve the lower bulk equation algebraically for $\underline w_{pp}$:
\[
 \underline w_{pp}=\frac{
  2A\underline w_q\underline w_{pq}-A^2\underline w_{qq}
  -\gamma(-p)A^3+\gamma(-p)a^{-3}(1+\underline w_q^2)
  +\varepsilon a^{-3}\underline w}{1+\underline w_q^2}.
\]
Every term on the right tends to zero uniformly.  The only
undifferentiated term is controlled by $\varepsilon|m|$ and
\eqref{eq5dw:osc-tail}.  The same conclusion holds in local
$C^\alpha$ norms by \eqref{eq:mean-endpoint-strip}.  The only point which needs care before differentiating this
identity in $p$ is the vorticity term.  Write
\[
 A=a^{-1}+\underline w_p,\qquad B=1+\underline w_q^2,
\]
and set
\[
 \mathcal V:=\gamma(-p)\bigl(a^{-3}B-A^3\bigr).
\]
Then
\[
 \begin{aligned}
 \mathcal V_p
 ={}&-\gamma'(-p)\bigl(a^{-3}B-A^3\bigr)\\
 &+\gamma(-p)\Bigl(
 (a^{-3})_p\underline w_q^2
 +2a^{-3}\underline w_q\underline w_{pq}
 -3(A^2-a^{-2})(a^{-1})_p
 -3A^2\underline w_{pp}
 \Bigr).
 \end{aligned}
\]
Here we used
\[
 (a^{-3})_p=3a^{-2}(a^{-1})_p.
\]
Moreover,
\[
 a^{-3}B-A^3
 =a^{-3}\underline w_q^2-(A^3-a^{-3}),
\]
so the coefficient of $\gamma'(-p)$ also contains only
$\underline w_q$ and $\underline w_p$.  Thus every occurrence of
$\gamma$ or $\gamma'$ in the $p$-differentiated formula is multiplied
by one of the already small quantities
\[
 \underline w_q,\quad \underline w_{pq},\quad
 \underline w_p,\quad \underline w_{pp}.
\]
This is the precise laminar cancellation: all these terms vanish when
$D\underline w=0$.  The globally bounded piecewise
$C^{1,\alpha}$ norm of $\gamma$ is therefore sufficient; no unproved
decay of its endpoint H\"older seminorm is being used.

Finally, unit-strip Schauder estimates for the homogeneous
differentiated $q$-equation give
\[
 \|\underline w_q\|_{C^{2,\alpha}(S_s)}
 \le C\|\underline w_q\|_{C^0(S_s^*)}\longrightarrow0.
\]
The algebraic formula then gives
$\|\underline w_{pp}\|_{C^\alpha(S_s)}\to0$; differentiating it in
$q$ controls the mixed third derivatives, while differentiating it
in $p$ expresses $\underline w_{ppp}$ in terms of the already
controlled derivatives and the laminar-cancelled vorticity terms.
Together with \eqref{eq:mean-endpoint-strip} this proves the full
endpoint estimate \eqref{eq5dw:derivative-tail}.
\end{proof}

\begin{remark}[Why the absolute zero-mode tail is not used]
\label{rem5dw:no-absolute-tail}
We do not assert
\[
 \sup_{0<\varepsilon\le\varepsilon_0}\sup_{p<-P}|m(p)|\to0.
\]
The scalar model $z''-\varepsilon z=0$ already has a stable mode
$e^{\sqrt\varepsilon p}$ whose decay length diverges as
$\varepsilon\downarrow0$.  What is uniform, and what the Euler
condition actually requires, is the decay of $m'$ together with
$\varepsilon m$.  This is the reason for the topology introduced next.
\end{remark}

\subsection{The natural deep-water topology and singular compactness}

Let $X_{\rm dw}$ be the completion of $X$ for the norm
\begin{equation}\label{eq5dw:Xnorm}
 \|(\overline w,\underline w)\|_{X_{\rm dw}}
 :=\|\overline w\|_{C^{3,\alpha}(\overline R_1)}
 +\|\underline w_p\|_{C^{2,\alpha}(\overline R_2)}
 +\|\underline w_q\|_{C^{2,\alpha}(\overline R_2)}.
\end{equation}
For clarity, this completion has the following concrete
characterization.  It is the closed subspace of triples
$(\overline w,U,V)$ in
\[
 C^{3,\alpha}_{\rm per}(\overline R_1)
 \times C^{2,\alpha}_{\rm per,0}(\overline R_2)
 \times C^{2,\alpha}_{\rm per,0}(\overline R_2)
\]
for which
\begin{equation}\label{eq5dw:completion-compatibility}
 U_q=V_p\quad\hbox{in }R_2,\qquad
 U(\cdot,p_0)=\overline w_p(\cdot,p_0),\qquad
 V(\cdot,p_0)=\overline w_q(\cdot,p_0),
\end{equation}
with the even/odd parity inherited from
$U=\underline w_p$ and $V=\underline w_q$.  These relations are
closed in the product norm, so this description is a Banach space.
To verify the converse identification explicitly, reconstruct
$\underline w$ by \eqref{eq5dw:reconstruct} and choose a smooth
cutoff $\chi_N(p)$ which is one for $p\ge-N$, zero for $p\le-2N$,
and satisfies $|\chi_N^{(j)}|\le C_jN^{-j}$.  Then
$\underline w_N:=\chi_N\underline w$ has the required endpoint
decay and retains both interface traces for large $N$.  Because
$U,V$ tend to zero in their full uniform local $C^{2,\alpha}$ norms
and $\underline w=o(|p|)$, the product rule gives
\[
 \|\partial_p\underline w_N-U\|_{C^{2,\alpha}}
 +\|\partial_q\underline w_N-V\|_{C^{2,\alpha}}
 \longrightarrow0.
\]
Thus $(\overline w,\underline w_N)\in X$ approximates the triple in
the norm \eqref{eq5dw:Xnorm}.  This proves that the concrete
description agrees with the abstract completion.
The lower component is reconstructed on every finite subcylinder from
the interface trace and $\underline w_p$ by
\begin{equation}\label{eq5dw:reconstruct}
 \underline w(q,p)=\overline w(q,p_0)
 -\int_p^{p_0}\underline w_p(q,s)\,\dd s.
\end{equation}
Under the above identification, $U=\underline w_p$ and
$V=\underline w_q$.  The compatibility and trace relations in
\eqref{eq5dw:completion-compatibility} ensure that the $q$-derivative
of \eqref{eq5dw:reconstruct} is exactly $\underline w_q$.  Since
$\underline w_p\to0$, every reconstructed
lower component satisfies $\underline w(q,p)=o(|p|)$.

For the statement of the physical theorem we also use the strong
space $X_{\rm s}$ with norm
\begin{equation}\label{eq5dw:Xstrong}
 \|(\overline w,\underline w)\|_{X_{\rm s}}
 :=\|\overline w\|_{C^{3,\alpha}(\overline R_1)}
 +\|\underline w\|_{C^{0,\alpha}(\overline R_2)}
 +\|\underline w_p\|_{C^{2,\alpha}(\overline R_2)}
 +\|\underline w_q\|_{C^{2,\alpha}(\overline R_2)}.
\end{equation}
Thus $X_{\rm s}$ records the absolute lower zero mode, whereas
$X_{\rm dw}$ is the topology in which the singular limit is first
compact.  The exact zero-mode lemma below recovers
\eqref{eq5dw:Xstrong} on the limiting continuum.

Define $\mathcal O_\delta^{\rm dw}\subset\mathbb R\times X_{\rm dw}$
by the same strict parameter, ellipticity and top-obliqueness inequalities
as $\mathcal O_\delta$.  Explicitly,
\[
\begin{split}
 \mathcal O_\delta^{\rm dw}:=\{(\lambda,\overline w,\underline w):\;&
 \lambda>-2\Gamma_{\inf}+\delta,\\
 &\inf_{\overline R_1}(a^{-1}+\overline w_p)>\delta,\quad
 \inf_{\overline R_2}(a^{-1}+\underline w_p)>\delta,\\
 &\sup_{q\in\mathbb S}\left(\overline w(q,0)
 -\frac{2\lambda-\delta}{4g}\right)<0\}.
\end{split}
\]
It is open in $\mathbb R\times X_{\rm dw}$.  The exact map $F$ is
continuous on this set: the lower bulk operator contains
$\underline w$ only through its derivatives, while the interface and
surface terms are controlled by the upper component and the traces in
\eqref{eq5dw:Xnorm}.

\begin{proposition}[Singular compactness in $X_{\rm dw}$]
\label{lem5.5-uniform}
Fix $\delta>0$ and $M>0$.  Let $\varepsilon_j\downarrow0$ and let
$z_j=(\lambda_j,\overline w_j,\underline w_j)$ be oriented regularized
solutions satisfying the $\mathcal O_\delta^{\rm dw}$ margins and
\[
 |\lambda_j|+\|(\overline w_j,\underline w_j)\|_{X_{\rm dw}}\le M.
\]
Then, after passage to a subsequence,
\[
 z_j\longrightarrow z=(\lambda,\overline w,\underline w)
 \quad\text{in }\mathbb R\times X_{\rm dw},
\]
where $F(z)=0$.  The same conclusion holds for a bounded sequence of
exact oriented solutions.  In particular, bounded subsets of the
oriented exact solution set are relatively compact in
$\mathbb R\times X_{\rm dw}$.
\end{proposition}
\begin{proof}
Passing to a subsequence, $\lambda_j\to\lambda$.  On every finite
transmission cylinder, \eqref{eq5dw:reconstruct} and the
$X_{\rm dw}$ bound give a uniform $C^0$ bound.  Local transmission,
interior and oblique Schauder estimates therefore yield, diagonally,
convergence in piecewise $C^{3,\beta}$ for every $0<\beta<\alpha$.
On each fixed core $\varepsilon_j\underline w_j\to0$ and
$\varepsilon_j\overline w_j\to0$ in the target norms, so the local
limit solves the exact equation.

We now recover the endpoint H\"older exponent without using the false
compact embedding of a bounded $C^{3,\alpha}$ set into itself.  Fix
$P<p_0-3$, put
\[
 C_P:=\overline R_1\cup
  \bigl(\mathbb S\times[P,p_0]\bigr),
\]
and let $C_P^*$ be the core obtained by moving the artificial lower
edge down by two units.  The uniform $X_{\rm dw}$ bound and
\eqref{eq5dw:reconstruct} imply a uniform piecewise
$C^{3,\alpha}(C_P^*)$ bound.  In particular, the local limit belongs
to piecewise $C^{3,\alpha}$ with the same bound.

Write $w_j=(\overline w_j,\underline w_j)$,
$w=(\overline w,\underline w)$, and $U_j=w_j-w$.  Define
\[
 \mathscr A_j:=\int_0^1D_wF
 \bigl(\lambda_j,w+t(w_j-w)\bigr)\,\dd t.
\]
On $C_P^*$ the exact mean-value identity is
\begin{equation}\label{eq:endpoint-difference-equation}
 \mathscr A_jU_j=R_j,
 \qquad
 R_j=\varepsilon_j\mathcal M(\lambda_j)w_j
 -\bigl(F(\lambda_j,w)-F(\lambda,w)\bigr),
\end{equation}
where
\[
 \mathcal M(\lambda)(\overline v,\underline v)
 :=(a^{-3}\overline v,a^{-3}\underline v,0).
\]
For a sequence of exact solutions the first term is absent.  Smooth
dependence on $\lambda$, the core $C^{3,\alpha}$ bounds, and
$\varepsilon_j\to0$ give convergence in the actual endpoint target:
\begin{equation}\label{eq:endpoint-residual-convergence}
 \|R_j\|_{C^{1,\alpha}(C_P^*)\times
 C^{2,\alpha}(\mathbb S)}
 \le C_P\bigl(\varepsilon_j+|\lambda_j-\lambda|\bigr)
 \longrightarrow0.
\end{equation}

The coefficients of $\mathscr A_j$ are uniformly bounded in
piecewise $C^{1,\alpha}$.  Convexity of the principal matrices and
the admissibility inequalities give a common ellipticity constant;
the coefficient of the top normal derivative stays uniformly
negative, and the perfect-transmission complementing constant is
uniform.  The local Agmon--Douglis--Nirenberg/transmission estimate
on the nested cores therefore yields
\[
 \|U_j\|_{C^{3,\alpha}(C_P)}
 \le C_P\left(
  \|R_j\|_{C^{1,\alpha}(C_P^*)\times C^{2,\alpha}}
  +\|U_j\|_{C^0(C_P^*)}
 \right)\longrightarrow0.
\]
Here the last term tends to zero by the previously obtained
$C^{3,\beta}_{\rm loc}$ convergence.  Thus endpoint convergence holds
on every fixed core.

It remains to control the cylindrical tail.  Lemmas~\ref{lem5.2} and
\ref{lem5dw:mean-tail} imply that, for every $\eta>0$, a level
$P<p_0$ can be chosen independently of $j$ so that
\[
 \sup_{s<P}\left(
 \|(\underline w_j)_p\|_{C^{2,\alpha}(S_s)}
 +\|(\underline w_j)_q\|_{C^{2,\alpha}(S_s)}\right)<\eta,
\]
and the same holds for the limit.  The unit-strip estimates control
the global H\"older seminorm on the whole tail as follows.  Pairs of
points at distance less than one lie in a common enlarged unit strip,
whereas pairs at distance at least one are controlled by twice the
supremum norm.  Hence the tail $X_{\rm dw}$ norm of the difference is
at most $C\eta$.  First choose $P$ so that this is small and then use
the endpoint core convergence above.  This proves
$z_j\to z$ in $\mathbb R\times X_{\rm dw}$.

Finally, \eqref{eq5dw:derivative-tail} gives, on the lower end,
\[
 \sup_{s<P}\varepsilon_j
 \|\underline w_j\|_{C^{1,\alpha}(S_s)}
 \le\widetilde\omega_{\delta,M}(|P|).
\]
On the fixed core
$\varepsilon_j\underline w_j\to0$ in $C^{1,\alpha}$ by
$\varepsilon_j\to0$ and the core bound; the upper regularization term
is treated directly on the bounded upper layer.  Thus the
regularization residual vanishes in the global target norm, and the
limit satisfies $F(z)=0$.
\end{proof}

\begin{lemma}[Recovery of the exact zero mode and the strong topology]
\label{lem5dw:exact-strong-recovery}
Fix $\delta>0$ and $M>0$.  Every oriented exact solution satisfying
the $\mathcal O_\delta^{\rm dw}$ margins and
\[
 |\lambda|+\|(\overline w,\underline w)\|_{X_{\rm dw}}\le M
\]
has
\begin{equation}\label{eq:exact-flux-zero}
 \mathfrak J\equiv0,\qquad m(0)=0,
\end{equation}
and there is a number $m_-=m_-(w)$ such that
\begin{equation}\label{eq:exact-mean-absolute-tail}
 |m_-|+\sup_{p<P_0}e^{-\sigma p}|m(p)-m_-|
 +\|(\overline w,\underline w)\|_{X_{\rm s}}
 \le C(\delta,M).
\end{equation}
Moreover, on such bounded exact families, convergence in
$\mathbb R\times X_{\rm dw}$ implies convergence in
$\mathbb R\times X_{\rm s}$.
\end{lemma}
\begin{proof}
For an exact solution, \eqref{eq5dw:Jprime} gives
$\mathfrak J'=0$.  On the lower end the first slow equation reads
\[
 \mathfrak J=\alpha(m'-f_1).
\]
The deep-water condition and \eqref{eq5dw:slow-coeff} give
$m'(p),f_1(p)\to0$ as $p\to-\infty$.  Hence the constant
$\mathfrak J$ is zero.  The top identity
\eqref{eq5dw:Jtop} then gives $m(0)=0$, proving
\eqref{eq:exact-flux-zero}.  On the lower layer we now have simply
\begin{equation}\label{eq:exact-mean-prime}
 m'=f_1,\qquad |f_1(p)|\le Ce^{\sigma p}.
\end{equation}
Consequently $m$ has a finite limit $m_-$ at the infinite bottom and
\[
 |m(p)-m_-|\le Ce^{\sigma p}.
\]
The interface trace is controlled by the upper $C^{3,\alpha}$ norm,
so integration of \eqref{eq:exact-mean-prime} also bounds $m_-$.
Together with Lemma~\ref{lem5.2} and
Lemma~\ref{lem5dw:mean-tail}, this proves
\eqref{eq:exact-mean-absolute-tail}.

It remains to verify the asserted topology, which is the point needed
in Section~6.  Suppose $z_j\to z$ in
$\mathbb R\times X_{\rm dw}$ in a family with the displayed bounds.
The oscillatory parts converge uniformly by Poincar\'e's inequality
and the convergence of $(w_j)_q$.  On every finite lower core,
$f_{1,j}\to f_1$ uniformly by the formula in
Lemma~\ref{lem5.3-tail}.  On the remaining tail,
\eqref{eq:exact-mean-prime} supplies the common integrable majorant
$Ce^{\sigma p}$.  Splitting the integral at a fixed deep level and
then letting that level tend to $-\infty$ gives
\[
 \sup_{p\le p_0}|m_j(p)-m(p)|\longrightarrow0.
\]
Thus $\underline w_j\to\underline w$ uniformly.  Since the first
derivatives already converge in $C^{2,\alpha}$ by the
$X_{\rm dw}$ convergence, the interpolation inequality
\[
 [u]_{\alpha}\le
 C\|u\|_{C^0}^{1-\alpha}\|Du\|_{C^0}^{\alpha}
\]
upgrades the uniform convergence to $C^{0,\alpha}$.  All remaining
terms of the $X_{\rm s}$ norm are already contained in the
$X_{\rm dw}$ norm.  This proves the last assertion.
\end{proof}

\begin{lemma}[Recovery of the auxiliary $X$-bound for fixed regularization]
\label{lem5dw:fixed-eps-X}
Fix $\varepsilon>0$, $\delta>0$ and $M>0$.  Every regularized solution
satisfying the admissibility margins and
\[
 |\lambda|+\|(\overline w,\underline w)\|_{X_{\rm dw}}\le M
\]
obeys
\[
 \|(\overline w,\underline w)\|_X\le C(\varepsilon,\delta,M).
\]
\end{lemma}
\begin{proof}
Let $m=\langle\underline w\rangle$.  By
Lemma~\ref{lem5dw:mean-tail}, for fixed $\varepsilon>0$ there is a
level $P<p_0$ such that
\[
 \sup_{p<P}\varepsilon |m(p)|\le C(\delta,M).
\]
Hence
\[
 \sup_{p<P}|m(p)|\le C(\varepsilon,\delta,M).
\]
On the finite strip $P\le p\le p_0$, the interface trace together
with the $X_{\rm dw}$ bound on $\underline w_p$ gives the same bound
for $m$.  Lemma~\ref{lem5.2} controls the oscillatory component
$\widetilde w=\underline w-m$, and therefore
\[
 \|\underline w\|_{C^0(\overline R_2)}
 \le C(\varepsilon,\delta,M).
\]
The $X_{\rm dw}$ norm already controls all lower-layer derivatives of
orders one through three in $C^{0,\alpha}$, while the upper component
is controlled in $C^{3,\alpha}$.  Since the regularized solution
belongs to $X$, the preceding $C^0$ estimate therefore yields
\[
 \|(\overline w,\underline w)\|_X
 \le C(\varepsilon,\delta,M).
\]
\end{proof}

\begin{lemma}[Exact $q$-independent states in $X_{\rm dw}$]
\label{lem5dw:qind-exact}
If an exact solution in $\mathcal O_\delta^{\rm dw}$ is independent of
$q$, then it is the trivial laminar state.
\end{lemma}
\begin{proof}
As in Lemma~\ref{lem:qind-exact}, the bulk equation gives
$h_p^{-2}=a^2$, hence $w_p=0$ and $w\equiv C$.  The top boundary
condition then gives
\[
 1+(2gC-\lambda)\lambda^{-1}=0,
\]
so $C=0$.
\end{proof}

\begin{lemma}[Uniform exclusion of bounded regularized shear states]
\label{lem5dw:qind-approx}
Fix $\delta>0$ and $M>0$.  There exists $\varepsilon_*>0$ such that,
for $0<\varepsilon<\varepsilon_*$, there is no nonzero
$q$-independent regularized solution satisfying the
$\mathcal O_\delta^{\rm dw}$ margins and
\[
 |\lambda|+\|(\overline w,\underline w)\|_{X_{\rm dw}}\le M.
\]
\end{lemma}
\begin{proof}
The proof of Lemma~\ref{lem:qind-approx} uses only the parameter bound,
the lower derivative bound, the top trace and the positive admissibility
margins, all of which are controlled by $X_{\rm dw}$.  For completeness,
suppose $\varepsilon_j\downarrow0$ and normalize a nonzero shear state
by $u_j=W_j/W_j(0)$.  The one-dimensional maximum principle gives
$0<u_j\le1$ and $u_j'\ge0$.  The coefficients in
\[
 u_j''+b_j(p)u_j'-\varepsilon_j a(p;\lambda_j)^{-3}u_j=0
\]
are locally precompact and $b_j$ is uniformly integrable at the lower
end.  Thus $u_j\to u$ locally, where
$u''+bu'=0$, $0\le u\le1$, $u'\ge0$ and $u(0)=1$.  Writing
$M=e^{\int_0^pb}$ gives $Mu'=C$; boundedness and nonnegativity on the
whole half-line force $C=0$.  Hence $u'(0)=0$.  On the other hand the
top boundary condition gives, by the mean-value theorem,
$u_j'(0)\ge c(\delta,M)>0$, a contradiction.
\end{proof}

\begin{lemma}[Uniform tail for the normalized transverse derivative]
\label{lem5.4-tail}
Let $(\lambda_j,\overline w_j,\underline w_j)$ be oriented exact
solutions bounded in $\mathbb R\times X_{\rm dw}$ and satisfying fixed
admissibility margins.  Whenever $(w_j)_q\not\equiv0$, put
\[
 v_j:=\frac{(w_j)_q}
 {\|(w_j)_q\|_{C^{2,\alpha}(\overline R_1)
 \times C^{2,\alpha}(\overline R_2)}}.
\]
Then there are $C,\sigma>0$, independent of $j$, such that
\[
 \sum_{i+j'\le2}|\partial_p^i\partial_q^{j'}v_j(q,p)|
 \le Ce^{\sigma p},\qquad p\le p_0-2.
\]
Thus the normalized mass cannot escape to the infinite bottom.
\end{lemma}
\begin{proof}
The normalized functions solve the homogeneous equation obtained by
differentiating the exact problem in $q$; they are odd in $q$ and
satisfy the corresponding homogeneous transmission and top oblique
conditions.  Because the family is bounded in
$\mathbb R\times X_{\rm dw}$ and has fixed admissibility margins, all
coefficients have uniform piecewise $C^{1,\alpha}$ bounds, while the
ellipticity, transmission-complementing, and top-obliqueness constants
have common positive lower bounds.

We first convert the $C^{2,\alpha}$ normalization into nonvanishing
$C^0$ mass.  On a fixed bounded core, the transmission/oblique
Schauder estimate gives
\[
 \|v_j\|_{C^{2,\alpha}(\text{core})}
 \le C\|v_j\|_{C^0(\text{enlarged core})},
\]
because both the interior and boundary right-hand sides vanish.  On
the lower cylinder, translation of the standard unit-strip interior
estimate gives
\[
 \|v_j\|_{C^{2,\alpha}(S_s)}
 \le C\|v_j\|_{C^0(S_s^*)},
 \qquad s<p_0-2,
\]
with the same constant $C$ for all $j$ and all translation levels.
Combining the core and tail estimates yields
\begin{equation}\label{eq:normalized-v-schauder}
 \|v_j\|_{C^{2,\alpha}(\overline R_1)
 \times C^{2,\alpha}(\overline R_2)}
 \le C\|v_j\|_{C^0(\overline R_1\cup\overline R_2)}.
\end{equation}
Since the left-hand side equals one,
\[
 \|v_j\|_{C^0(\overline R_1\cup\overline R_2)}\ge C^{-1}.
\]
Thus the normalization cannot be carried only by derivatives or by a
H\"older seminorm.

If a definite part of this $C^0$ mass could escape to the infinite
bottom, there would be points $(q_j,p_j)$ with $p_j\to-\infty$ and
$|v_j(q_j,p_j)|\ge c_0>0$.  Translating the lower equations by $p_j$
and using the endpoint coefficient convergence gives, after a
subsequence,
\[
 v_j(q,p+p_j)\longrightarrow V(q,p)
 \quad\text{locally in }C^{2,\beta}(\mathbb S\times\mathbb R),
 \qquad 0<\beta<\alpha,
\]
where $V$ is bounded, odd in $q$, nonzero, and solves
\[
 V_{pp}+a_\infty^{-2}V_{qq}=0
 \quad\text{on }\mathbb S\times\mathbb R.
\]
Writing $V(q,p)=\sum_{k\ge1}V_k(p)\sin(kq)$ gives
\[
 V_k''-k^2a_\infty^{-2}V_k=0.
\]
No nonzero solution of this equation is bounded on the whole real
line, so $V\equiv0$, a contradiction.  Hence a fixed finite core
captures a definite amount of the normalized $C^0$ mass.

Finally apply the Phragm\'en--Lindel\"of barrier used in
Lemma~\ref{lem5.2} to the homogeneous equation for $v_j$ on the lower
semi-cylinder.  The zero Fourier mode is absent because $v_j$ is odd,
and the limiting transverse gap is bounded below uniformly.  Since
$\|v_j\|_{C^0}\le1$ and the data on one fixed cross-section are
uniformly bounded, the barrier gives
\[
 |v_j(q,p)|\le Ce^{\sigma p},\qquad p\le p_0-2,
\]
with $C,\sigma$ independent of $j$.  Translation-invariant unit-strip
Schauder estimates then give the same exponential bound for all
derivatives of total order at most two.  This proves the lemma.
\end{proof}

\subsection{Whyburn continuation in $X_{\rm dw}$}

Let $\mathcal C_\delta^{\rm dw}\subset\mathcal O_\delta^{\rm dw}$ be
the set consisting of $(\lambda_*^0,0,0)$ and all nontrivial exact
solutions satisfying either of the two oriented nodal patterns
$\mathcal N_0^\pm$.

\begin{lemma}[Closedness at the trivial line]
\label{lem5.4}
The set $\mathcal C_\delta^{\rm dw}$ is relatively closed in
$\mathcal O_\delta^{\rm dw}$.  In particular, if nontrivial oriented
exact solutions converge in $\mathbb R\times X_{\rm dw}$ to
$(\lambda,0,0)$, then $\lambda=\lambda_*^0$.
\end{lemma}
\begin{proof}
A nontrivial $q$-dependent limit retains the strict nodal inequalities
by the maximum principle, the two-sided transmission Hopf lemma, the
surface oblique boundary lemma and Serrin's corner lemma.  If the limit
is $q$-independent, Lemma~\ref{lem5dw:qind-exact} makes it trivial.
It remains to identify the parameter of a trivial limit.  Normalize
$v_j$ as in Lemma~\ref{lem5.4-tail}.  That lemma prevents its mass from
escaping to $p=-\infty$, so a subsequence converges locally to a
nonzero odd solution of the exact linearized problem.  The first sine
coefficient has a fixed sign and corresponds to generalized eigenvalue
$-1$.  If the principal eigenvalue were strictly smaller than $-1$,
its positive eigenfunction would be weighted-orthogonal to this
fixed-sign coefficient, a contradiction.  Hence
$\mu^0(\lambda)=-1$, and Lemma~\ref{lem4.3} gives
$\lambda=\lambda_*^0$.  Proposition~\ref{lem5.5-uniform} supplies the
sequential compactness needed for relative closedness.
\end{proof}

\begin{lemma}[Whyburn continuation principle]
\label{thm5.1}
Let $S$ be a relatively closed subset of an open set in a metric
space, contain $z_0$, and have relatively compact bounded subsets.  If
every bounded open neighborhood $U$ of $z_0$ with closure contained in
the admissible open set satisfies $\partial U\cap S\ne\varnothing$,
then the maximal connected subset of $S$ containing $z_0$ is unbounded
or its closure meets the boundary of the admissible set.
\end{lemma}
\begin{proof}
This is the standard Whyburn separation argument; see
\cite[Theorem A6]{AT81} and, for the same deep-water continuation
mechanism, \cite{Hur11}.  If the component through $z_0$ were bounded
and had closure contained in the admissible open set, relative
compactness would make that closure compact.  Separation of compact
metric spaces would then produce a bounded open neighborhood $U$ of
$z_0$, with admissible closure, whose boundary misses $S$, a
contradiction.
\end{proof}

\begin{lemma}[Boundary-point property in the deep-water topology]
\label{lem5.6}
Let $U$ be a bounded open neighborhood of $(\lambda_*^0,0,0)$ in
$\mathbb R\times X_{\rm dw}$ with
$\overline U\subset\mathcal O_\delta^{\rm dw}$.  For all sufficiently
small $\varepsilon>0$, the positive oriented regularized route contains
a point on $\partial U$.
\end{lemma}
\begin{proof}
Choose $\delta'\in(0,\delta)$ and use the distinguished regularized
route constructed in $\mathcal O_{\delta'}$.  For small
$\varepsilon$ its local germ is the Crandall--Rabinowitz germ through
$(\lambda_*^\varepsilon,0,0)$.  Decrease $\varepsilon$, if necessary,
so that this point lies in $U$.  Since $U$ is bounded in
$\mathbb R\times X_{\rm dw}$ and $\overline U$ has the
$\mathcal O_\delta^{\rm dw}$ margins, there is an $M$ independent of
$\varepsilon$ such that every point of $\overline U$ satisfies the
bound in Lemma~\ref{lem5dw:qind-approx}.  That lemma therefore
excludes every nonzero $q$-independent regularized solution from
$\overline U$ for all sufficiently small $\varepsilon$.

Suppose that the Crandall--Rabinowitz positive half-route never leaves
$U$.  We first establish its orientation, before invoking any estimate
whose proof uses that orientation.  Lemma~\ref{lem4.7} puts an initial
punctured segment in $\mathcal N_\varepsilon^+$.  If this property
failed at a first finite parameter value, Lemma~\ref{lem4.8} would
make the limiting point $q$-independent.  A nonzero such point has
just been excluded.  If the point were trivial, Lemma~\ref{lem4.9}
would identify it as $(\lambda_*^\varepsilon,0,0)$.  The intervening
open segment is nontrivial and lies in $\mathcal N_\varepsilon^+$,
so this would be the first positive return of the route to its
bifurcation point.  Lemma~\ref{lem:no-retracing-return} says that the
incoming germ at that return is the opposite local half-germ, which
belongs to $\mathcal N_\varepsilon^-$ by Lemma~\ref{lem4.7}.  This
contradicts the positive orientation immediately before the return.
Hence every nontrivial point of the trapped positive half-route
belongs to $\mathcal N_\varepsilon^+$.

We may now use Lemma~\ref{lem5dw:fixed-eps-X}.  Its hypotheses hold
uniformly along the trapped route, with the orientation just proved,
and it gives, for this fixed $\varepsilon$, a uniform bound in the
original Banach space $X$.  Combining this bound with the parameter
bound and the closed $\mathcal O_\delta^{\rm dw}$ margins of
$\overline U$, we can enclose the whole positive half-route in a
bounded closed subset $Q$ of $\mathcal O_{\delta'}$; the strict
inequality $\delta'<\delta$ preserves ellipticity and obliqueness on
the $X$-closure.  The escape alternative in
Theorem~\ref{thm4.6} is therefore impossible, so the distinguished
route would have to be periodic.

Periodicity makes the entire route equal to the image of the trapped
positive half-route.  Its $X$-closure is a bounded closed subset of
$\mathcal O_{\delta'}$, and, because the image remains in
$\overline U$, it contains no nonzero $q$-independent solution.
Lemma~\ref{lem:global-nodal-loop} now excludes the resulting periodic
closed loop.  This contradiction proves that the positive half-route
leaves $U$.  At its first exit all preceding nontrivial points lie in
$\mathcal N_\varepsilon^+$; the exclusion of $q$-independent states
and the sequential closedness in Lemma~\ref{lem4.8} place the exit
point in the same positive orientation.
\end{proof}

Let $\mathcal K_\delta^{\rm dw}$ be the maximal connected subset of
$\mathcal C_\delta^{\rm dw}$ containing $(\lambda_*^0,0,0)$.

\begin{theorem}[Exact global continuum]
\label{thm5.7}
For every sufficiently small $\delta>0$,
$\mathcal K_\delta^{\rm dw}$ is nontrivial and either is unbounded in
$\mathbb R\times X_{\rm dw}$ or
\[
 \overline{\mathcal K_\delta^{\rm dw}}^{\,\mathbb R\times X_{\rm dw}}
 \cap\partial\mathcal O_\delta^{\rm dw}\ne\varnothing.
\]
\end{theorem}
\begin{proof}
Let $U$ be any bounded open neighborhood of $(\lambda_*^0,0,0)$ with
$\overline U\subset\mathcal O_\delta^{\rm dw}$.  By
Lemma~\ref{lem5.6}, for a sequence $\varepsilon_j\downarrow0$ there
are oriented regularized solutions on $\partial U$.
Proposition~\ref{lem5.5-uniform} gives a subsequence converging in
$\mathbb R\times X_{\rm dw}$ to an exact solution on $\partial U$,
and Lemma~\ref{lem5.4} places this limit in
$\mathcal C_\delta^{\rm dw}$.  The boundary-intersection hypothesis of
Lemma~\ref{thm5.1} is therefore satisfied.
\end{proof}

Choose $\delta_n\downarrow0$ so that the bifurcation point belongs to
every $\mathcal O_{\delta_n}^{\rm dw}$.  These sets are nested, and
maximality of the connected components through the common bifurcation
state gives
\[
 \mathcal K_{\delta_n}^{\rm dw}\subset
 \mathcal K_{\delta_{n+1}}^{\rm dw}.
\]
Hence
\begin{equation}\label{eq5dw:Kglobal}
 \mathcal K_{\rm dw}:=\bigcup_{n\ge n_0}\mathcal K_{\delta_n}^{\rm dw}
\end{equation}
is connected.

\begin{corollary}[Strong connectedness of the exact continuum]
\label{cor:strong-exact-continuum}
The same set $\mathcal K_{\rm dw}$ is contained in
$\mathbb R\times X_{\rm s}$ and is connected in the
$\mathbb R\times X_{\rm s}$ topology.
\end{corollary}
\begin{proof}
Fix a point of $\mathcal K_{\rm dw}$.  It belongs to some
$\mathcal K_{\delta_n}^{\rm dw}$ and has a bounded neighborhood in
$\mathbb R\times X_{\rm dw}$ with slightly smaller positive
admissibility margins.  Lemma~\ref{lem5dw:exact-strong-recovery}
shows sequentially, hence continuously in these metric spaces, that
the identity map from the exact solution set with the
$X_{\rm dw}$ topology into $X_{\rm s}$ is continuous near that
point.  It is therefore continuous on $\mathcal K_{\rm dw}$.
Since the latter is connected in $X_{\rm dw}$, its image, which is
the same set, is connected in $X_{\rm s}$.
\end{proof}

\section{Proof of Theorem \ref{thm2.1}}
In this section, we mainly focus on the proof of Theorem \ref{thm2.1}.
Before that, we first establish a key lemma.

\begin{lemma}[Strong a priori bound under fixed admissibility margins]
\label{lem6.1}
Fix $\delta>0$ and $M>0$.  There is a constant
$C=C(\delta,M,g,p_0,\gamma,\alpha)$ with the following property.
Every exact oriented solution
$(\lambda,\overline w,\underline w)\in
\mathcal O_\delta^{\rm dw}$ satisfying
\begin{equation}\label{eq6:lemma61-hypothesis}
 |\lambda|+\|\overline w\|_{C^0(\overline R_1)}
 +\|(\overline w_p,\underline w_p)\|_{C^0(\overline R_1)
 \times C^0(\overline R_2)}\le M
\end{equation}
satisfies
\begin{equation}\label{eq6:lemma61-conclusion}
 |\lambda|+
 \|(\overline w,\underline w)\|_{X_{\rm s}}\le C.
\end{equation}
Consequently, every exact oriented family obeying the displayed
hypothesis uniformly is bounded in $\mathbb R\times X_{\rm s}$.
\end{lemma}
\begin{proof}
All constants below depend only on the quantities displayed in the
statement.  We write $w$ for the piecewise perturbation and put
\[
 A:=a(p;\lambda)^{-1}+w_p,
 \qquad v:=w_q.
\]
The parameter bound, the boundedness of $\Gamma$, the fixed
$\mathcal O_\delta^{\rm dw}$ margins, and
\eqref{eq6:lemma61-hypothesis} give constants
$0<\kappa_0<\kappa_1$ such that
\begin{equation}\label{eq6:lemma61-A-bounds}
 \kappa_0\le a(p;\lambda),\ a(p;\lambda)^{-1},\ A(q,p)
 \le\kappa_1
 \qquad\text{on }\overline D.
\end{equation}

\smallskip
\noindent\emph{Step 1: a uniform gradient bound.}
By orientation, $v$ has one sign on each open half-period and
vanishes on $q=0,\pi$.  Differentiating the bulk equation in $q$
gives \eqref{eq:vk-nonlinear} with $\varepsilon=0$.  Apply the weak
maximum principle on
$(0,\pi)\times(-N,0)$, then let $N\to\infty$.  The end condition gives
$v\to0$ at the artificial lower boundary.  A strict extremum cannot
occur at $p=p_0$: continuity of $v$ and $v_p$ would contradict the
two one-sided Hopf inequalities.  Hence $\|v\|_{C^0(\overline D)}$
is controlled by its trace at $p=0$.  On that boundary the exact
surface condition gives, with $A=\lambda^{-1/2}+w_p$ there,
\[
 v^2=(\lambda-2gw)A^2-1.
\]
The right-hand side is uniformly bounded by
\eqref{eq6:lemma61-hypothesis} and
\eqref{eq6:lemma61-A-bounds}.  Thus
\begin{equation}\label{eq6:lemma61-gradient}
 \|Dw\|_{C^0(\overline D)}\le C.
\end{equation}
In particular, $w$ has a uniform $C^0$ bound on every fixed core:
the upper bound is part of \eqref{eq6:lemma61-hypothesis}, the
interface trace is common, and the lower value on a finite interval
is recovered by integrating $w_p$.

For later use, the principal matrix of either non-divergence-form
bulk equation is
\[
 \mathcal A(Dw)=
 \begin{pmatrix}
  1+v^2&-Av\\
  -Av&A^2
 \end{pmatrix},
 \qquad \det\mathcal A=A^2.
\]
The bounds for $A$ and $v$ make this matrix uniformly elliptic.  At
$p=p_0$, the function $a$ is continuous, and the transmission
conditions give common traces of $w_p$ and $w_q$.  Thus the two
principal matrices have the same trace at the interface, while
\[
 [w]=[w_p]=0.
\]
This is the perfect-transmission structure, with a complementing
constant bounded away from zero.  On $p=0$ write
\[
 G(w,Dw):=1+(2gw-\lambda)(\lambda^{-1/2}+w_p)^2+w_q^2.
\]
If $B=\lambda^{-1/2}+w_p$, then, using $G=0$,
\[
 G_{w_p}=2(2gw-\lambda)B
 =-\frac{2(1+w_q^2)}{B}.
\]
Consequently the surface operator is uniformly oblique; all its
first derivatives range in a fixed compact set determined by
\eqref{eq6:lemma61-hypothesis} and $\delta$.

\smallskip
\noindent\emph{Step 2: the uniform starting
$C^{1,\beta}$ estimate, including the vorticity interface.}
This is the point at which one cannot simply combine a smooth
quasilinear boundary estimate with a linear transmission estimate.
Instead, we use the weak divergence form across the interface.  For
$h=H+w$, define, for $\xi_2>0$,
\[
 \mathbf A(p,\xi)
 :=\left(\frac{\xi_1}{\xi_2},
 \Gamma(p)-\frac{1+\xi_1^2}{2\xi_2^2}\right).
\]
The weak formulation following \eqref{eq3.2} is exactly
\begin{equation}\label{eq6:lemma61-divergence}
 \operatorname{div}_{(q,p)}\mathbf A(p,Dh)=0
 \qquad\text{in }\mathbb S\times(-\infty,0).
\end{equation}
Although $\gamma$ jumps at $p=p_0$, the primitive $\Gamma$ is
globally Lipschitz.  Moreover,
\[
 D_\xi\mathbf A(p,\xi)=
 \begin{pmatrix}
  \xi_2^{-1}&-\xi_1\xi_2^{-2}\\
  -\xi_1\xi_2^{-2}&(1+\xi_1^2)\xi_2^{-3}
 \end{pmatrix}.
\]
By \eqref{eq6:lemma61-A-bounds} and
\eqref{eq6:lemma61-gradient}, this matrix satisfies
\begin{equation}\label{eq6:lemma61-div-ellipticity}
 \nu|\zeta|^2\le
 D_\xi\mathbf A(p,Dh)\zeta\cdot\zeta
 \le\nu^{-1}|\zeta|^2
\end{equation}
with one $\nu>0$ for the whole family.

We give the complete difference-quotient argument because it is what
removes the apparent gap at the transmission streamline.  Let
$B_{2\rho}$ be a disk or flat-interface disk whose closure does not
meet $p=0$, and let $\zeta\in C_c^1(B_{2\rho})$ equal one on
$B_\rho$.  The horizontal quotient $z_\tau=\delta_q^\tau h$ satisfies
\[
 \operatorname{div}(A_\tau Dz_\tau)=0,
 \qquad
 A_\tau=\int_0^1D_\xi\mathbf A
 (p,Dh+t\tau D\delta_q^\tau h)\,\dd t,
\]
and testing with $\zeta^2z_\tau$ gives
\begin{equation}\label{eq6:lemma61-horizontal-energy}
 \int_{B_\rho}|Dz_\tau|^2
 \le C\rho^{-2}\int_{B_{2\rho}}|z_\tau|^2
 \le C|B_{2\rho}|.
\end{equation}
Here the second inequality follows from the uniform Lipschitz bound
for $h$, and every constant is independent of $\tau$, including when
the disk crosses $p=p_0$: the weak equation is global and no
interface boundary term is created.

For the vertical quotient $y_\tau=\delta_p^\tau h$ we have
\begin{equation}\label{eq6:lemma61-p-difference}
 \operatorname{div}(\widetilde A_\tau Dy_\tau)
 =-\operatorname{div}
 \left(0,\delta_p^\tau\Gamma\right).
\end{equation}
Both averaged matrices obey \eqref{eq6:lemma61-div-ellipticity}, and
$\|\delta_p^\tau\Gamma\|_{L^\infty}\le\|\gamma\|_{L^\infty}$,
even when the quotient crosses $p=p_0$.  Testing
\eqref{eq6:lemma61-p-difference} with $\zeta^2y_\tau$ and using
Young's inequality gives
\begin{equation}\label{eq6:lemma61-vertical-energy}
 \int_{B_\rho}|Dy_\tau|^2
 \le C\rho^{-2}\int_{B_{2\rho}}|y_\tau|^2
   +C\int_{B_{2\rho}}|\delta_p^\tau\Gamma|^2
 \le C|B_{2\rho}|.
\end{equation}
Thus the two families of difference quotients have uniform local
$W^{1,2}$ bounds, not merely uniform supremum bounds.

We make the regularity input precise.  If
$\operatorname{div}(A Du)=\operatorname{div}F$ in a planar ball,
$A$ is merely measurable and satisfies
\eqref{eq6:lemma61-div-ellipticity}, and $F\in L^q$ for some
$q>2$, then the De Giorgi--Campanato estimate
\cite[Theorem 8.24]{GilbargT} gives numbers
$\beta=\beta(\nu,q)\in(0,1)$ and $C=C(\nu,q)$ such that
\begin{equation}\label{eq6:lemma61-explicit-campanato}
 R^\beta[u]_{C^{0,\beta}(B_{R/2})}
 \le C\left(
 R^{-1}\|u\|_{L^2(B_R)}
 +R^{1-2/q}\|F\|_{L^q(B_R)}
 \right).
\end{equation}
Apply \eqref{eq6:lemma61-explicit-campanato} to the quotient equations.
For $z_\tau$ the datum is zero; for $y_\tau$ take $q=4$ and
$F=(0,-\delta_p^\tau\Gamma)$, whose $L^4$ norm is uniformly
bounded because $\Gamma$ is globally Lipschitz.  The energy bounds
\eqref{eq6:lemma61-horizontal-energy}--
\eqref{eq6:lemma61-vertical-energy}, together with the uniform
$L^2$ bounds for the quotients, therefore supply
$\beta\in(0,1)$ and
\begin{equation}\label{eq6:lemma61-quotient-campanato}
 [z_\tau]_{C^{0,\beta}(B_{\rho/2})}
 +[y_\tau]_{C^{0,\beta}(B_{\rho/2})}\le C.
\end{equation}
The exponent and constant depend only on the ellipticity ratio and
the displayed bounds.  Because the quotient equations were obtained
from the single global weak equation
\eqref{eq6:lemma61-divergence}, this application is an interior one
even when the disk crosses $p=p_0$; in particular there is no
unverified interface boundary condition in
\eqref{eq6:lemma61-explicit-campanato}.  Taking weak $W^{1,2}$ and
uniform limits as $\tau\to0$ identifies the limits with $h_q$ and
$h_p$, respectively, and gives their $C^{0,\beta}$ bounds.  We
decrease $\beta$, if necessary, so that $0<\beta<\alpha$.

Near $p=0$ there is no vorticity interface.  We verify explicitly the
boundary regularity hypotheses.  On a fixed upper half-ball write
the height equation and its boundary operator as
\[
 \begin{aligned}
 \mathcal F(p,z,\xi,r)
 &=(1+\xi_q^2)r_{pp}-2\xi_q\xi_p r_{pq}
   +\xi_p^2r_{qq}+\gamma(-p)\xi_p^3,\\
 \mathcal G(z,\xi)&=1+\xi_q^2+2gz\xi_p^2.
 \end{aligned}
\]
The bounds \eqref{eq6:lemma61-A-bounds} and
\eqref{eq6:lemma61-gradient} confine
$(h,Dh)$ to one compact set of jets on which $\mathcal F$ is
uniformly elliptic and all of its first structural derivatives have
the natural bounded growth required in
\cite[Theorem 4.1]{LiebermanTrudinger}.  The upper-layer vorticity is
$C^{1,\alpha}$, the boundary is flat, and $h$ and $Dh$ are uniformly
bounded on the fixed upper half-ball.
Moreover, on $p=0$,
\[
 \frac{\partial}{\partial h_p}
 (1+h_q^2+2ghh_p^2)
 =4ghh_p=-\frac{2(1+h_q^2)}{h_p}\le-c_0<0.
\]
Thus $|\mathcal G_{\xi}\cdot e_p|\ge c_0$ and all first and second
derivatives of $\mathcal G$ are uniformly bounded on the same jet
set.  These are precisely the uniform obliqueness and boundary
structure constants in the cited theorem.  Its local boundary
gradient estimate therefore supplies the same $C^{0,\beta}$ bound
for $Dh$ on upper half-balls, with a constant independent of the
solution.  A finite cover handles a core containing both $p=0$ and
$p=p_0$.  The lower end is covered by translated enlarged unit
strips.  Their constants are independent of the translation; one
subtracts a constant from $h$ on each strip, and
\eqref{eq6:lemma61-gradient} controls the resulting oscillation.
Consequently
\begin{equation}\label{eq6:lemma61-C1beta}
 \|Dw\|_{C^{0,\beta}(\overline R_1)}
 +\|Dw\|_{C^{0,\beta}(\overline R_2)}\le C.
\end{equation}
For completeness, the global seminorm follows from the local ones:
points at distance less than one lie in a common enlarged patch,
whereas points at distance at least one are controlled by twice the
supremum norm in \eqref{eq6:lemma61-gradient}.

\smallskip
\noindent\emph{Step 3: a non-circular divergence/transmission
bootstrap.}
Put $v=h_q=w_q$ and
\[
 \mathbb A:=D_\xi\mathbf A(p,Dh)=
 \begin{pmatrix}
  h_p^{-1}&-h_qh_p^{-2}\\
  -h_qh_p^{-2}&(1+h_q^2)h_p^{-3}
 \end{pmatrix}.
\]
Taking a horizontal difference quotient in the single weak equation
\eqref{eq6:lemma61-divergence}, and then passing to the limit, gives
the global weak equation
\begin{equation}\label{eq6:lemma61-v-divergence}
 \operatorname{div}(\mathbb A Dv)=0
 \qquad\text{in }\mathbb S\times(-\infty,0).
\end{equation}
There is no source involving $\gamma$, because the primitive
$\Gamma$ is independent of $q$.  In particular, the interface is
not a boundary for \eqref{eq6:lemma61-v-divergence}.

The top condition also linearizes without introducing an unknown
coefficient of insufficient regularity.  Differentiating
$1+h_q^2+2ghh_p^2=0$ tangentially gives
\[
 2h_qv_q+2gh_p^2v+4ghh_pv_p=0.
\]
Using the undifferentiated condition to eliminate $h$ transforms
this identity exactly into
\begin{equation}\label{eq6:lemma61-v-robin}
 (\mathbb A Dv)\mathbin{\cdot}e_p-gv=0
 \qquad\text{on }p=0.
\end{equation}
Thus the boundary operator is the conormal derivative of the same
matrix appearing in the bulk equation, followed by the constant
Robin term $-gv$.  This identity is the reason no circular
$C^{1,\beta}$ assumption on a mean-value oblique coefficient is
needed.

By \eqref{eq6:lemma61-C1beta}, $\mathbb A$ is globally
$C^{0,\beta}$, including across $p=p_0$, and it satisfies the fixed
ellipticity bounds \eqref{eq6:lemma61-div-ellipticity}.  The local
divergence-form boundary Schauder estimate for
\eqref{eq6:lemma61-v-divergence}--
\eqref{eq6:lemma61-v-robin}, together with the interior estimate,
therefore gives on nested cores $K\Subset K^*$
\begin{equation}\label{eq6:lemma61-v-C1beta-estimate}
 \|v\|_{C^{1,\beta}(K)}
 \le C\|v\|_{C^0(K^*)}.
\end{equation}
This is the $C^{1,\beta}$ conormal estimate in
\cite[Section 16]{Lady}; its hypotheses here are precisely
$\mathbb A\in C^{0,\beta}$, a flat boundary, uniform ellipticity,
and the fixed Robin coefficient.  Translated unit strips give the
same estimate at the lower end, with one constant because
$\|v\|_{C^0}$ is globally bounded.

It follows that $h_{qq}=v_q$ and $h_{pq}=v_p$ are globally
$C^{0,\beta}$.  Solving the height equation for the remaining
second derivative gives, separately in the two vorticity layers,
\begin{equation}\label{eq6:lemma61-hpp-algebraic}
 h_{pp}=
 \frac{2h_qh_ph_{pq}-h_p^2h_{qq}-\gamma(-p)h_p^3}
 {1+h_q^2}.
\end{equation}
Hence $h_{pp}$ is piecewise $C^{0,\beta}$, and
\begin{equation}\label{eq6:lemma61-C2beta}
 \|w\|_{C^{2,\beta}(\overline R_1)}
 +\|Dw\|_{C^{1,\beta}(\overline R_2)}\le C.
\end{equation}

We record explicitly why the next interface estimate is uniform.
After \eqref{eq6:lemma61-C2beta}, the matrix $\mathbb A$ is
piecewise $C^{1,\beta}$ and has a common trace $A_0$ at $p=p_0$.
The weak equation supplies
\[
 [v]=0,\qquad
 [(\mathbb A Dv)\mathbin{\cdot}e_p]=0.
\]
For a nonzero tangential frequency $\xi$, the two decaying roots for
the frozen matrix $A_0=(a^{ij})$ are
\[
 \tau_+=\frac{\sqrt{\det A_0}|\xi|+ia^{qp}\xi}{a^{pp}},
 \qquad
 \tau_-=\frac{\sqrt{\det A_0}|\xi|-ia^{qp}\xi}{a^{pp}}.
\]
On these modes the two frozen conormal traces are respectively
$-\sqrt{\det A_0}|\xi|$ and
$+\sqrt{\det A_0}|\xi|$.  The transmission symbol therefore has
determinant of modulus
\begin{equation}\label{eq6:lemma61-complementing-symbol}
 2\sqrt{\det A_0}|\xi|\ge c|\xi|.
\end{equation}
This proves the uniform complementing condition.  The piecewise
divergence-form transmission/conormal Schauder estimate
\cite{Agmon,Lady} now yields
\[
 \|v\|_{C^{2,\beta}\text{ piecewise}}\le C.
\]
All third derivatives of $h$ containing a $q$-derivative are thus
$C^{0,\beta}$.  Differentiating
\eqref{eq6:lemma61-hpp-algebraic} in $p$ expresses $h_{ppp}$ in
terms of those derivatives and $\gamma'$; it does not reproduce
$h_{ppp}$ on the right.  Consequently
\begin{equation}\label{eq6:lemma61-C3beta}
 \|w\|_{C^{3,\beta}(\overline R_1)}
 +\|Dw\|_{C^{2,\beta}(\overline R_2)}\le C.
\end{equation}
All estimates on the lower cylinder are global: adjacent points lie
in a common enlarged unit strip, and pairs at distance at least one
are controlled by the already obtained supremum norms.

\smallskip
\noindent\emph{Step 4: upgrade to the prescribed exponent
$\alpha$.}
The bound \eqref{eq6:lemma61-C3beta} makes $\mathbb A$ uniformly
piecewise $C^{1,\alpha}$ for every fixed $\alpha<1$: on each unit
patch its first derivatives are Lipschitz, with norm controlled by
the displayed $C^{2,\beta}$ bound for $Dh$.  The common trace,
the symbol bound \eqref{eq6:lemma61-complementing-symbol}, and the
conormal--Robin identity \eqref{eq6:lemma61-v-robin} are unchanged.
Applying the same transmission/conormal Schauder estimate with
exponent $\alpha$ gives
\begin{equation}\label{eq6:lemma61-v-C2alpha}
 \|w_q\|_{C^{2,\alpha}\text{ piecewise}}\le C.
\end{equation}
Thus every third derivative with a $q$-derivative is
$C^{0,\alpha}$.  Formula
\eqref{eq6:lemma61-hpp-algebraic} first gives piecewise
$C^{1,\alpha}$ control of $h_{pp}$; differentiating it in $p$ then
gives piecewise $C^{0,\alpha}$ control of $h_{ppp}$.  Here the only
new coefficient is $\gamma'$, which is piecewise
$C^{0,\alpha}$ by hypothesis.  Therefore
\[
 \|w_{ppp}\|_{C^{0,\alpha}\text{ piecewise}}\le C.
\]
Together with \eqref{eq6:lemma61-v-C2alpha}, this gives
\begin{equation}\label{eq6:lemma61-Xdw-bound}
 |\lambda|+
 \|(\overline w,\underline w)\|_{X_{\rm dw}}\le C.
\end{equation}

\smallskip
\noindent\emph{Step 5: recovery of the absolute lower zero mode.}
The solution is exact and oriented and retains the fixed
$\mathcal O_\delta^{\rm dw}$ margins.  We may therefore apply
Lemma~\ref{lem5dw:exact-strong-recovery} with the bound
\eqref{eq6:lemma61-Xdw-bound}.  It gives a uniform
$C^{0,\alpha}(\overline R_2)$ bound for the reconstructed lower
height, including its absolute zero Fourier mode.  Combining this
with the derivative bounds above and the definition
\eqref{eq5dw:Xstrong} proves \eqref{eq6:lemma61-conclusion}.
\end{proof}

\begin{lemma}[Uniform separation from the lower ellipticity edge]
\label{lem6:deep-tail-separation}
Let $I\Subset(-2\Gamma_{\inf},\infty)$ and $M>0$.  There exist
$P>0$ and $\kappa>0$ with the following property.  If an exact
oriented solution satisfies
\[
 \lambda\in I,\qquad
 \|\overline w\|_{C^0(\overline R_1)}\le M,
 \qquad
 0<h_p\le M,
\]
then
\begin{equation}\label{eq:uniform-deep-hp-lower}
 h_p(q,p)\ge\kappa
 \qquad(q\in\mathbb S,\ p<-P).
\end{equation}
\end{lemma}
\begin{proof}
For an exact solution, Lemma~\ref{lem5dw:exact-strong-recovery}
gives $\mathfrak J\equiv0$.  Hence
\begin{equation}\label{eq:section6-exact-flux}
 \frac1{2\pi}\int_{-\pi}^{\pi}
 \frac{1+h_q(q,p)^2}{h_p(q,p)^2}\,\dd q
 =a(p;\lambda)^2.
\end{equation}
Since $I$ is compact and $\Gamma$ is bounded, the right-hand side is
uniformly bounded.  Jensen's inequality gives
\[
 \langle h_p\rangle(p)\ge a(p;\lambda)^{-1}\ge c_I>0.
\]
Moreover, using $h_p\le M$ in \eqref{eq:section6-exact-flux},
\[
 \int_{-\pi}^{\pi}|h_q(q,p)|^2\,\dd q
 \le M^2\int_{-\pi}^{\pi}
 \frac{h_q(q,p)^2}{h_p(q,p)^2}\,\dd q
 \le C_I M^2.
\]
Hence the horizontal oscillation of $h(\cdot,p)$ is uniformly
bounded.  The upper trace is uniformly bounded by the hypothesis on
$\overline w$, while
\[
 \frac{\dd}{\dd p}\langle h\rangle=\langle h_p\rangle\ge c_I.
\]
Therefore
\[
 \sup_{q\in\mathbb S}h(q,p)\longrightarrow-\infty
 \qquad\text{uniformly as }p\to-\infty.
\]
Thus sufficiently deep streamlines lie a fixed Euclidean distance
below both the free surface and the transmission interface.

Suppose that \eqref{eq:uniform-deep-hp-lower} fails.  Then there are
solutions as above and points $(q_j,p_j)$ with $p_j\to-\infty$ and
\[
 h_{j,p}(q_j,p_j)\longrightarrow0.
\]
Let $(x_j,y_j)=(q_j,h_j(q_j,p_j))$ and set
\[
 r_j:=c_j-u_j=-\psi_{j,y}=h_{j,p}^{-1}>0.
\]
For all large $j$, a fixed ball $B_{2\rho}(x_j,y_j)$ is contained in
the lower smooth-vorticity layer.  In this ball,
\[
 \Delta r_j+\gamma'(\psi_j)r_j=0.
\]
Since $\gamma'$ is uniformly bounded, the interior Harnack inequality
gives
\[
 r_j(x,y)\ge C^{-1}r_j(x_j,y_j)
 \qquad\text{on }B_\rho(x_j,y_j),
\]
with $C$ independent of $j$.

Take the connected component of the streamline $p=p_j$ inside
$B_\rho(x_j,y_j)$ containing $(x_j,y_j)$.  It meets
$\partial B_\rho$, so its Euclidean arclength is at least $\rho$.
Along this streamline,
\[
 v_j^2=r_j^2h_{j,q}^2,\qquad
 \dd s=\sqrt{1+h_{j,q}^2}\,\dd q.
\]
Consequently,
\[
 \int_{-\pi}^{\pi}(r_j^2+v_j^2)
       (q,h_j(q,p_j))\,\dd q
 \ge
 \int_{\text{component}}r_j^2\,\dd s
 \ge C^{-2}\rho\,r_j(x_j,y_j)^2
 \longrightarrow\infty.
\]
But \eqref{eq:section6-exact-flux} gives
\[
 \int_{-\pi}^{\pi}(r_j^2+v_j^2)
       (q,h_j(q,p_j))\,\dd q
 =2\pi a(p_j;\lambda_j)^2,
\]
which is uniformly bounded for $\lambda_j\in I$.  This contradiction
proves the lemma.
\end{proof}

\begin{lemma}[Pressure bound across the vorticity interface]
\label{lem6:pressure}
Let $(\lambda,\overline w,\underline w)$ be an exact admissible
solution and let $(c,u,v,\psi)$ be the corresponding water-wave
solution.  Set
\[
 M_\gamma:=\max\left\{0,\sup_{s\ge0}\gamma(s)\right\}.
\]
Then
\begin{equation}\label{eq6:pressure-bound}
 \frac12\bigl((c-u)^2+v^2\bigr)+gy-\Gamma(-\psi)
 -\frac{M_\gamma}{2}\psi\le0
 \qquad\text{in }\Omega_\eta.
\end{equation}
\end{lemma}
\begin{proof}
In each smooth-vorticity layer, the calculation in
\cite[Lemma 5.2]{Hur11} applies to
\[
 \mathscr B=
 \frac12|\nabla\psi|^2+gy-\Gamma(-\psi)
 -\frac{M_\gamma}{2}\psi
\]
and gives the required elliptic differential inequality.  It remains
to check the internal streamline $p=p_0$.  Since $\psi=-p$, in height
variables
\[
 \mathscr B=
 \frac{1+h_q^2}{2h_p^2}+gh-\Gamma(p)
 +\frac{M_\gamma}{2}p.
\]
The transmission conditions imply continuity of $\mathscr B$ across
$p=p_0$.  Using the height equation in either layer,
\[
 \mathscr B_p=
 \frac{h_{qq}}{h_p}
 -\frac{h_qh_{pq}}{h_p^2}
 +gh_p+\frac{M_\gamma}{2}.
\]
The traces of $h_q,h_{qq},h_{pq}$ and $h_p$ agree from the two sides,
and therefore
\[
 [\mathscr B_p]_{p=p_0}=0.
\]
The tangential derivative $\mathscr B_q$ is continuous as well.
If $n$ is either fixed choice of unit normal to the streamline, then
the change of variables gives, up to the common positive normalizing
factor,
\[
 \partial_n\mathscr B
 =\frac{1+h_q^2}{h_p}\mathscr B_p-h_q\mathscr B_q.
\]
Hence the physical normal derivative has the same trace from both
sides.  A positive interface maximum would instead give opposite
strict one-sided normal derivatives by the two Hopf lemmas.  Thus an
interface maximum is impossible.

For completeness, apply the maximum principle first on the periodic
truncation of the physical fluid above $y=-R$, separately in the two
layers.  On the free surface $\mathscr B=0$ by Bernoulli's law.  The
deep-water conditions give $|\nabla\psi|\to c$ and
$\psi(x,y)=-cy+o(|y|)$ as $y\to-\infty$; since $g>0$ and
$M_\gamma\ge0$, this implies
\[
 \limsup_{R\to\infty}\sup_{y=-R}\mathscr B\le0.
\]
There are no lateral boundary terms because of periodicity.  Letting
$R\to\infty$ proves \eqref{eq6:pressure-bound}.  This is the argument
of \cite[Lemma 5.2]{Hur11}, supplemented by the transmission check
above.
\end{proof}

We next isolate the exhaustion argument which converts the abstract
Whyburn alternative into the concrete alternatives used in the main
theorem.  This also makes explicit where the natural deep-water
topology enters the final proof.

\begin{proposition}[Exhaustion of the global alternatives]
\label{prop:section6-exhaustion}
Let $\mathcal K=\mathcal K_{\rm dw}$ be the nested union in
\eqref{eq5dw:Kglobal}.  At least one of the following alternatives
occurs:
\begin{itemize}
  \item[(1)] there exists a sequence
$(\lambda_k,\overline w_k,\underline w_k)\in\mathcal K$ such that
either $\lambda_k\to\infty$ or
$\lambda_k+2\Gamma_{\inf}\to0$;
  \item[(2)] there exists a sequence in $\mathcal K$ such that
$\|\overline w_k\|_{C^0(\overline R_1)}\to\infty$;
  \item[(3)] there exists a sequence in $\mathcal K$ such that
$\|\underline w_k\|_{C^0(\overline R_2)}\to\infty$;
  \item[(4)] there exists a sequence in $\mathcal K$ such that
$\max_{\overline R_1}\partial_p\overline w_k\to\infty$;
  \item[(5)] there exists a sequence in $\mathcal K$ such that
$\sup_{\overline R_2}\partial_p\underline w_k\to\infty$;
  \item[(6)] there are $z_k=(\lambda_k,\overline w_k,
\underline w_k)\in\mathcal K$, numbers $\eta_k\downarrow0$, and
points in $\overline R_1$ at which
$a^{-1}(p;\lambda_k)+\overline w_{k,p}=\eta_k$;
  \item[(7)] there are $z_k\in\mathcal K$, numbers
$\eta_k\downarrow0$, and points in $\overline R_2$ at which
$a^{-1}(p;\lambda_k)+\underline w_{k,p}=\eta_k$;
  \item[(8)] there are $z_k\in\mathcal K$, numbers
$\eta_k\downarrow0$, and points on $p=0$ at which
$2\lambda_k-4g\overline w_k=\eta_k$.
\end{itemize}
\end{proposition}
\begin{proof}
Suppose that none of \textup{(1)}--\textup{(8)} occurs.  The
negation of \textup{(1)} gives a compact interval
$I\Subset(-2\Gamma_{\inf},\infty)$ containing every parameter in
$\mathcal K$.  The negations of \textup{(2)} and \textup{(3)} give
uniform upper- and lower-layer $C^0$ bounds.  The negations of
\textup{(4)} and \textup{(5)} give a uniform upper bound for both
components of $w_p$.  Since $h_p=a^{-1}+w_p>0$ and $a^{-1}$ is
uniformly bounded for $\lambda\in I$, these are in fact uniform
two-sided $C^0$ bounds for $w_p$.

If the infimum of $h_p$ over $\mathcal K\times\overline R_i$ were
zero for either layer, points at which it is less than $1/k$ would,
after taking their positive values as $\eta_k$, produce
alternative \textup{(6)} or \textup{(7)}.  Hence their negations
give one common lower bound for $h_p$ in both layers.  The same
argument applied to the positive surface quantity
$2\lambda-4g\overline w(\cdot,0)$ shows, by the negation of
\textup{(8)}, that it has a common positive lower bound.  We may
therefore choose $\eta>0$ so small that
\begin{equation}\label{eq:section6-global-fixed-margin}
 \mathcal K\subset\mathcal O_{4\eta}^{\rm dw}.
\end{equation}

The upper $C^0$ bound and the two-sided $C^0$ bound for $w_p$ verify
the hypotheses of Lemma~\ref{lem6.1}, uniformly on $\mathcal K$.
Together with \eqref{eq:section6-global-fixed-margin}, that lemma
gives
\begin{equation}\label{eq:section6-global-Xdw-bound}
 \sup_{(\lambda,\overline w,\underline w)\in\mathcal K}
 \bigl(|\lambda|+
 \|(\overline w,\underline w)\|_{X_{\rm s}}\bigr)<\infty.
\end{equation}
In particular, $\mathcal K$ is bounded in
$\mathbb R\times X_{\rm dw}$.

Choose $n$ so large that $\delta_n<\eta$.  Then
$\mathcal K_{\delta_n}^{\rm dw}\subset\mathcal K$ is bounded.
Moreover, convergence in $X_{\rm dw}$ preserves the parameter,
derivative, and surface traces occurring in
\eqref{eq:section6-global-fixed-margin}; consequently
\[
 \overline{\mathcal K_{\delta_n}^{\rm dw}}
 ^{\,\mathbb R\times X_{\rm dw}}
 \subset\mathcal O_{\delta_n}^{\rm dw}.
\]
Thus this component is neither unbounded nor has closure meeting
$\partial\mathcal O_{\delta_n}^{\rm dw}$, contradicting
Theorem~\ref{thm5.7}.  The contradiction proves that the list is
exhaustive.
\end{proof}

Now we are in position to give the proof of our main Theorem
\ref{thm2.1}. Let us set
$\mathcal{K}:=\mathcal K_{\rm dw}$ as in \eqref{eq5dw:Kglobal}.
By Corollary~\ref{cor:strong-exact-continuum}, this is a connected
set in the strong physical topology.  Proposition
\ref{prop:section6-exhaustion} gives the exhaustive alternatives
\textup{(1)}--\textup{(8)}.  We now show that each of them gives one
of the two conclusions in Theorem~\ref{thm2.1}.

First,
$c^2=\lambda+2\Gamma_\infty$.  If the first part of alternative
\textup{(1)} occurs, then $c_k\to\infty$.  Suppose instead that
\[
 \lambda_k+2\Gamma_{\inf}\longrightarrow0.
\]
Put
\[
 \eta_k:=\lambda_k+2\Gamma_{\inf}\downarrow0
\]
and choose $p_k\le0$ such that
\[
 \Gamma(p_k)<\Gamma_{\inf}+\eta_k.
\]
Then
\[
 a(p_k;\lambda_k)^2
 =\lambda_k+2\Gamma(p_k)
 <3\eta_k\longrightarrow0.
\]
By \eqref{eq:section6-exact-flux} and $1+h_{k,q}^2\ge1$,
\[
 \frac{1}{\bigl(\sup_q h_{k,p}(q,p_k)\bigr)^2}
 \le
 \frac1{2\pi}\int_{-\pi}^{\pi}
 \frac{1}{h_{k,p}(q,p_k)^2}\,\dd q
 \le a(p_k;\lambda_k)^2\longrightarrow0.
\]
Consequently
\[
 \sup_q h_{k,p}(q,p_k)\longrightarrow\infty,
\]
which is alternative~\textup{(b)}.  Thus the lower parameter edge
$\lambda=-2\Gamma_{\inf}$ introduces no additional outcome, and no
assumption $\Gamma_{\inf}=\Gamma_\infty$ is needed.

For the rest of the proof suppose that \textup{(a)} and
\textup{(b)} both fail.  After passage to subsequences, $\lambda$
therefore ranges in a compact interval
$I\Subset(-2\Gamma_{\inf},\infty)$ and
\begin{equation}\label{eq:section6-hp-upper}
 0<h_p\le M\qquad\text{on }\overline D.
\end{equation}
In particular $a^{-1}$ is uniformly bounded and
$w_p=h_p-a^{-1}$ is bounded from below.  Alternatives
\textup{(4)} and \textup{(5)} are impossible, since an unbounded
positive maximum of $w_p$ would make $h_p$ unbounded.  Alternative
\textup{(8)} is also impossible: at the points in that alternative,
the surface condition gives
\[
 h_p^2=\frac{1+w_q^2}{\lambda-2gw}
 =\frac{2(1+w_q^2)}{\delta_k}\longrightarrow\infty.
\]

We next establish the surface and upper-layer bounds needed for
alternatives \textup{(2)} and \textup{(3)}, without choosing a bottom
level separately for each solution.  The
exact flux identity \eqref{eq:exact-flux-zero} gives $m(0)=0$.
On the surface the nodal pattern makes the maximum and minimum occur
at $q=0$ and $q=\pi$, while
\begin{equation}\label{eq:surface-section6}
 h_p^2=\frac{1+w_q^2}{\lambda-2gw}.
\end{equation}
If the surface $C^0$ norm were unbounded, its maximum would remain
bounded above by $\lambda/(2g)$, so its minimum would tend to
$-\infty$.  Since the surface mean is zero, the one-dimensional mean
value theorem gives
\[
 \|w_q(\cdot,0)\|_{C^0}\ge c\|w(\cdot,0)\|_{C^0}.
\]
Evaluating \eqref{eq:surface-section6} at a point of maximal
$|w_q|$, and using
$\lambda-2gw\le C(1+\|w(\cdot,0)\|_{C^0})$, gives
$\sup h_p\to\infty$, a contradiction.  Hence the
surface trace is bounded.  Because the upper layer has finite height,
\eqref{eq:section6-hp-upper} and the boundedness of $a^{-1}$ show
that an unbounded $\|\overline w\|_{C^0}$ would force an unbounded
positive maximum of $\overline w_p$, which has already been excluded.
Thus alternative~\textup{(2)} cannot occur.

We next exclude \textup{(6)} and \textup{(7)}.  At a point
where $h_{k,p}=\delta_k\to0$, the physical horizontal velocity
satisfies
\begin{equation}\label{eq:section6-fast-velocity}
 u_k=c_k-h_{k,p}^{-1}\longrightarrow-\infty.
\end{equation}
The surface and upper-layer bounds just proved, together with
\eqref{eq:section6-hp-upper}, make every physical image of a fixed
$p$-core uniformly bounded.  If the points in alternatives
\textup{(6)} or \textup{(7)} remain in such a core, the pressure bound \eqref{eq6:pressure-bound} contradicts
\eqref{eq:section6-fast-velocity}: the kinetic term diverges while
$gy$, $\Gamma(-\psi)$, and $M_\gamma\psi$ are uniformly bounded.
If the points in alternative~\textup{(7)} satisfy $p_k\to-\infty$,
Lemma~\ref{lem6:deep-tail-separation}, applied with
\eqref{eq:section6-hp-upper} and the already established upper-layer
bound, gives $h_{k,p}(q_k,p_k)\ge\kappa>0$, again a contradiction.
Thus alternatives \textup{(6)} and \textup{(7)} are impossible
when \textup{(a)} and \textup{(b)} fail.

The preceding argument is uniform: otherwise a sequence at which
$h_p\to0$ would fall either in a fixed $p$-core or in the deep tail,
and one of the two contradictions just obtained would apply.
Together with the already excluded surface-obliqueness loss and the
compact parameter interval, it therefore supplies one fixed
admissibility margin $\delta>0$ for the family under consideration.
Lemma~\ref{lem6.1} now applies.  Its first four steps give the
$X_{\rm dw}$ derivative bound without assuming an absolute lower
zero-mode bound, and its final step invokes
Lemma~\ref{lem5dw:exact-strong-recovery} to recover precisely that
missing bound.  Consequently the family is bounded in $X_{\rm s}$;
in particular,
$\|\underline w\|_{C^0(\overline R_2)}$ cannot diverge, so
alternative~\textup{(3)} is impossible.  This replaces the invalid
selection of one common truncation level from merely pointwise
deep-water convergence.

All eight alternatives have now been exhausted.  Hence the connected
set $\mathcal K$ contains a sequence for which either $c_k\to\infty$
or $\sup_{\overline D}(h_k)_p\to\infty$, completing the proof of
Theorem~\ref{thm2.1}.

\bibliographystyle{amsplain}
\makeatletter
\def\@biblabel#1{#1.~}
\makeatother


\begin{thebibliography}{10}

\bibitem{AT81}
C. Amick and J. Toland,
\newblock {On solitary water-waves of finite amplitude},
\newblock {\em Arch. Rational Mech. Anal.} 76 (1981), 9-95.

\bibitem{Agmon}
S. Agmon, A. Douglis, L. Nirenberg,
\newblock {Estimates near the boundary for
solutions of elliptic partial differential equations satisfying
general boundary conditions I},
\newblock {\em Comm.Pure Appl. Math.} 12 (1959),
623-727.


\bibitem{BuffoniT}
B. Buffoni and J. Toland,
\newblock {Analytic Theory of Global Bifurcation: An
Introduction, Princeton University Press}, 2003.

\bibitem{CaoQ}
D. Cao, G. Qin, C. Zou,
\newblock {Existence of steady gravity-capillary water waves with concentrated vortex sheet},
\newblock {\em Commun. Inf. Syst.} 25 (2025), 205-231.


\bibitem{Constantin01}
A. Constantin,
\newblock {On the deep water wave motion},
\newblock {\em J. Phys. A}, 34 (2001),
1405-1417.


\bibitem{Constantin}
A. Constantin,
\newblock {Nonlinear water waves with applications to
wave-current interactions and tsunamis, CBMS-NSF Conference Series
in Applied Mathematics 81, SIAM, Philadelphia}, 2011.


\bibitem{ConstantinS}
A. Constantin, W. Strauss,
\newblock {Exact steady periodic water waves with vorticity},
\newblock {\em Comm. Pure Appl. Math.}, 57 (2004), 481-527.


\bibitem{ConstantinS1}
A. Constantin, W. Strauss,
\newblock {Periodic traveling gravity water waves
with discontinuous vorticity},
\newblock {\em Arch. Ration. Mech. Anal.}, 202 (2011), 133-175.

\bibitem{ConstantinSV}
A. Constantin, W. Strauss, E.
V$\breve{a}$rv$\breve{a}$ruc$\breve{a}$,
\newblock {Global bifurcation of
steady gravity water waves with critical layers},
\newblock {\em Acta Math.}, 217 (2016), 195-262.


\bibitem{ConstantinV}
A. Constantin, E. V$\breve{a}$rv$\breve{a}$ruc$\breve{a}$,
\newblock {Steady
periodic water waves with constant vorticity: regularity and local
bifurcation},
\newblock {\em Arch. Ration. Mech. Anal.}, 199 (2011), 33-67.


\bibitem{Cordoba}
A. C\'{o}rdoba, D.A. C\'{o}rdoba, F. Gancedo,
\newblock {Interface evolution: the Hele--Shaw and Muskat problems},
\newblock {\em Ann. Math.}, 173(1) (2011), 477-542.

\bibitem{Cheng}
C.H. Cheng, R. Granero-Belinch\'{o}n, S. Shkoller,
\newblock {Well-posedness of
the Muskat problem with $H^2$ initial data},
\newblock {\em Adv. Math.}, 286 (2016), 32-104.

\bibitem{ChuDE}
J. Chu, Q. Ding, J. Escher,
\newblock {Variational formulation of rotational steady water waves in two-layer flows},
\newblock {\em J. Math. Fluid Mech.}, 23 (2021), 23pp.


\bibitem{ChuE}
J. Chu, J. Escher,
\newblock {Steady periodic equatorial water waves with vorticity},
\newblock {\em Discrete Contin. Dyn. Syst.}, 39 (2019), 4713-4729.

\bibitem{ChuW}
J. Chu, L. Wang,
\newblock { Analyticity of rotational traveling gravity two-layer waves},
\newblock {\em  Stud. Appl. Math.}, 146 (2021), 605-634.

\bibitem{ChuY}
J. Chu, Y. Yang,
\newblock {Constant vorticity water flows in the equatorial $\beta$-plane approximation with centripetal forces},
\newblock {\em J. Differential Equations}, 269 (2020), 9336-9347.

\bibitem{DaiFZ}
G. Dai, T. Feng, Y. Zhang,
\newblock {The existence and geometric structure of periodic solutions to rotational electrohydrodynamic waves problem},
\newblock {\em J. Geom. Anal.}, 35 (2025), 23pp.

\bibitem{DaiLZ}
G. Dai, F. Li, Y. Zhang,
\newblock { Bifurcation structure and stability of steady gravity water waves with constant vorticity},
\newblock {\em J. Differential Equations}, 332 (2022), 306-332.

\bibitem{DaiXZ}
G. Dai, F. Xu, Y. Zhang,
\newblock { The dynamics of periodic traveling interfacial electrohydrodynamic waves: bifurcation and secondary bifurcation},
\newblock {\em J. Nonlinear Sci.}, 34 (2024), 31pp.


\bibitem{Dancer}
E. N. Dancer,
\newblock {Bifurcation theory for analytic operators},
\newblock {\em Proc. London
Math. Soc.}, 26 (1973), 359-384.

\bibitem{Dubreil}
M. L. Dubreil-Jacotin,
\newblock {Sur la d\'etermination rigoureuse des ondes
permanentes p\'erodiques d'ampleur finie},
\newblock {\em J. Math. Pures Appl.}, 13
(1934), 217-291.

\bibitem{Escher}
J. Escher, P. Laurenot, B.v. Matioc,
\newblock {Existence and stability of weak
solutions for a degenerate parabolic system modelling two-phase
flows in porous media},
\newblock {\em Ann. Inst. H. Poincar\'{e} C Anal. Non
Lin\'{e}aire}, 28 (2011), 583-598.

\bibitem{Gilbarg}
D. Gilbarg,
\newblock {The Phragm\'en-Lindel\"{o}f theorem for elliptic partial
differential equations},
\newblock {\em J. Rational Mech. Anal.}, 1 (1952), 411-417.

\bibitem{GilbargT}
D. Gilbarg and N. Trudinger,
\newblock {Elliptic Partial Differential Equations
of Second Order, Springer, Berlin}, 2001.

\bibitem{Haziot}
S.V. Haziot,
\newblock {Stratified large-amplitude steady periodic water waves with critical layers},
\newblock {\em Comm. Math. Phys.}, 381 (2021), 765-797.

\bibitem{HenryBM}
D. Henry, B.V. Matioc,
\newblock {On the existence of steady periodic
capillary-gravity stratified water waves},
\newblock {\em Ann. Sc. Norm. Super. Pisa
CI. Sci.}, 12 (2013), 955-974.

\bibitem{HealeySimpson}
T. Healey and H. Simpson,
\newblock {Global continuation in nonlinear elasticity},
\newblock {\em Arch. Rational Mech. Anal.}, 143 (1998), 1--28.

\bibitem{Hur06}
V. M. Hur,
\newblock {Global bifurcation of deep-water waves with vorticity},
\newblock {\em SIAM J. Math. Anal.}, 37 (2006), 1482-1521.

\bibitem{Hur11}
V. M. Hur,
\newblock {Stokes waves with vorticity},
\newblock {\em J. Anal. Math.}, 113 (2011),
331-386.

\bibitem{Kato}
T. Kato,
\newblock {Perturbation Theory for Linear Operators, Springer-Verlag,
New York}, 1967.

\bibitem{Kielhofer}
H. Kielh\"{o}fer,
\newblock {Bifurcation theory. An introduction with
applications to partial differential equations. Second edition.
Applied Mathematical Sciences, Springer, New York}, 2012.

\bibitem{Krylov}
N. V. Krylov,
\newblock {Lectures on Elliptic and Parabolic Equations in
H\"{o}lder Spaces, Amer. Math. Soc., Providence, RI}, 1996.

\bibitem{Lady}
O.A. Ladyzhenskaya, N.N. Ural'tseva,
\newblock {Linear and quasilinear elliptic
equations. Translated from the Russian by Scripta Technica, Inc.
Translation editor: Leon Ehrenpreis. Academic Press, New
York-London}, 1968. xviii+495 pp.

\bibitem{LiebermanTrudinger}
G. M. Lieberman and N. S. Trudinger,
\newblock {Nonlinear oblique boundary value problems for nonlinear
elliptic equations},
\newblock {\em Trans. Amer. Math. Soc.}, 295 (1986), 509--546.


\bibitem{Matioc}
A. V. Matioc,
\newblock {Steady internal water waves with a critical layer
bounded by the wave surface},
\newblock {\em J. Nonlinear Math. Phys.}, 19 (2012), 21
pp.

\bibitem{MatiocM}
A. V. Matioc, B. V. Matioc,
\newblock {A new reformulation of the Muskat problem
with surface tension},
\newblock {\em J. Differential Equations}, 350 (2023), 308-335.

\bibitem{MartinM}
C. I. Martin, B. V. Matioc,
\newblock {Existence of capillary-gravity water
waves with piecewise constant vorticity},
\newblock {\em J. Differential Equations}, 256 (2014), 3086-3114.

\bibitem{Rabinowitz}
P. H. Rabinowitz,
\newblock {Some global results for nonlinear eigenvalue
problems},
\newblock {\em J. Funct. Anal.}, 7 (1971), 487-513.

\bibitem{Sinambela}
D. Sinambela,
\newblock {Large-amplitude solitary waves in two-layer density
stratified water},
\newblock {\em SIAM J. Math. Anal.}, 53 (2021), 4812-4864.

\bibitem{Varholm}
K. Varholm,
\newblock { Global bifurcation of waves with multiple critical layers},
\newblock {\em SIAM J. Math. Anal.}, 52 (2020), 5066-5089.

\bibitem{Var}
E. V$\breve{a}$rv$\breve{a}$ruc$\breve{a}$,
\newblock {On some properties of travelling water waves with vorticity},
\newblock {\em SIAM J. Math. Anal.}, 39 (5) (2008) 1686-1692.

\bibitem{Wahlen1}
E. Wahl\'{e}n,
\newblock {Steady periodic capillary-gravity waves with
vorticity},
\newblock {\em  SIAM J. Math. Anal.}, 38 (2006), 921-943.

\bibitem{Wahlen2}
E. Wahl\'{e}n,
\newblock {Steady water waves with a critical layer},
\newblock {\em J. Differential Equations}, 246 (2009), 2468-2483.

\bibitem{Erik}
E. Wahl\'{e}n, J. Weber,
\newblock {Global bifurcation of capillary-gravity
water waves with overhanging profiles and arbitrary vorticity},
\newblock {\em  Int.
Math. Res. Not.}, (2023), 17377-17410.

\bibitem{Erik1}
E. Wahl\'{e}n, J. Weber,
\newblock {Large-amplitude steady gravity water waves
with general vorticity and critical layers},
\newblock {\em Duke Math. J.}, 173
(2024), 2197-2258.

\bibitem{Walsh09}
S. Walsh,
\newblock {Stratified steady periodic water waves},
\newblock {\em SIAM J. Math.
Anal.}, 41 (2009), 1054-1105.

\bibitem{Walsh14a}
S. Walsh,
\newblock {Steady stratified periodic gravity waves with surface
tension I: local bifurcation},
\newblock {\em Discrete Contin. Dyn. Syst. Ser. A.}, 34 (2014), 3241-3285.


\bibitem{WangXZ}
J. Wang, F. Xu, Y. Zhang,
\newblock {The existence of steady stratified two-mode water waves with stagnation points},
\newblock {\em J. Math. Fluid Mech.}, 27 (2025), 21 pp.




\end{thebibliography}

\providecommand{\bysame}{\leavevmode\hbox
to3em{\hrulefill}\thinspace}
\providecommand{\MR}{\relax\ifhmode\unskip\space\fi MR }
\providecommand{\MRhref}[2]{%
  \href{http://www.ams.org/mathscinet-getitem?mr=#1}{#2}
} \providecommand{\href}[2]{#2}

\end{document}